\documentclass[11pt]{article}
\usepackage{color,xcolor}
\usepackage{graphicx,mathrsfs}
\usepackage{amssymb,bbm}
\usepackage{amsthm}
\usepackage{amsmath}
\usepackage[margin=1in]{geometry}
\usepackage{caption,numbysec}
\usepackage{subcaption}
\usepackage{float}
\usepackage{mdframed,pdfsync,mathrsfs} 
\usepackage{esint}
\usepackage{dsfont}
\usepackage[normalem]{ulem}

\usepackage[colorlinks=true,pdfpagemode=UseNone,urlcolor=blue,linkcolor=blue,citecolor=blue,breaklinks=true]{hyperref}

\usepackage{cleveref}

\newtheorem{thm}{Theorem}[section]
\newtheorem{lemma}[theorem]{Lemma}
\newtheorem{proposition}[theorem]{Proposition}
\newtheorem{prop}[theorem]{Proposition}

\newtheorem{cor}[theorem]{Corollary}

\newcommand{\bea}{\begin{eqnarray*}}
\newcommand{\eea}{\end{eqnarray*}}
\newcommand{\ben}{\begin{eqnarray}}
\newcommand{\een}{\end{eqnarray}}
\newcommand{\beq}{\begin{equation}}
\newcommand{\eeq}{\end{equation}}

\newcommand{\cW}{\mathcal W}

\newcommand{\C}{\ensuremath{\mathbb{C}}}

\newcommand{\R}{\ensuremath{\mathbb{R}}}

\newcommand{\Rm}{{\mathbb R}}

\newcommand{\1}{\mathds{1}}

\newcommand{\cD}{\mathcal{D}}

\newcommand{\cL}{\mathcal{L}}

\newcommand{\bal}{\begin{aligned}}
\newcommand{\enbal}{\end{aligned}}
\newcommand{\be}{\begin{equation}}
\newcommand{\ee}{\end{equation}}

\renewcommand{\hat}[1]{\widehat{#1}}

\newcommand{\eps}{\varepsilon}

\newcommand{\farc}{\frac}

\newcommand{\pdr}[2]{\frac{\partial{#1}}{\partial{#2}}}

\renewcommand{\d}{\partial}

\newcommand{\one}{\mathbbm{1}}

\newcommand{\under}{\underline}

\newcommand{\vphi}{\varphi}
\newcommand{\doverline}[1]{\overline{\overline{#1}}}
\newcommand{\dy}{\, dy}
\newcommand{\dx}{\, dx}
\newcommand{\disp}{\displaystyle}

\begin{document}
\title{Pushed, pulled and pushmi-pullyu fronts of the Burgers-FKPP equation}
\author{Jing An\footnote{Max-Planck Institute for Mathematics in the Sciences, Inselstr. 22, 04103 Leipzig, Germany;
jing.an@mis.mpg.de}  
\and Christopher Henderson\footnote{Department of Mathematics, University of Arizona, 
Tucson, AZ 85721, USA;
ckhenderson@math.arizona.edu} \and 
Lenya Ryzhik\footnote{Department of Mathematics, Stanford University, Stanford, CA 94305, USA;
ryzhik@stanford.edu}}
\maketitle
\numberbysection

\begin{abstract}
We consider the long time behavior of the solutions to the Burgers-FKPP equation with advection of a strength $\beta\in\Rm$.
This equation exhibits a transition from pulled to pushed front behavior at $\beta_c=2$.  
We prove convergence of the solutions to a traveling wave in a reference frame centered at a position $m_\beta(t)$ and study the asymptotics of the front location $m_\beta(t)$.  When $\beta < 2$, it has the same form as for the standard Fisher-KPP equation established by Bramson~\cite{Bramson1,Bramson2}: $m_\beta(t) = 2t - (3/2)\log(t) + x_\infty + o(1)$ as $t\to+\infty$.  This form is typical of pulled fronts.  When $\beta > 2$, the front is located at the position~$m_\beta(t)=c_*(\beta)t+x_\infty+o(1)$ with~$c_*(\beta)=\beta/2+2/\beta$, which is the typical form of pushed fronts.  However, at the critical value~$\beta_c = 2$, the expansion changes to $m_\beta(t) = 2t - (1/2)\log(t) + x_\infty + o(1)$, reflecting the ``pushmi-pullyu'' nature of the front.
The arguments for $\beta<2$ rely on a new weighted Hopf-Cole transform that allows to control the
advection term, when combined with additional steepness comparison arguments. The case~$\beta>2$ relies on standard pushed front techniques.  The proof in the case $\beta=\beta_c$ is much more
intricate and involves arguments not usually encountered in the study of the Bramson correction.
It relies on a somewhat hidden viscous conservation law
structure of the Burgers-FKPP equation at $\beta_c=2$ 
and utilizes a dissipation inequality, which comes from a relative entropy type computation, together with a weighted Nash inequality involving dynamically changing weights.
\end{abstract}

\section{Introduction}
We consider the long time behavior of the solutions to the Burgers-FKPP equation  
\begin{equation}\label{burgerskpp}
\begin{aligned}
&u_t + \beta uu_x = u_{xx} +u-u^2,~~t>0,~x\in\Rm.
\end{aligned}
\end{equation}
Here, $\beta \in \R$ is a parameter that measures the strength of the advection effect.
The relevance of this type of nonlinear advection-reaction-diffusion model in biological and chemical applications is discussed in 
Chapter 13.4 of~\cite{murray2007mathematical}.  We refer the interested reader to this book, as well as the references therein, rather than attempting to replicate the discussion here. 

Our main interest is in the study of the transition from the ``pulled'' to ``pushed'' nature of the Burgers-FKPP equation that happens
at $\beta_c=2$ and its effect on the long time behavior of the solutions.  In order to motivate the 
basic questions addressed in the present paper, and before introducing this transition phenomenon, we recall the
well-known results for the classical Fisher-KPP equation as well as the notions of ``pulled'' and ``pushed'' fronts.

\subsubsection*{The classical results for the Fisher-KPP equation}

When $\beta=0$, the Burgers-FKPP equation (\ref{burgerskpp}) reduces to the classical Fisher-KPP equation 
\begin{equation}\label{nov2502}
u_t=u_{xx}+u-u^2,
\end{equation}
that dates back to the seminal work of Fisher~\cite{Fisher} 
and Kolmogorov, Petrovskii and Piskunov~\cite{kolmogorov1937etude}, and has been studied extensively since. This equation arises in numerous applications   
in the physical and biological sciences, also discussed  in~\cite{murray2007mathematical}.
Beyond this, it has been used to study the fine properties of branching Brownian motion after McKean~\cite{McK} discovered its connection to~\eqref{nov2502}; 
see, for instance,~\cite{ABK1,ABK2,BD1,BD2,MRR} and references therein. 

Both of the original papers~\cite{Fisher} and~\cite{kolmogorov1937etude}   
showed that~(\ref{nov2502}) admits traveling wave solutions 
of the form~$u(t,x)=U_c(x-ct)$ for all $c\ge c_*=2$, which necessarily satisfy the ODE
\begin{equation}\label{nov3004}
-cU_c'=U_c''+U_c-U_c^2,
	~~~U_c(-\infty)=1,
	~~~\text{and}
	~~~U_c(+\infty)=0.
\end{equation}
We denote  by $U_*(x)$ the traveling wave moving with the minimal speed $c_*=2$.
Given $c \ge c_*$, solutions to~(\ref{nov3004}) are unique up to  translation in $x$.  
One may fix a normalization for~$U_*(x)$, for example, by requiring that $U_*(0)=1/2$. 

It was observed already in~\cite{Fisher,kolmogorov1937etude} (albeit at a different level of
mathematical rigor) that the
solution to (\ref{nov2502}) with an initial condition
that is a step function
\begin{equation}\label{dec302}
\begin{aligned}
 u(0,x)
 	= \1(x \leq 0)
	:= \begin{cases}
1, ~~~&\text{if } x\le 0,\\
\displaystyle 0,
&\text{if } x>0,
\end{cases}
\end{aligned}
\end{equation}
 converges in shape
to a traveling wave: there exists a reference frame $m(t)$ so that
\begin{equation}\label{nov2504}
	u(t,x+m(t))\to U_*(x),~~\hbox{ as $t\to+\infty$, uniformly in $x\in\Rm$.}
\end{equation}
It was also  argued informally in~\cite{Fisher} and proved in
\cite{kolmogorov1937etude} that      
\begin{equation}\label{nov2506}
m(t)=2t+o(t)\hbox{ as $t\to+\infty$.}
\end{equation}
We refer to $m(t)$ as the location of the front at time $t>0$ since, roughly, it separates the regions~$\{u(t,x) \approx 1\}$ for $x \ll m(t)$ and 
$\{u(t,x) \approx 0\}$ for $x \gg m(t)$.

This result was refined in the pioneering works by Bramson~\cite{Bramson1,Bramson2}, who 
used probabilistic techniques and the 
connection of~\eqref{nov2502} to branching Brownian motion to analyze the front position~$m(t)$.
In particular, Bramson showed 
that 
there exists a constant~$x_\infty$ that depends 
on the initial condition~$u_{\rm in}(x)$ for~(\ref{nov2502}), so that 
\begin{equation}\label{nov2508}
m(t)=2t-\farc{3}{2}\log t+x_\infty+o(1)\hbox{ as $t\to+\infty$.}
\end{equation}
This asymptotic expansion holds as long as the initial conditions $u_{\rm in}(x)$ decay sufficiently fast to zero as $x\to+\infty$. 
We also mention the related work by Uchiyama~\cite{Uchiyama} and Lau~\cite{Lau}, which are based on PDE techniques, and more recent refinements and alternative
proofs of Bramson's result in~\cite{berestycki2018new,Graham,HNRR,NRR1,NRR2,Roberts},
that use both probabilistic and PDE methods. Very recently, spectral techniques have been applied 
to study the ``Bramson shift'' in~\cite{AG,AS1,AS2}, including some problems that do not obey the comparison principle. 
 
In particular, it was shown in~\cite{Graham,NRR2} that
convergence in (\ref{nov2504}) is only algebraic in time:
\begin{equation}\label{nov2511}
|u(t,x+m(t))-U_*(x)|\le \farc{C(1+|x|)e^{-x}}{\sqrt{t}},~~~~x>0.
\end{equation}
The more precise results in~\cite{Graham} show that the convergence rate can not be improved to an exponential in time rate: even
convergence in shape to a traveling
wave does not hold beyond the order~$O(t^{-1})$. 
The algebraic rate of convergence of the solution to (\ref{nov2502}) to a shift of a traveling wave is closely related to the fact that FKPP fronts are ``pulled''
 -- that is, the long time behavior of the solutions is governed by the behavior far ahead of the front, where $u$ is small. 
Hence, the problem is not compact in a certain sense, making the algebraic-in-time (rather than exponential) rate of convergence 
 natural. 
 
The ``pulled'' behavior should be contrasted with the class of equations of the 
form similar to~(\ref{nov2502}):
\begin{equation}\label{nov2512}
v_t=v_{xx}+f(v)
\end{equation}
but with solutions that behave as ``pushed fronts" -- that is, the long time behavior of the solutions is governed by the behavior at the front, where $u$ is neither small nor approximately $1$.  An example of such $f$ is a bistable nonlinearity of the form $f(v)=v(1-v)(\theta-v)$,
with some~$\theta\in(0,1)$. Unlike the Fisher-KPP equation, in the pushed cases 
solutions to the initial value problem for~(\ref{nov2512}) with a rapidly decaying initial condition~$v(0,x)=v_{\rm{in}}(x)$
converge to a shift of the traveling wave~$U_f(x)$ for (\ref{nov2512}) exponentially fast in time: there exists $\omega>0$ such that
\begin{equation}\label{nov2510}
|v(t,x+m_f(t))-U_f(x)|\le Ce^{-\omega t}.
\end{equation}
Moreover,
the front location has the asymptotics
\be\label{aug1202}
m_f(t)=c_ft+x_\infty,~~\hbox{as $t\to+\infty$},
\ee
without any logarithmic in time correction. 
The rate of convergence in (\ref{nov2510})  is exponential in time precisely because the fronts are ``pushed,'' so that the long time behavior
is determined by what happens in a compact region around the front, and not by the tail behavior as $x\to+\infty$.
We refer to the monograph~\cite{Leach-Needham} and extensive presentations 
in~\cite{Ebert-vanSaarlos,vanSaarlos}
for, respectively, an applied mathematics and a physics perspective on pushed
and pulled fronts, and to~\cite{garnier2012inside} for a recent mathematical analysis, among
many other references.

\subsubsection*{The pulled to pushed fronts transition in the Burgers-FKPP equation}
  
As we have mentioned, the long time behavior of the reaction-diffusion equations that admit
either pulled or pushed fronts is reasonably well-understood, at least on an intuitive level. 
An interesting aspect of the 
Burgers-FKPP equation~(\ref{burgerskpp}) is that it exhibits a transition from the pulled 
to pushed behavior at $\beta_c=2$. 

The behavior of traveling waves for (\ref{burgerskpp}) already illustrates the change in behavior at $\beta_c = 2$.   For a given $\beta\in\Rm$, the Burgers-FKPP equation~(\ref{burgerskpp}) admits 
traveling wave solutions for all $c\ge c_*(\beta)$, with the minimal speed
\begin{equation}\label{speed-bis}
\begin{aligned}
    c_* = \begin{cases}
2, ~~~&\text{if } \beta \le 2,\\
\displaystyle\frac{\beta}{2}+\frac{2}{\beta},
&\text{if }\beta\geq 2.
\end{cases}
\end{aligned}
\end{equation}
The minimal speed traveling wave $\phi_\beta$ satisfies
\begin{equation}\label{nov2514}
	-c_*\phi_\beta'+\beta\phi_\beta\phi_\beta'
		=\phi_\beta''+\phi_\beta-\phi_\beta^2,
	~~~~\phi_\beta(-\infty)=1,
	~~~~\text{and}~~~~
	\phi_\beta(+\infty)=0.
\end{equation}
Once again, the solution to (\ref{nov2514}) is unique only up to a translation in $x$. We fix
the normalization by requiring that $\phi_\beta(0)=1/2$. It happens that the traveling wave profile for $\beta\ge 2$ is explicit.  Indeed, one can check by direct computation that
\be\label{jul1502}
\phi_{\beta}(x)=\farc{1}{1+e^{\beta x/2}},~~\hbox{ for $\beta\ge 2$.}
\ee
On the other hand, when $\beta < 2$, the profile of the minimal speed traveling wave is, to the best of our knowledge, not explicit, and the asymptotics of $\phi_\beta$ as $x\to+\infty$ are no longer purely exponential, being given by
\be\label{jul1602}
\phi_\beta(x)\sim (Ax+B)e^{-x},~~\hbox{ as $x\to+\infty$, for $\beta<2$,}
\ee 
with some $A>0$ and $B\in\Rm$ that depend on $\beta$.  
This was shown, for instance, in~\cite{murray2007mathematical} by a phase plane analysis. It is discussed further in Appendix~\ref{sec:phase-plane}.

As explained in Remark 1 of~\cite{garnier2012inside},
a quantitative mathematical criterion for a traveling wave profile $U_c(x)$, that moves with a
speed $c\ge 0$, to be pushed is that 
\be\label{aug1102}
U_c(x)e^{cx/2}\in L^2(\Rm).
\ee
Otherwise, a traveling wave is pulled. We will see the motivation behind the criterion
(\ref{aug1102}) in the discussion of the long time behavior of the
solutions to (\ref{burgerskpp}) for $\beta>2$, which is contained in \Cref{sec:beta>2}.  
According to this classification, the Burgers-FKPP
traveling waves are pushed for $\beta> 2$, and pulled for $\beta \leq 2$, as can be seen
from (\ref{jul1502}) and (\ref{jul1602}). 

While the case $\beta=2$ obviously fails the pushed front criterion (\ref{aug1102}), it has an additional property distinguishing it from the case $\beta < 2$.
The traveling wave, still given by
\be\label{aug1204}
\phi_2(x)=\farc{1}{1+e^x}
\ee
due to~(\ref{jul1502}), satisfies the weaker condition (cf.~\eqref{aug1102}, recalling that $c_* = 2$ when $\beta =2$) that
\be\label{aug1104}
	\phi_{2}(x) e^{x}
		= \farc{e^x}{1+e^x}
		\in L^\infty(\R).
\ee
As we will see, this reflects a very different long time
behavior of the solutions to (\ref{burgerskpp}) at $\beta=2$ compared both to $\beta<2$ and
$\beta>2$. Borrowing the terminology of~\cite{Lofting}, we will refer to such ``dual nature''
fronts at $\beta=\beta_c$ as ``pushmi-pullyu" fronts.

\subsubsection*{The large time behavior of the solutions }
  
We now describe the main results of this paper on the large time behavior of the solutions to 
the Burgers-FKPP  equation~(\ref{burgerskpp}) with $\beta\neq 0$ and rapidly decaying initial
conditions, generalizing the results on the standard Fisher-KPP equation (\ref{nov2502}) discussed above. 
The main new feature is the aforementioned transition from the pulled to pushed behavior at $\beta_c=2$.  The analysis at $\beta=\beta_c$ turns out to be surprisingly delicate.   

The study of the large time behavior in the present paper relies, in particular, on the notion of steepness of the solution.  While such arguments
date back to the original KPP paper~\cite{kolmogorov1937etude}, we were to a large extent motivated by 
the definition in the recent paper of Giletti and Matano~\cite{GM}.  As we will only need it for smooth functions, we can formulate
their notion as follows. Let us denote by~${\cal W}$ the class of $C^1$, monotonically decreasing function $u(x)$, $x\in \Rm$, such that
\be\label{jul1504}
\lim_{x\to-\infty} u(x)=1,~~\lim_{x\to+\infty}u(x)=0.
\ee
Let $\overline \cW$ be its closure in $L^1$.
%
Given two functions $u_{1},u_2\in{\cal W}$, we say that $u_1$ is steeper than $u_2$ if  
\be\label{jul1506}
|u_1'(u_1^{-1}(z))|>|u_2'(u_2^{-1}(z))|,~~\hbox{for all $z\in(0,1)$.} 
\ee
In other words, the graph of $u_1(x)$ is steeper than the graph of $u_2(x)$ when compared 
at each fixed level~$z\in(0,1)$, rather than at a fixed point $x\in\Rm$. This notion is translation
invariant; if $u_1$ is steeper than $u_2$, it is also steeper than any translate $u_2(\cdot+h)$,
with a fixed $h\in\Rm$.   For $u_1, u_2 \in \overline \cW$, we say that $u_1$ is steeper than $u_2$, if $u_1$ and $u_2$ can be approximated by $C^1$-functions $u_{1,\eps}$ and $u_{2,\eps}$ as $\eps \to0$ such that $u_{1,\eps}$ is steeper than $u_{2,\eps}$ for all $\eps$.

Equation (\ref{burgerskpp}) has the following important property.

\begin{prop}\label{prop-jul1502}
Let $u_1(t,x)$ and $u_2(t,x)$ be the solutions to (\ref{burgerskpp})  with the corresponding initial conditions $u_{10},u_{20}\in  \overline \cW$. 
If  $u_{10}$ is steeper than $u_{20}$, then 
$u_1(t,\cdot)$ is steeper than~$u_2(t,\cdot)$ for all $t>0$. 
\end{prop} 
This result was essentially proved for the classical Fisher-KPP equation (\ref{nov2502}) in the original KPP paper~\cite{kolmogorov1937etude}.
For the convenience of the reader, we present the proof for the Burgers-FKPP case in Section~\ref{sec:shape}.  
An interesting aspect here is that while  the steepness comparison is used throughout the paper, most of its applications are not in the spirit of the elegant intersection number type of arguments but to produce
estimates.
 A novel element of the use of these ideas in the present paper
is their quantitative application.  In contrast to the qualitative arguments used in previous works, we use ``steepness'' to obtain estimates on $u$ and integral quantities involving it.
An exception is the proof of Proposition~\ref{prop-jul1502} itself  and its consequences discussed in Section~\ref{sec:shape}.

The main result of the present paper is the following: 
\begin{thm}\label{thm:main} Let $u(t,x)$ be the solution to (\ref{burgerskpp}) with  the initial condition $u_{\rm in} \in \overline \cW$ such that $u_{\rm in}$ is steeper than the minimal speed traveling wave $\phi_\eta$.
Then, for each $\beta\leq2$, there exists a constant~$x_\infty$ that depends on $\beta$ and the initial
condition $u_{\rm in}$   so that
\begin{align}\label{nov2516}
    \lim_{t\to+\infty} u(t,x+m_\beta(t)) = \phi_\beta(x),
\end{align}
with the function $m_\beta(t)$ given by
\begin{align}\label{nov2518}
    m_\beta(t) = 2t-\frac{3}{2}\log (t+1) - x_{\infty}+o(1), ~\text{as} ~t\to +\infty,
\end{align}
if $\beta<2$, and, otherwise, by
\begin{align}\label{nov2520}
    m_{\beta=2}(t) = 2t - \frac{1}{2}\log (t+1)-x_{\infty}+o(1), ~\text{as} ~t\to +\infty.
\end{align}

For $\beta>2$, 
there exists $\omega>0$, which depends on $\beta$ but not on $u_{\rm in}$,  and $K>0$, which depends both on $\beta$ and $u_{\rm in}$, such that
\begin{align}\label{nov2522}
    \sup_{x\in\R}|u(t,x)-\phi_{\beta}(x-c_*t-x_\infty)|<Ke^{-\omega t}.
\end{align}
\end{thm}

 We note that the class of initial data considered in \Cref{thm:main} includes the Heaviside initial data: $u_{\rm in } = \1\left(x < 0\right)$.  Indeed, it is easy to see that $u_{\rm in}$ is the limit in $L^1$ as $\eps \to 0$ of 
\[
	u_{{\rm in},\eps} = \begin{cases}
			1	\qquad &\text{ if } x \leq -\eps\\
			\phi_\beta\big(\frac{x}{\eps^2 - x^2}\big)	\qquad&\text{ if } x\in (-\eps, \eps)\\
			0	\qquad&\text{ if } x \geq \eps,
		\end{cases}
\]
and $u_{{\rm in},\eps}$ is clearly steeper than $\phi_\beta$ for $\eps < 1$.

\Cref{thm:main} reflects the different nature of the Burgers-FKPP fronts we have discussed above for various values of $\beta\in\Rm$. 
For $\beta<2$, the solution is pulled and the front location has the same asymptotics (\ref{nov2518})
as (\ref{nov2508}) for the standard Fisher-KPP equation (\ref{nov2502}). For $\beta>2$, the solution
is pushed and the exponential-in-time convergence to the traveling wave (\ref{nov2522}) agrees with what we have seen
in (\ref{nov2510}) for pushed fronts.  The new asymptotics~(\ref{nov2520}) for the ``pushmi-pullyu" solutions at 
$\beta=\beta_c$ is different from both of these cases.

One may ask if the asymptotics (\ref{nov2518}) and (\ref{nov2520}) can be 
refined, as was done in~\cite{berestycki2017exact,berestycki2018new,Ebert-vanSaarlos,Graham,NRR2}
for the standard Fisher-KPP equation (\ref{nov2502}). 
It seems that this is possible using the fascinating, if formal, technique of~\cite{berestycki2018new}.  Arguing in the manner of~\cite{berestycki2018new}, it turns out that the $o(1)$ terms in the expansion of $m_\beta(t)$ are exactly the same as in the Fisher-KPP case when $\beta < 2$; however, they are different when $\beta = 2$.  In this ``pushmi-pullyu'' case, one finds
\begin{align}\label{aug1302}
m_2(t)=2t-\frac{1}{2}\log (t)-x_{\infty} -\frac{\sqrt{\pi}}{2\sqrt{t}}
+\frac{(1-\log 2)}{4} \frac{\log t}{t}+o\Big(\farc{\log t}{t}\Big).
\end{align} 
The first terms in this 
expansion above agree with the formal computation in~\cite{Ebert-vanSaarlos}, 
performed by a completely different method, 
where it was computed up to the order
$O(t^{-1/2})$. While the error in~\cite{Ebert-vanSaarlos} is given as $O(1/t)$, and
the~$O(\log t/t)$ correction seems to be 
missing, 
it is tempting to conjecture that~\eqref{aug1302} is universal for the ``pushmi-pullyu'' situations, 
such as, for instance, 
the reaction-diffusion equation~\eqref{aug1602} below.  The application 
of the formal technique of~\cite{berestycki2018new} to derive~\eqref{aug1302} 
is contained in \Cref{sec:higher-order}.

We should mention that the convergence of the solution to the Burgers-FKPP equation (\ref{burgerskpp})
to a traveling wave was studied using the matched asymptotic expansions in \cite{LeachHanac2016}. Their formal
results agree with Theorem~\ref{thm:main} in the cases~$\beta<2$ and $\beta>2$ but unfortunately overlook some of the details in the case $\beta=2$,
leading to an incorrect prediction. On the other hand, as we have mentioned, the
matched asymptotics analysis in~\cite{Ebert-vanSaarlos,Leach-Needham} 
predicts the shift~$(1/2)\log t$ for the ``pushmi-pullyu" transition
situations, such as (\ref{aug1602}) below, mimicking, in a certain sense, 
the Burgers-FKPP equation at $\beta=2$. 

%
%
%

\subsubsection*{Connection to other pulled to pushed transition problems}

The Burgers-FKPP equation at $\beta=2$ is not the only one example of a ``pushmi-pullyu" 
front. 
A well known instance is a particular generalized Fisher nonlinearity considered in~\cite{HadelerRothe}:
\be\label{aug1112}
u_t=u_{xx}+u(1-u)(1+au).
\ee
Here, the situation is qualitatively and even quantitatively 
extremely similar to the Burgers-FKPP picture: traveling waves exist for all $c\ge c_*(a)$, with $c_*(a)=2$ for all $0\le a\le 2$, and
\be\label{aug1114}
c_*(a)=\sqrt{\farc{a}{2}}+\sqrt{\farc{2}{a}},~~\hbox{for $a>2$,}
\ee
as in (\ref{speed-bis}). 
Even the traveling wave profile for $a=2$ is given exactly by $\phi_{\beta=2}$ in  (\ref{aug1204}).
 
 As explored in~\cite{Ebert-vanSaarlos},~\eqref{aug1112} is a special case of the more general class of equations
\be\label{aug1602}
	u_t=u_{xx}+u(1-u^n)(1+a u^n),
\ee
in which the pulled to pushed transition occurs at $a=n+1$: when $a \leq n+1$, the minimal speed is $c_* = 2$ and when $a> n+1$, $c_*  = \sqrt{\frac{a}{n+1}}+\sqrt{\frac{n+1}{a}}>2$.  In the critical case $a=n+1$, the corresponding traveling wave is explicit,
\be\label{aug1604}
\phi(x)=\big(1+e^{nx}\big)^{-1/n},
\ee
with the same purely exponential decay as in (\ref{aug1204}).

%

 Numerous other examples, including the pushed-pulled transition for systems of reaction-diffusion
equations arising in chemistry and biology are discussed in Section 3.13 of~\cite{vanSaarlos}.
We also mention the repulsive Keller-Segel-FKPP equation that was recently used to model the spread 
of a population in which individuals reproduce and diffuse, influenced by a preference 
for low population density regions with strength $|\chi|$. This is seen, for instance, in slime 
molds~\cite{slimemolds} and bacteria~\cite{bacteria}. Origins of some other similar chemotaxis models 
have been discussed in~\cite{FuGrietteMagal1, FuGrietteMagal2}.
The pushed-pulled transition
at the level of traveling waves has been shown in~\cite{henderson_chemotaxis}: 
there exists a threshold $\chi_0>0$ such that, when~$|\chi| < \chi_0$, the fronts are pulled, 
while, as $\chi \to -\infty$, fronts are pushed.  In the same vein is a system of equations considered in~\cite{BramburgerHenderson}: one equation is a Fisher-KPP equation for $T$ with drift $u$ and the other, which governs $u$, satisfies a Burgers type equation with a Boussinesq-type forcing depending on $T$.  Actually, this specializes to~\eqref{burgerskpp} for a certain choice of parameters.  As with the Keller-Segel-FKPP, a pushed-pulled transition occurs.  In both models, the pushed-pulled analysis in~\cite{BramburgerHenderson}
is performed only at the level of the traveling wave, which is a much simpler setting.  The Cauchy problem, on the other hand, appears to require new ideas.

A matched asymptotic analysis in~\cite{Ebert-vanSaarlos,Leach-Needham}
predicts the results as in Theorem~\ref{thm:main} to hold for the long time behavior 
of the solutions of  equations at the pushed-pulled transition, 
such as~(\ref{aug1602}). However, to the best of our knowledge, there 
are no rigorous results with this precision  in any of such critical cases. 
The best result in this direction
seems to be the very recent paper~\cite{Giletti},  which
shows that when $a=2$, the level sets of the solutions to (\ref{aug1112}) are located at a position
\[
m(t)=2t-\farc{1}{2}\log t+o(\log t),~~\hbox{ as $t\to+\infty$.}
\]
As   discussed in detail in~\cite{AS1} and~\cite{AS2}, 
the threshold cases that separate the pulled and pushed fronts 
seem to be also outside the scope of the currently available spectral methods.

\subsubsection*{Comments on the proofs}

Let us comment on the strategy of the proof of the three cases in
Theorem~\ref{thm:main}, in the order of increasing difficulty and intricacy. The case~$\beta>2$ falls into the category of pushed fronts, and the
proof follows the classical strategy of~\cite{rothe1981convergence,sattinger1976stability,sattinger1977weighted}, 
with appropriate modifications. 

For $\beta<2$ in Theorem \ref{thm:main}, we use an extension of the arguments in~\cite{HNRR,NRR1} for the Fisher-KPP equation, 
approximating the dynamics by the Dirichlet problem for the
linear heat equation on a half line. The Burgers drift term causes a difficulty, since the linearization strategy 
used, for instance, in \cite{Graham,henderson2016,HNRR,NRR1,NRR2} taking out the spatial exponential decay in the solution
seems insufficient for~$0<\beta<2$. To overcome this, we pass to the moving frame $x\to x-2t$, setting
\be\label{aug1214}
\hat u(t,x)=u(t,x+2t),
\ee
and introduce a weighted Hopf-Cole transform
\begin{align}\label{intro:trans}
v(t,x) = \exp\Big(x +\frac{\beta}{2}\int_x^{\infty} \hat u(t,y) dy\Big) \hat u(t,x),
\end{align}
combining the standard Hopf-Cole transform for the heat equation and the exponential weight used in the standard Fisher-KPP arguments. It turns out that 
if the initial condition $u_{\rm in}$ for (\ref{burgerskpp}) is steeper than the traveling wave $\phi_\beta$, then $v(t,x)$ satisfies a differential inequality
\be\label{aug1116}
v_t\le v_{xx},
\ee
as in Proposition~\ref{prop-jul1902} below. We stress that the steepness comparison of the initial condition to a traveling wave plays a crucial role in 
showing that $v(t,x)$ satisfies the differential inequality (\ref{aug1116}).

With (\ref{aug1116}) in hand, we are able to construct upper and lower barriers in the self-similar variables for the linearized equation for $v(t,x)$ on the half line, 
and then the convergence in the tail implies the convergence in the bulk due to the pulled-front nature of the dynamics, as in~\cite{NRR1}. 
Interestingly, this last step also utilizes the assumption that the initial condition, and hence the solution, 
is steeper than the minimal speed traveling wave, in an explicit quantitative way.  Qualitatively, the case $\beta<2$ is similar to the standard Fisher-KPP equation,
and the weighted Hopf-Cole transform gives a tool to see that. However, the repeated use of the steepness comparison is something new in this argument
for Burgers-FKPP equation.
 
The weighted Hopf-Cole transform also indicates the technical reason for the transition at~$\beta=2$:  while it is easy to see  from (\ref{intro:trans})
that~$v(t,x)\to 0$ as~$x\to-\infty$ for $\beta<2$ (recall that $\hat u \leq 1$), in the case
$\beta=2$ the function $v(t,x)$ approaches a positive constant as~$x\to-\infty$. This modifies the boundary condition for the linearized problem
for the upper and lower barriers in the
self-similar variables, and ultimately leads to the change in the logarithmic shift
from $(3/2)\log t$ to~$(1/2)\log t$ at $\beta=2$.

Let us now discuss the ingredients of the the proof of Theorem~\ref{thm:main} in the critical case~$\beta=2$, which is   remarkably different from 
the approach for the standard Fisher-KPP equation.  This analysis is probably the most novel part of the present paper. 
The first key observation is that when~$\beta=2$, the Burgers-FKPP equation (\ref{burgerskpp}) has a special structure: the function
\be\label{jul1510}
p(t,x)=e^{x}\hat u(t,x),
\ee
satisfies a spatially inhomogeneous conservation law:
\begin{align}\label{intro:para}
p_t+(e^{-x}p^2)_x = p_{xx}.
\end{align}
Here, $\hat u(t,x)$ is defined in (\ref{aug1214}). An immediate consequence of (\ref{intro:para}) is a conservation law for the exponential moment of $\hat u(t,x)$:
\be\label{aug1210}
\int e^{x}\hat u(t,x)dx=\int e^x u_{\rm in}(x)dx,~~\hbox{ for all $t>0$.}
\ee
This conservation law eventually leads to a lower bound for $m_2(t)$ of the form
\be\label{aug1312}
m_2(t)\ge 2t-\farc{1}{2}\log t +O(1),~~\hbox{as $t\to+\infty$,}
\ee
see the proof of Lemma~\ref{lem-may2502} in Section~\ref{s:upper} below. 

A matching upper bound for $m_2(t)$  is related to the behavior of $p(t,x)$. 
Note that, together with the explicit expression (\ref{aug1204}) for the profile $\phi_2(x)$, the convergence to a traveling wave in shape in~(\ref{nov2516}) yields, roughly,
\be\label{aug1216}
p(t,0)\approx \exp\big(-(2t-m_{2}(t))\big).
\ee
Thus, an upper bound of the form
\be\label{aug1314}
m_2(t)\le 2t-\farc{1}{2}\log t +O(1),~~\hbox{as $t\to+\infty$,}
\ee
would follow from an $L^\infty$-bound on $p(t,x)$ of the form
\be\label{aug1218}
p(t,x)\le \farc{C}{\sqrt{t}}.
\ee
Such decay, while natural to expect in view of~\eqref{intro:para},  is not automatic for solutions of mass-conserving
advection-diffusion equations, even if the advection is bounded: see the end of Section~\ref{sec:beta2_outline} for simple 
examples of such equations with solutions 
that do not decay in time.

The proof of (\ref{aug1218}), presented in Section~\ref{sec:twp-bounds}, turns out to be rather 
intricate. While (\ref{intro:para}) looks like a degenerate viscous conservation law,
we were unable to adapt the methods of~\cite{Carlen1996optimal} or \cite{howard2002} to~(\ref{intro:para}) and  instead 
take a different approach. The first step is a relative entropy computation inspired by~\cite{Const,MPP} where it was used for linear advection-diffusion 
equations.  An unusual twist is that
we compute the relative entropy not with respect to another solution but to a super-solution to (\ref{intro:para}). This leads to a weighted dissipation inequality for the function
\be\label{jul1512}
\vphi(t,x)=\farc{p(t,x)}{\rho(t,x)},
	\quad\text{ where }\rho(t,x)={1-u(t,x+2t)},
\ee
of the form
\be\label{jul1514}
\frac{d}{dt}\int \ \vphi^2(t,x)\rho(t,x)dx\le -2\int\vphi_x^2(t,x)\rho(t,x)dx. 
\ee
The dissipation identity (\ref{jul1514}) is similar to that for the standard heat equation, where it takes the form
\be\label{aug1118}
\frac{d}{dt}\int \ \vphi^2(t,x) dx\le -2\int\vphi_x^2(t,x) dx,
\ee
that is, as in~\eqref{jul1514} but without the weight $\rho(t,x)$. 
In the latter case, (\ref{aug1118}) combined with the Nash inequality and a standard duality argument directly
leads to the temporal decay rate $t^{-1/2}$ in~$\Rm$.  Here, the time-dependent weight $\rho(t,x)$ that appears
in (\ref{jul1514}) is degenerate as $x\to-\infty$, so the standard Nash inequality can not be used.
Instead, we obtain a Nash-type inequality for weighted spaces for a certain class of  degenerate weights: see Proposition~\ref{lem-may2508} below. 
The weights need to satisfy certain quantitative assumptions,
and we need to verify that the dynamics do not take the weight~$\rho(t,x)$, defined in (\ref{jul1512}), 
out of the class of the admissible weights or make the constants in the 
weighted Nash inequality in \Cref{lem-may2508} degenerate as $t\to+\infty$. 
Applying the weighted Nash inequality in~\eqref{jul1514},  leads to the appropriate decay of $\vphi(t,x)$ in a weighted~$L^2$-space.  
However, the nature of the Nash inequality leads to an extra delay 
in time, before the decay sets in, that depends on the initial
condition. This, among other issues, prevents us from using the duality argument to establish the  weighted~$L^\infty$-decay of $\vphi(t,x)$
from the $L^2$-estimate. 
Instead, it comes from additional ad hoc arguments, first establishing the bound at a large time $t'$ less than, but comparable to, $t$ and then extending it to the final time
$t$.  The decay of $p(t,x)$ follows from that of $\vphi(t,x)$ and bounds on the weights.

An extra technical complication is that the
function $\rho(t,x)$ appears in the denominator in the definition (\ref{jul1512}) of the function~$\vphi(t,x)$ but vanishes as~$x\to-\infty$.  
As a result, the $L^2$-norm of~$\vphi(t,x)$ is actually infinite for a large class of interesting initial conditions and extra approximations have to be used
to deal with this issue.
 At this stage, the assumption that $u_{\rm in}(x)=1$
for $x\le L_1$ is actually not a simplification but a complication that can not be avoided if one wants to include the Heaviside function into the class 
of admissible initial conditions.  The need to control the behavior in the back of the front is another reflection of the ``pushmi-pullyu" nature of the
solution. 

 

\subsubsection*{Organization of the paper}

This paper is organized as follows. Section~\ref{sec:shape} uses the steepness comparison arguments to prove
Proposition~\ref{prop-jul1502} and convergence of the solution to a traveling wave in shape. Section~\ref{sec:hopf-cole}
describes the weighted Hopf-Cole transform leading to the differential inequality (\ref{aug1116}). The proof of
Theorem~\ref{thm:main} for $\beta<2$ is contained in Section~\ref{sec:beta<2}. 
Section~\ref{sec:twp-bounds} is devoted to the 
decay estimates of the solutions to the inhomogeneous viscous conservation law (\ref{intro:para}) that appears in the case $\beta=2$. It is here
that we prove the aforementioned $L^\infty$-decay of the solutions to (\ref{intro:para}).   Section~\ref{sec:proof-beta=2} uses these
results to prove Theorem~\ref{thm:main} for $\beta=2$, bootstrapping the bounds (\ref{aug1312}) and (\ref{aug1314}) 
to the precise asymptotics~(\ref{nov2520}).
 The case $\beta>2$ is considered in Section~\ref{sec:beta>2}. Section~\ref{sec:higher-order}
uses the techniques of~\cite{berestycki2018new} to obtain, by formal arguments, further corrections to the front location asymptotics
given in Theorem~\ref{thm:main} for $\beta\le 2$. The result for $\beta<2$ is identical to the standard Fisher-KPP equation but is different for $\beta=2$. Finally, Appendix~\ref{sec:phase-plane} contains some basic facts about the traveling  waves for the 
Burgers-FKPP equation.

{\bf Acknowledgment.} JA was supported by Joe Oliger Fellowship. CH was partially
supported by NSF grant DMS-2003110. LR was partially supported by 
NSF grant DMS-1910023, and ONR grant N00014-17-1-2145. We are grateful to John Leach
for a very friendly and helpful discussion of the results in~\cite{LeachHanac2016} and~\cite{Leach-Needham}.

 \section{Convergence to a traveling wave in shape}\label{sec:shape} 
 
In this section, as a preliminary step to the proof of Theorem~\ref{thm:main},
we use a strategy inspired by 
the original KPP paper \cite{kolmogorov1937etude} to show convergence in shape 
of a solution to the Burgers-FKPP equation to a traveling wave. 
As the first step, we prove Proposition~\ref{prop-jul1502}.

\subsubsection*{The proof of Proposition~\ref{prop-jul1502}}

 It is enough to assume that $u_{10}, u_{20} \in \cW$ by standard density arguments. 
Let $u_1(t,x)$ and $u_2(t,x)$ be the solutions to (\ref{burgerskpp})  with the respective 
initial conditions $u_{10},u_{20}$ such that $u_{10}$ is steeper than $u_{20}$. 
First, we note that since the initial conditions are decreasing, both $u_1(t,x)$ and
$u_2(t,x)$ are decreasing and have the left and right limits as in (\ref{jul1504}),
so that both~$u_1(t,\cdot)$ and $u_2(t,\cdot)$ lie in $\cal W$.  

Recall the definition of ``steeper''~\eqref{jul1506}.  To show that $u_1(t,\cdot)$ is steeper than $u_2(t,\cdot)$ for any $t>0$, consider the functions 
\[
w(t,x;k_0)=u_1(t,x)-u_2(t,x+k_0),~~q(t,x)=u_1(t,x)+u_2(t,x+k_0),
\]
for a fixed $k_0\in\Rm$. 
The function $w(t,x;k_0)$ satisfies 
\begin{equation}\label{jul1606}
w_t+\farc{\beta}{2}(qw)_x=w_{xx}+w-qw,
\end{equation}
with the initial condition
\begin{equation}\label{jul1608}
w(0,x;k_0) = u_{10}(x)-u_{20}(x+k_0).
\end{equation}
Since $u_{10}$ is steeper than $u_{20}$, it is also steeper than $u_{20}(\cdot+k_0)$.
Therefore, there exists $x_0$ so that 
\[
\hbox{$w(0,x;k_0)>0$ for all $x<x_0$},
\]
and 
\[
\hbox{$w(0,x;k_0)<0$ for all $x>x_0$.}
\]
As $w(t,x;k_0)$ is a solution to the parabolic equation (\ref{jul1606}),  a consequence of \cite[Theorems A and B]{Angenent} is 
that the  function $w(t,x;k_0)$ has exactly one zero $y(t;k_0)$ for all~$t>0$.   Indeed, \cite{Angenent} shows that the number of zeros is nonincreasing and that a zero may only disappear at a time $t_0>0$ when two zeros ``collide.''  Hence,
$w(t,x;k_0)>0$ for all~$x<y(t;k_0)$ and~$w(t,x;k_0)<0$ for all $x>y(t;k_0)$, with $y(0;k_0)=x_0$. In addition, we   have 
\begin{equation}\label{jul1610}
\partial_x u_1(t,y(t;k_0))<\partial_x u_2(t,y(t;k_0)).
\end{equation}
Since this is true for all $k_0\in\Rm$, it follows that $u_1(t,\cdot)$ is steeper than 
$u_2(t,\cdot)$.~$\Box$

A standard approximation argument shows the following.
\begin{cor}\label{cor-jul1602}
Let $v(t,x)$ and $u(t,x)$ be the solutions to (\ref{burgerskpp}) with the 
respective initial conditions~$v_{\rm in}\in{\cal W}$ and $u_{\rm in}(x)=\one (x\le 0)$.
Assume that $v_{\rm in}(x)$ is steeper than the minimal speed traveling wave $\phi_\beta(x)$.
Then, for any $t>0$ the solution $v(t,\cdot)$ is steeper than $\phi_\beta$, and 
is less steep than $u(t,x)$.
\end{cor}

\subsubsection*{Convergence in shape}

We now establish convergence of the solution in shape to a traveling wave. 
\begin{prop}\label{prop:aug262}
Let $u(t,x)$ be the solution to (\ref{burgerskpp}) 
with the initial condition $u_{\rm in}\in{\cal W}$ that is steeper 
than the minimal speed traveling wave $\phi_\beta(x)$, or with $u_{\rm in}(x)=\one(x\le 0)$. 
Then, there exists a function $m_\beta(t)$ such that~$\dot m_\beta(t)\to c_*(\beta)$ 
as $t\to +\infty$ and
\begin{equation}\label{dec130}
	u(t,x+m_\beta(t))\to\phi_\beta(x)~~\hbox{as $t\to+\infty$, uniformly on $\Rm$.}
\end{equation}
Here, $\phi_\beta(x)$  is a solution to (\ref{nov2514}) with the minimal speed $c_*=c_*(\beta)$.  
\end{prop}
Corollary~\ref{cor-jul1602} shows that it suffices to consider the solution $u(t,x)$
to (\ref{burgerskpp}) with the initial condition 
$u(0,x)=\one(x\le 0)$. Note that for any $\tau>0$, the function $u^{(\tau)}(t,x)=u(t+\tau,x)$
is the solution to (\ref{burgerskpp}) with the initial condition $u^{(\tau)}(0,x)=u(\tau,x)$
that is less steep than $u(0,x)$. It follows that for any $t>0$ and $\tau>0$ the function
$u(t,\cdot)$ is steeper than $u(t+\tau,\cdot)$. In addition,~$u(t,\cdot)$ is steeper
than the minimal speed traveling wave $\phi_\beta(x)$ for all $t>0$. 
Hence, if for each~$v\in(0,1)$ and $t>0$, 
we let $x(t,v)$ be the unique point such that~$u(t,x(t,v))=v$, 
then, the function 
\be\label{jul1612}
E(t,v)=u_x(t,x(t,v))< 0,
\ee
is increasing in $t$ for all $v\in(0,1)$, 
and
\begin{equation}\label{dec132}
E(t,v)\le \bar E(v):=\phi_\beta'(\phi_\beta^{-1}(v)).
\end{equation}

Let now $m_\beta(t)$ be the position such that $u(t,m_\beta (t)) = 1/2$ for all $t>0$, and consider the translate 
\[
	\tilde u(t,x)=u(t,x+m_\beta(t)),
\] 
as well as  the corresponding inverse $\xi(t,v)$ defined by
$\tilde u(t,\xi(t,v))=v$, for $0<v<1$. Observe that~$\xi(t,1/2)=0$ for all $t>0$.
We see from (\ref{jul1612})-(\ref{dec132}) 
that the function $E(t,v)$ is negative and increasing in time. Thus, it has a limit
\begin{equation}\label{dec926}
	E(t,v)
		\to E_\infty(v)
		\le \bar E(v)=\phi_\beta'(\phi_\beta^{-1}(v))<0,~~\hbox{as $t\to+\infty$.}
\end{equation}
Hence 
\[
\pdr{\xi(t,v)}{v}=\farc{1}{E(t,v)}\to \farc{1}{E_\infty(v)},~~\hbox{as $t\to+\infty$,}
\]
and
\begin{equation}\label{dec310}
	\xi(t,v)
		=\int_{1/2}^v\pdr{\xi(t,v')}{v'}dv'
		\to \int_{1/2}^v\frac{dv'}{E_\infty (v')}
		:=\xi_\infty (v).
\end{equation}
As a consequence, the   function $\tilde u(t,x)$ also converges uniformly on compact sets
to a limit~$\tilde u_\infty(x)$:
\begin{equation}\label{dec314}
	\tilde u(t,x)\to\tilde u_\infty(x)\hbox{ as $t\to+\infty$,}
\end{equation}
with $\tilde u_\infty(x)$ determined by
\be\label{apr902}
	\xi_\infty(\tilde u_\infty(x))=x.
\ee
Moreover, due to~\eqref{dec926}, we have
\be\label{apr904}
	|\xi_\infty(v)|
		=\int_{1/2}^v\frac{dv'}{|E_\infty(v')|}
		\le \int_{1/2}^v\frac{dv'}{|\bar E(v')|}:=\bar\xi(v).
\ee
This yields the correct behavior of the limits $x\to\pm\infty$:
\begin{equation}\label{dec318}
	\tilde u_\infty (-\infty)=1,
	~~ \tilde u_\infty(+\infty)=0.
\end{equation}
Indeed, considering, for example the behavior $x\to +\infty$, we have
\[
	+\infty
		= \lim_{x\to+\infty} x
		= \lim_{x\to+\infty} \xi_\infty(\tilde u_\infty(x))
		\leq \bar \xi(\tilde u_\infty(x)).
\]
and $\bar \xi^{-1}(\infty) = 0$.  The argument for $x\to-\infty$ is similar.

Furthermore, as $u(t,x)$ is strictly decreasing in $x$ and $u_x(t,m_\beta(t))<0$, the function $m_\beta(t)$ is differentiable in $t$. Hence, $\tilde u(t,x)$ satisfies
\begin{equation}\label{dec312}
	\tilde u_t-\dot m_\beta(t)\tilde u_x+\beta \tilde u \tilde u_x=\tilde u_{xx}+\tilde u (1 - \tilde u).
\end{equation}
Notice that
\[
	\dot m_\beta(t)
		= -\frac{u_t(t,m_\beta(t))}{u_x(t,m_\beta(t))}
		= -\frac{u_t(t,m_\beta(t))}{\tilde u_x(t,0)}.
\]
By parabolic regularity theory, the numerator is bounded and, 
 by~\eqref{dec132} with $v=1/2$, the denominator is bounded away from zero.  It follows that $\dot m_\beta$ is bounded uniformly in $t$.   Hence, for any sequence $t_n \to \infty$, there is a subsequence $t_{n_k}\to\infty$ and a real number $c\in \R$ such that $\dot m_\beta(t_{n_k}) \to c$.  
Using then the convergence (\ref{dec314}), we deduce 
that 
\begin{equation}\label{dec316}
-c \partial_x \tilde u_\infty + \beta \tilde u_\infty \partial_x \tilde u_\infty
	= \partial_x^2 \tilde u_\infty+\tilde u_\infty(1 - \tilde u_\infty),
\end{equation}
where we have switched to $\partial$ notation to avoid the awkward double subscript. 

From~\eqref{dec316}, we see that $\tilde u_\infty(x)$ is a traveling wave solution to (\ref{burgerskpp}) moving
with the speed $c$. It remains to show that $c=c_*(\beta)$. The key point is that the 
steepness comparison argument above 
applies to any traveling wave solution to  
\begin{equation}\label{apr906}
	-c\phi_x+\beta \phi\phi_x=\phi_{xx}+\phi-\phi^2.
\end{equation}
In other words, if we set 
\[
	E_\phi(v)=\phi'(\phi^{-1}(v)),~~\hbox{ for $0<v<1$},
\]
then we know that 
\[
	E_\infty (v)\le E_\phi(v),
\]
for any $\phi$ that satisfies (\ref{apr906}) with some $c\ge c_*(\beta)$. Therefore, 
the limit $\tilde u_\infty(x)$ is the 
traveling wave that is the steepest among all traveling wave solutions. Lemma~\ref{lem-apr902} implies
that $\tilde u_\infty(x)=\phi_\beta(x)$ is the minimal speed traveling wave and, thus, $c=c_*(\beta)$.  By the arbitrariness of the sequence $t_n$, it follows that $\dot m_\beta(t) \to c_*(\beta)$ as $t\to\infty$. This finishes the proof of Proposition~\ref{prop:aug262}.~$\Box$

\section{The weighted Hopf-Cole transform}\label{sec:hopf-cole}

In this section, we discuss a weighted Hopf-Cole transform that will play a key role in
the analysis of the Burgers-FKPP equation for $\beta\le 2$.  
Let us recall that the standard Burgers equation 
\begin{equation}\label{dec323}
u_t+\beta uu_x=u_{xx}
\end{equation}
can be linearized by means of the Hopf-Cole transform. Namely, if $u$ is
a solution to (\ref{dec323}) then the function 
\begin{equation}\label{dec324}
v(t,x)=\exp\Big(\farc{\beta}{2}\int_x^{+\infty} u(t,y)dy\Big)
\end{equation}
satisfies the heat equation
\begin{equation}\label{dec325}
v_t=v_{xx}.
\end{equation} 
The second simple  observation is that if $\hat u(t,x)$ is the solution to the standard
Fisher-KPP equation in a frame moving with the speed $c_*=2$:
\begin{equation}\label{dec326}
\hat u_t-2\hat u_x=\hat u_{xx}+\hat u-\hat u^2,
\end{equation} 
then the function
\begin{equation}\label{dec328}
v(t,x)=e^{x}\hat u(t,x)
\end{equation}
satisfies
\be\label{jul1908}
v_t=v_{xx}-e^{-x}v^2.
\ee
The nonlinear term in (\ref{jul1908}) is negligible for $x$ very large and positive but plays the role of a large absorption 
for $x$ very negative. Therefore, the solution to (\ref{jul1908}) should be well approximated by the solution of the heat equation
on a half-line $x>0$ with the Dirichlet boundary condition:
\be\label{jul1910}
\bal
&v_t=v_{xx},~~x>0,\\
&v(t,0)=0.
\enbal
\ee
This simple idea is what is driving the convergence to a traveling wave in~\cite{Graham,HNRR,NRR1,NRR2}. 

The weighted Hopf-Cole transform that we discuss below allows us to adapt this intuition to the Burgers-FKPP equation (\ref{burgerskpp})
with $\beta\le 2$, and also shows why the transition from pulled to pushed fronts happens at $\beta=2$. 

We will consider the solution to~(\ref{burgerskpp}) in the reference frame
\begin{equation}\label{dec414}
\tilde u(t,x)=u(t,x+m_\beta(t)),
\end{equation}
centered at 
\begin{equation}\label{dec416}
m_\beta(t)=2t-\farc{r(\beta)}{2}\log(t+1).
\end{equation}
Here, we take
\be\label{jul1904}
r(\beta)=\left\{\begin{matrix} 3,~~\hbox{ if $\beta<2$},\cr 1,~~\hbox{ if $\beta=2$},\cr \end{matrix}\right.
\ee
in accordance with the different behavior in Theorem~\ref{thm:main} in these two cases. 
In the above reference frame,   (\ref{burgerskpp})
takes the form
\begin{equation}\label{dec330}
\tilde u_t-\Big(2-\farc{r(\beta)}{2(t+1)}\Big)\tilde u_x+\beta \tilde u\tilde u_x=\tilde u_{xx}+
\tilde u-\tilde u^2.
\end{equation}
Motivated by (\ref{dec324}) and (\ref{dec328}), we introduce the  weighted Hopf-Cole transform
\begin{align}\label{dec332}
    v(t,x) = \exp(\Gamma(t,x))\tilde u(t,x),~~
    \Gamma(t,x)=x+\frac{\beta}{2}\int_x^{+\infty} \tilde u(t,y) dy, 
    \end{align}
that is a combination of (\ref{dec324}) and (\ref{dec328}). 

Let us point out one way to see why $\beta = 2$ is the critical value by using $v$.  Since $\tilde u$ is steeper than $\phi$ and $\phi$ converges exponentially to $1$ as $x\to-\infty$, we find
\[
	- C_t - \frac{\beta}{2} x
		\leq \frac{\beta}{2} \int_x^\infty \tilde u(t,y) dy
		\leq C_t - \frac{\beta}{2} x.
\]
The constant $C_t$ may depend on $t$, but not on $x$.  From this, we see that, as $x\to-\infty$,
\[
	v(t,x) = e^{O(1) + x \left(1 - \frac{\beta}{2}\right)}
		\to \begin{cases}
			0 \qquad &\text{ if } \beta < 2,\\
			e^{O(1)} \qquad &\text{ if } \beta = 2,\\
			\infty \qquad &\text{ if } \beta > 2.
		\end{cases}
\]
These differences reflect the three different behaviors in \Cref{thm:main}: when $\beta < 2$, the nonlinear term (integral of $\tilde u$) does not dominate, and when $\beta > 2$, the nonlinear term dominates.
 
The main result of this section is
the following analogue of (\ref{jul1908}) in the standard Fisher-KPP case.  It will allow us to adapt an
approximation similar to the linear Dirichlet boundary problem~(\ref{jul1910}) for~$\beta<2$ in Section~\ref{sec:beta<2}
and with a different boundary condition for~$\beta=2$ in Section~\ref{sec:proof-beta=2}. This will be extremely
important for the proof of Theorem~\ref{thm:main} for $\beta\le 2$. 
 \begin{prop}\label{prop-jul1902}
Let $u(t,x)$ be the solution to (\ref{burgerskpp}) with $\beta\le 2$ and the initial
condition $u(0,x)$ as in Theorem~\ref{thm:main}. Then, the function $v(t,x)$ defined 
in (\ref{dec332}) satisfies the differential inequality
\begin{equation}\label{dec404}
v_t - v_{xx}+\frac{r(\beta)}{2(t+1)}(v_x-v) \le 0. 
\end{equation}
\end{prop}
As we will see in the proof, it is here, among other places, that the steepness assumption on the initial condition $u(0,x)$
plays a crucial role, together with propagation of steepness in Proposition~\ref{prop-jul1502}. 

\subsubsection*{Proof of Proposition~\ref{prop-jul1902}}

We claim that the function $v(t,x)$
satisfies an equation of the form
\begin{equation}\label{dec333}
v_t - v_{xx}+\frac{r(\beta)}{2(t+1)}(v_x-v) = -G(t,x;\tilde u)v,
\end{equation}
where
\begin{equation}\label{dec334}
G(t,x;\tilde u)=\tilde u(t,x)-\frac{\beta}{2}\int_x^{+\infty} \tilde u(t,y)(1-\tilde u(t,y))dy.
\end{equation}
Let us verify that (\ref{dec333}) holds. We compute
\begin{equation}\label{dec340}
v_t = e^{\Gamma}\tilde u_t + \frac{\beta v}{2}\int_x^{+\infty} \tilde u_t(t,y)dy
	\quad\text{ and }\quad
v_x = e^{\Gamma}\tilde u_x + \Big(1-\frac{\beta}{2}\tilde u\Big)v,
\end{equation}
so that 
\begin{equation}\label{dec336}
\begin{aligned}
v_{xx} &= e^{\Gamma}\tilde u_{xx}+\Big(1-\frac{\beta}{2}\tilde u \Big)e^{\Gamma}\tilde u_x 
- \frac{\beta}{2}\tilde u_x v + \Big(1-\frac{\beta}{2}\tilde u\Big)v_x \\
&= e^{\Gamma}\tilde u_{xx} + 2\Big(1-\frac{\beta}{2}\tilde u\Big)v_x 
-\Big(1-\frac{\beta}{2}\tilde u\Big)^2v -\frac{\beta}{2}\tilde u_x v.
\end{aligned}
\end{equation}
Using these identities in (\ref{dec330}) gives
\begin{equation}\label{july2100bis}
\begin{aligned}
v_t &-\frac{\beta v}{2}\int_x^{+\infty} \tilde u_t(t,y) dy 
-\Big(2-\frac{r(\beta)}{2(t+1)}\Big)\Big(v_x- \Big(1-\frac{\beta}{2}\tilde u\Big)v\Big) + \beta v\tilde u_x\\
&=v_{xx}-2\Big(1-\frac{\beta}{2}\tilde u\Big)v_x +\Big(1-\frac{\beta}{2}\tilde u\Big)^2 v  
+ \frac{\beta}{2}v\tilde u_x+ v -v\tilde u.
\end{aligned}
\end{equation}
To simplify this equation, we integrate (\ref{dec330}) from $x$ to $\infty$ to get
\begin{align}\label{dec338}
\int_x^{+\infty} \tilde u_t(t,y) dy = 
- \Big(2-\frac{r(\beta)}{2(t+1)}\Big) \tilde u+\frac{\beta}{2}\tilde u^2 -\tilde u_x 
+\int_x^{+\infty} \tilde u(1-\tilde u)dy.
\end{align}
Substituting this back into (\ref{july2100bis})  gives, after some algebra
\begin{equation}\label{dec339}
\begin{aligned}
v_t &-v_{xx} +\frac{r(\beta)}{2(t+1)}\big(v_x-v\big)  + \beta v\tilde u_x
\\&= \frac{\beta v}{2} \int_x^{+\infty} \tilde u(1-\tilde u)dy 
+ {\beta} \tilde u v_x +\frac{\beta^2}{2}\tilde u^2 v -\beta \tilde uv -v\tilde u. 
\end{aligned}
\end{equation}
Using the second identity in  (\ref{dec340}) in the right side gives
\begin{equation}\label{dec341}
\begin{aligned}
v_t &-v_{xx} +\frac{r(\beta)}{2(t+1)} \big(v_x- v\big) =  
\frac{\beta v}{2} \int_x^{+\infty} \tilde u(1-\tilde u)dy 
-v\tilde u, 
\end{aligned}
\end{equation}
which is exactly (\ref{dec333})-(\ref{dec334}).

Here is the key observation.
\begin{lemma}\label{lem:neg}
If $\beta\le 2$ and the initial condition $u(0,x)$ is as in Theorem~\ref{thm:main},
then $G(t,x;\tilde u)\ge 0$ for all $t>0$ and $x\in\Rm$.
\end{lemma}
{\bf Proof.} The claim is trivially true for 
 $\beta\leq 0$, so we only consider the case $0<\beta\le 2$. Let us first show that a traveling wave $\phi_\beta(x)$ satisfies
the inequality
\be\label{aug2310bis}
\farc{\beta}{2}\int_{x}^{+\infty}\phi_\beta(y)(1-\phi_\beta(y))dy\le \phi_\beta(x),~~\hbox{ for all $x\in\Rm$.}
\ee
We will prove (\ref{aug2310bis}) for $\beta<2$, and deduce the conclusion for $\beta=2$ by continuity.
Recall that for~$\beta<2$ the traveling wave has the asymptotics (\ref{asymp1}):
\[
\phi_\beta(x)\sim (Ax+B)e^{-x},~~\hbox{ as $x\to+\infty$,}
\]
with some constants $A$ and $B$. As $\beta<2$, it follows that there exists $L>0$ so that (\ref{aug2310bis}) holds for all $x\ge L$,
and we only need to prove this inequality for $x< L$. 

Next, we integrate the traveling wave equation, with $c_*(\beta)=2$:
 \begin{equation}\label{july1100bis}
\begin{aligned}
    -&\phi_\beta'' -2\phi_\beta' + \beta\phi_\beta\phi_\beta' = \phi_\beta - \phi_\beta^2,\\
    &\phi_\beta(-\infty) = 1, ~~\phi_\beta(+\infty) = 0,
\end{aligned}
\end{equation}
from $x$ to $+\infty$ to get
\begin{align}\label{sep75}
    \phi_\beta'(x)+2\phi_\beta(x)-\frac{\beta}{2}\phi_\beta^2(x) = \int_x^{+\infty} \phi_\beta(y)(1-\phi_\beta(y)) dy.
\end{align}
Passing to the limit $x\to-\infty$ in (\ref{sep75}), 
and keeping in mind that $\phi_\beta(-\infty)= 1$, gives 
\begin{align*}
    \frac{\beta}{2}\int_{-\infty}^{+\infty}\phi_\beta(y)(1-\phi_\beta(y)) dy -\phi_\beta(-\infty) = -\frac{\beta^2}{4}+\beta -1 = -(\frac{\beta}{2}-1)^2<0.
\end{align*}
Therefore, there exists $L_1<0$ so that (\ref{aug2310bis}) holds for all $x<L_1$. In order to show that this inequality also holds for
$L_1<x<L$, consider the function
\begin{equation}\label{dec344}
G(x,\phi_\beta)=\phi_\beta(x)-\frac{\beta}{2}\int_x^{+\infty}\phi_\beta(y)(1-\phi_\beta(y))dy.
\end{equation}
Note that
\[
G'(x,\phi_\beta)=\phi_\beta'(x)+\frac{\beta}{2} \phi_\beta(1-\phi_\beta).
\]
Using (\ref{july1100bis}), we obtain
\[
G''(x,\phi_\beta) = \phi_\beta''+\frac{\beta}{2}\phi_\beta' -\beta \phi_\beta \phi_\beta' =\phi_\beta^2-\phi_\beta-2\phi_\beta'
+\frac{\beta}{2}\phi_\beta'.  
\]
Therefore, if $x_0$ is a critical point of $G(x,\phi_\beta)$, then
\[
\phi_\beta'(x_0)=-\frac{\beta}{2} \phi_\beta(x_0)(1-\phi_\beta(x_0)),
\]
and
\[
G''(x_0,\phi_\beta)=\farc{2}{\beta}\phi_\beta'-2\phi_\beta'
+\frac{\beta}{2}\phi_\beta'=\Big(\farc{2}{\beta}+\farc{\beta}{2}-2\Big)\phi_\beta'<0.
\]
Hence, the only possible critical points of $G(x,\phi_\beta)$ are local maxima. 
As $G(x)>0$ for $x>L$ and~$x<L_1$ we deduce that $G(x,\phi_\beta)>0$ for all
$x\in\Rm$, and (\ref{aug2310bis}) holds.

To finish the proof of Lemma~\ref{lem:neg}, consider the integral
\be\label{aug2304}
I(t,x)=\int_{x}^\infty \tilde u(t,y)(1-\tilde u(t,y))dy. 
\ee
To write this integral differently, define $\tilde x$ by $\tilde u(t,\tilde x(t,v)) = v$ and note that
\be
	\tilde u_x(t,\tilde x(t,v))=E(t,v),~~\hbox{ for $0<v<1$},
\ee
with $E(t,v)$ as in (\ref{jul1612}). 
%
Then, making the change of variables $y \mapsto v$ via $y=\tilde x(t,v)$, we find 
\be\label{aug2304bis}
\begin{aligned}
I(t,x)&=\int_{x}^\infty \tilde u(t,y)(1-\tilde u(t,y))dy
=\int_0^{\tilde u(t,x)}v(1-v)\farc{dv}{|E(t,v)|}\le 
\int_0^{\tilde u(t,x)}v(1-v)\farc{dv}{|\bar E(v)|}\\
&=
\int_{\phi_\beta^{-1}(\tilde u(t,x))}^\infty\phi_\beta(y)(1-\phi_\beta(y))dy,
\end{aligned}
\ee
with $\bar E(v)$ defined in (\ref{dec132}).
We used (\ref{dec132}) in the inequality in the first line in (\ref{aug2304bis}),
employing the assumption that the initial condition is steeper than the minimal
speed traveling wave. 
Using~(\ref{aug2310bis}) in (\ref{aug2304bis}) we see that
\be\label{dec347}
\begin{aligned}
I(t,x)&\le\int_{\phi_\beta^{-1}(\tilde u(t,x))}^\infty\phi_\beta(y)(1-\phi_\beta(y))dy\le \frac{2}{\beta} \phi_\beta(\phi_\beta^{-1}(\tilde u(t,x)))= \frac{2}{\beta} \tilde u(t,x).
\end{aligned}
\ee
We conclude that
\[
G(t,x;\tilde u)=\tilde u(t,x)-\frac{\beta}{2} I(t,x)\ge 0,
\]
finishing the proof of Lemma~\ref{lem:neg}, and hence that of Proposition~\ref{prop-jul1902} as well.~$\Box$

\section{The proof of Theorem~\ref{thm:main} for $\beta<2$}\label{sec:beta<2}

\subsubsection*{Outline of the proof}

In this section, we prove  Theorem~\ref{thm:main} for $\beta<2$. The overall
strategy of the proof is similar to~\cite{NRR1} which considers the classical
Fisher-KPP equation with $\beta=0$ but there are also non-trivial differences worth
mentioning.

The first step is to get control of the solution on the spatial scales~$x\sim O(\sqrt{t})$. 
This is done using  the self-similar variables. 
The estimates are precise enough to include the tail behavior
of~$\tilde u(t,x)$ on the intermediate 
scales~$x\sim O(t^\gamma)$, with $\gamma\in(0,1/2)$ that are
between the traveling wave scale~$x\sim O(1)$ and the diffusive scale $x\sim O(\sqrt{t})$. 
Unlike in~\cite{NRR1}, in this step we rely crucially on the weighted Hopf-Cole transform and
Proposition~\ref{prop-jul1902} to construct upper and lower barriers for the
solution.  The main estimate is the following:
\begin{lemma}\label{lem:close}
For $\beta<2$, let $\tilde u(t,x)$ be the solution to (\ref{dec330}) 
with the initial condition as in Theorem~\ref{thm:main}. 
There exist $\alpha_\infty>0$ and $\eps_0>0$ so that, for any~$0<\gamma<1/2$
and $\eps\in(0,\eps_0)$, 
there exists~$T_{\eps,\gamma}>0$ such that
\begin{align}\label{aug279}
	|\tilde u(t,x_{\gamma})
		- \alpha_{\infty}x_{\gamma}e^{-x_{\gamma}}|
	\le  {\eps} 
	 x_{\gamma}e^{-x_{\gamma}}
		\qquad \text{ for all } t > T_{\eps,\gamma},
\end{align}
where $x_{\gamma} = (t+1)^{\gamma}$.  
\end{lemma}

The second step is to use the pulled nature of the problem to show that the control of $\tilde u$ given in \Cref{lem:close} at $x=x_\gamma(t)$ 
induces convergence to a traveling wave 
on the scales $x\sim O(1)$.  Before discussing the modifications required for $\beta\neq 0$, we recall the argument in~\cite{NRR1}.  It proceeds by constructing 
solutions~$\tilde u_\alpha$ to the same equation 
as satisfied by $\tilde u$ (which, in~\cite{NRR1} is~(\ref{dec330}) with~$\beta = 0$), considered 
on the half-line $(-\infty,(t+1)^\gamma)$, with the boundary condition
\[
\tilde u_\alpha(t, (t+1)^{\gamma}) = \alpha x_\gamma e^{-x_\gamma},
\]
and with $\alpha = \alpha_\infty \pm \eps$ (cf.~\cite[Section~4]{NRR1}).   Then, the analogue of (\ref{aug279}) for $\beta=0$, which is Lemma~5.1 in~\cite{NRR1}, 
and the comparison principle imply that
\[
\tilde u_{\alpha_\infty-\eps} \leq \tilde u \leq \tilde u_{\alpha_\infty+\eps}.
\]
One concludes after proving the convergence of $\tilde u_\alpha(t,x)$, as $t\to+\infty$, to a suitable shift $\vphi_\alpha(x)$ of the traveling wave, showing that 
the shifts for the waves $\vphi_{\alpha_\infty\pm\eps}$ 
are $O(\eps)$ apart, and, finally, letting~$\eps\to 0$. 

Unfortunately, a direct attempt to do this in our setting fails for several reasons.   
The convergence of $\tilde u_\alpha$ to $\vphi_\alpha$ in~\cite{NRR1} is achieved by the construction of an explicit super-solution for the equation 
for the difference $e^x(\tilde u_\alpha - \vphi_\alpha)$.  However, the Burgers term in our context leads to a growth term in this equation when $\beta > 0$ 
(see, in contrast,~\eqref{sep77} when $\beta \leq 0$).  We bypass this issue by working with the weighted Hopf-Cole transforms of $\tilde u$ and the traveling wave.  
One might be tempted to define $\tilde u_\alpha$ as above and then take its weighted Hopf-Cole transform. However,  there is no apparent reason for 
$\tilde u_\alpha$ to be steeper 
than the traveling wave, meaning that \Cref{lem:neg} does not apply, and we cannot deduce 
the key differential inequality~\eqref{dec404}.  To bypass this difficulty, we work at the level of the weighted Hopf-Cole transform, defining $v_\alpha$ solving~\eqref{dec333} 
treated as a linear equation with the $\tilde u$ terms serving as given coefficients (see~\eqref{jul2806}--\eqref{jul2814} below), with a suitable boundary condition at $x=(t+1)^\gamma$.  
Proceeding as in~\cite{NRR1}, we obtain an upper bound on $v_\alpha$ given by a shift of the weighted Hopf-Cole transform of the traveling wave and a decaying term.  
The lower bound is obtained using yet another steepness comparison.  Afterwards, an additional argument is needed to upgrade this to the convergence of $\tilde u$ due to the fact 
that $\tilde u$ and $v$ are connected in a nonlocal fashion.

Below, we first prove Lemma~\ref{lem:close}, then apply it to show closeness of the weighted Hopf-Cole transforms of $\tilde u$ 
and the traveling wave $\phi_\beta$ in \Cref{lem:93}, and, finally, deduce the closeness of $\tilde u$ and~$\phi_\beta$ from this. 
For the sake of concreteness, we take $\tilde u(0,x)=\one(x\le 0)$. The argument for general initial conditions as in Theorem~\ref{thm:main} is nearly verbatim.

\subsection*{Analysis in the self-similar variables }

We start with equations  (\ref{dec333})-(\ref{dec334})
and pass to self-similar variables: let
\be\label{dec702}
	\tau = \log (t+1),
	~~\eta = \frac{x}{\sqrt{t+1}}
    	~~\text{and}~~\omega(\tau, \eta)=v\big(e^\tau-1,\eta e^{\tau/2}\big)e^{-\tau/2}.
\ee
In these variables, (\ref{dec333})-\eqref{dec334} becomes
\begin{align}\label{jun21757}
\omega_{\tau}+\cL\omega + \frac{3}{2}e^{-\tau/2}\omega_{\eta} 
=e^{\tau}\omega \Big(\frac{\beta}{2}\int_{\eta e^{\tau/2}}^{\infty} \tilde u(1-\tilde u) dy-
\tilde u(e^\tau-1, \eta e^{\tau/2})\Big),~~\tau>0,~\eta\in\Rm,
\end{align}
with the operator $\cL$ defined by 
\begin{align}\label{dec1128}
    \cL\omega := -\omega_{\eta\eta}-\frac{\eta}{2}\omega_{\eta} -\omega.
\end{align}
Note that when $\beta=0$, which is the classical Fisher-KPP equation,
(\ref{jun21757}) reduces to a local equation
\begin{align}\label{dec712}
\omega_{\tau}+\cL\omega + \frac{3}{2}e^{-\tau/2}\omega_{\eta} 
=-e^{3\tau/2-\eta \exp(\tau/2)}\omega^2,~~\tau>0,~\eta\in\Rm.
\end{align}
In that case, the nonlinear term in the right side is very small for $\eta>0$ but plays the role
of a large absorption for $\eta<0$. This was used in~\cite{NRR1} to show that (\ref{dec712})
is well-approximated by the linear problem 
\begin{equation}\label{dec714}
\bar\omega_\tau+\cL\bar\omega + \frac{3}{2}e^{-\tau/2}\bar\omega_{\eta}=0,~~\eta>0,
\end{equation}
augmented with the Dirichlet boundary condition $\bar\omega(\tau,0)=0$. This is the main intuition behind the proof of the long-time asymptotics of $\omega$ 
 in~\cite{NRR1}.

We now collect the ingredients that would allow us to use the arguments of~\cite{NRR1}.  
First, Lemma~\ref{lem:neg} implies that, as long as the initial condition $u(0,x)$ is steeper than
the minimal speed traveling wave, the right hand side of~\eqref{jun21757} is non-positive.
Hence, the solution to~\eqref{dec714} is still  a super-solution to~\eqref{jun21757}.

Second, the coefficient in the parenthesis in the right side of~\eqref{jun21757} can be bounded below as follows: by \Cref{lem:neg}, we have that
\[
	\int_{\eta e^{\tau/2}}^\infty \tilde u(1-\tilde u) dy
		\leq \frac{2}{|\beta|} \tilde u,
\]
leading to 
\be\label{e.c7302}
	e^\tau \Big( \frac{\beta}{2} \int_{\eta e^{\tau/2}}^\infty \tilde u(1-\tilde u) dy - \tilde u\Big)
		\geq - 2e^\tau \tilde u.
\ee
It follows that the solution to
\be\label{e.c7304}
	\under\omega_\tau + \cL \under\omega + \frac{3}{2} e^{-\tau/2}\under\omega_\eta
		= - 2 e^\tau \tilde u \under\omega
\ee
is a sub-solution to~\eqref{jun21757}.

The third observation is that 
the right side of~\eqref{jun21757} is very small for $\eta \gg e^{-\tau/2}$.  To see this, we first show that, if $\beta < 2$, then there exist $A>1$ and~$L>0$ so that
\begin{equation}\label{dec321}
u(t,x)\leq \bar u(t,x):=\frac{A}{1+e^{x-2t-L}},~~\hbox{ for all $x\in\Rm$}.
\end{equation} 
Indeed, the function $\bar u(t,x)$ satisfies (we set $L=0$ momentarily to simplify the notation)
\[
\bar u_t = \frac{2Ae^{x-2t}}{(1+e^{x-2t})^2},~~~
\bar u_x = -\frac{Ae^{x-2t}}{(1+e^{x-2t})^2},~~~
\bar u_{xx} = -\frac{Ae^{x-2t}(1-e^{x-2t})}{(1+e^{x-2t})^3}.
\]
Hence, as long as $\beta < 2$, we may choose $A\in(1,2/\beta)$ so that  
\begin{align}\label{dec704}
\bar u_t+\beta \bar u\bar u_x-\bar u_{xx}-\bar u(1-\bar u ) = A\frac{(2-\beta A)e^{x-2t}}{(1+e^{x-2t})^3}\geq 0,
\end{align}
With this choice of $A$, the function 
$\bar u(t,x)$ is a super-solution to (\ref{burgerskpp}) for any~$L\in\Rm$. 
Since $A>1$, we can choose 
$L>0$ sufficiently large so that~$\one(x\le 0)\le \bar u(0,x)$ for all~$x\in\Rm$. Then, the comparison principle for (\ref{burgerskpp}) implies that (\ref{dec321}) holds. 
As a consequence, we see that the right side of~\eqref{jun21757} can be bounded by
\begin{equation}\label{dec708}
\tilde u(t,x)\leq \min\Big\{1, \frac{A}{1+e^{x-({3}/{2})\log (t+1)-L}}\Big\},
\end{equation}
whence 
\be\label{calb}
	\int_x^\infty \tilde u(t,y)dy\le \begin{cases}
		C + \farc{3}{2}\log (t+1) - x
			\qquad &\text{ if } x \leq \frac{3}{2}\log(t+1)\\
		C(t+1)^{3/2} e^{-x}
			\qquad &\text{ if } x \geq \frac{3}{2}\log(t+1)
		\end{cases}
\ee
which are double-exponentially small on the scales $\eta = O( e^{-(\frac{1}{2}-\gamma)\tau})$ for any $\gamma \in (0,1/2)$ (recall that $x = \eta e^{\tau/2}$).

Finally, the asymptotics in~\eqref{calb} yield the approximate Dirichlet boundary condition for $\omega$ at~$\eta = -e^{-(1/2-\gamma)\tau}$.  Indeed, we have
\be\label{e.c7305}
	\begin{split}
	\omega(\tau, -e^{-(1/2-\gamma)\tau})
		&= e^{-\tau/2} \tilde u(e^\tau-1, -e^{\gamma \tau})
			\exp\left\{- e^{\gamma\tau} + \frac{\beta}{2}\int_{-e^{\gamma\tau}}^\infty \tilde u dy\right\}
		\\&\leq C \exp\left\{\tau - e^{\gamma\tau}\Big(1 - \frac{\beta}{2}\Big)\right\}.
	\end{split}
\ee

In order to analyze the long-time behavior of $\omega$, the main point is the following.  The above arguments show that $\omega(\tau,\eta)$ should be well-approximated by
the solution to the linear problem~(\ref{dec714}) with the Dirichlet boundary condition. 
The linear operator $\cL$, given by~(\ref{dec1128}), is compact
and self-adjoint  on~$H^1_0(e^{\eta^2/4}d\eta;\Rm_+)$.
Its spectrum consists of the eigenvalues $0,1,2,\dots$, and its principal eigenfunction is $\eta e^{-\eta^2/4}$  (in general, the eigenfunctions are given by the odd Hermite polynomials).  
Hence, the dominant behavior for $\omega(\tau,\eta)$ as $\tau\to+\infty$
should be given by $\alpha_\infty \eta e^{-\eta^2/4}$, for some $\alpha_\infty$ depending on the initial data.  This simple picture 
is complicated by the error terms in~\eqref{jun21757} and the ``not quite zero" boundary condition~\eqref{e.c7305}.  However, they have ``fast'' decay, so, with careful analysis, 
they can be suitably controlled.  This is the argument in a nutshell, even though the details of the proof are more intricate.

In order to carry out this strategy, the authors of~\cite{NRR1} require exactly the four ingredients listed above: equations for the super- and sub-solutions (given in
our case by~\eqref{dec714} and~\eqref{e.c7304}), double-exponential decay of the coefficients in the right side of~\eqref{jun21757} 
(see~\eqref{dec708} and~\eqref{calb}), and the approximate Dirichlet boundary condition, as in~\eqref{e.c7305}.  Thus, the strategy of that paper, which involves constructing successively more precise sub- and super-solutions of $\omega$ using the spectral properties of $L$, can be applied without alteration.  As such, we state the following lemma giving the asymptotics of $\omega$ and omit the details.

\begin{lemma}\label{lem:omega}
Let $\omega(\tau, \eta)$ be the solution of (\ref{jun21757}) on $\R$, with the 
initial condition~$\omega_0(\eta)$ such that~$\omega_0(\eta) = 0$ for all $\eta>A$ for some $A>0$, 
and $\omega_0(\eta) = O(e^{(1-\beta/2)\eta})$ for $\eta<0$. Then, for~$\beta<2$, there exists a constant $\alpha_{\infty}>0$ and functions $h$ and $R$ such that 
\begin{align}\label{dec1127}
	\omega(\tau, \eta)
		= (\alpha_{\infty} + h(\tau)) \eta e^{-\eta^2/4}
			+ R(\tau,\eta)e^{-\eta^2/6},
			~~~~ \eta\geq 0,
\end{align}
and, for any $\gamma'  \in (0,1/2)$,
\[
	\lim_{\tau\to\infty} h(\tau) = 0
	\quad\text{ and }\quad
	|R(\tau,\eta)|
		\leq C_{\gamma'} e^{- (1/2 - \gamma') \tau}.
\]
%
\end{lemma}


\subsubsection*{The proof of Lemma~\ref{lem:close}}

We now establish \Cref{lem:close} using~\Cref{lem:omega}.  
Let us recall, from~\eqref{dec332} and~\eqref{dec702}, that
\begin{align*}
v(t,x)=\exp\Big(x+\frac{\beta}{2}\int_x^{\infty}\tilde u(t,y) dy \Big) \tilde u(t,x) 
= \sqrt{t+1}\omega\Big(\log (t+1), \frac{x}{\sqrt{t+1}}\Big).
\end{align*}
It follows from (\ref{calb}) that  
we can take $T_\eps$ sufficiently large, so that when $t>T_\eps$, we have
\be\label{e.c6293}
	\begin{split}
		&e^{-C (t+1)^{3/2} e^{-(t+1)^\gamma}}\sqrt{t+1}\omega\Big(\log (t+1), \frac{x_\gamma}{\sqrt{t+1}}\Big)
	\leq e^{- \frac{\beta}{2} \int_{x_\gamma}^\infty \tilde u dy} v(t, x_\gamma)
	=
		e^{x_{\gamma}} \tilde u(t,x_{\gamma})
	\end{split}
\ee
and 
\be\label{e.c6294}
	\begin{split}
		e^{x_{\gamma}}\tilde u(t,x_{\gamma})
	= e^{- \frac{\beta}{2} \int_{x_\gamma}^\infty \tilde u dy} v(t,x_\gamma)
	\leq e^{C (t+1)^{3/2} e^{-(t+1)^\gamma}}
		\sqrt{t+1}\omega\Big(\log (t+1), \frac{x_\gamma}{\sqrt{t+1}}\Big),
	\end{split}
\ee
recalling that $x_\gamma=(t+1)^\gamma$.
Using~\eqref{e.c6293}-\eqref{e.c6294} first, and then recalling \Cref{lem:omega} yields
\be\label{aug1002}
\begin{aligned} 
	&|e^{x_{\gamma}}\tilde u(t,x_{\gamma}) -\alpha_{\infty}x_{\gamma}| 
		\\&\leq \Big|e^{x_{\gamma}}\tilde u(t,x_{\gamma}) - \sqrt{t+1} \omega\Big( \log(t+1), \frac{x_\gamma}{\sqrt{t+1}}\Big)\Big|
		+ \sqrt{t+1} \omega\Big( \log(t+1), \frac{x_\gamma}{\sqrt{t+1}}\Big)
			\Big|e^{\frac{x_\gamma^2}{4(t+1)}} - 1\Big|
		\\&\qquad\qquad
			 + \Big|\sqrt{t+1} \omega\Big( \log(t+1), \frac{x_\gamma}{\sqrt{t+1}}\Big)e^{\frac{x_\gamma^2}{4(t+1)}} - \alpha_\infty x_\gamma 
		 \Big|
		\\&\leq
			C\Big((t+1)^{2}e^{-(t+1)^\gamma} + \frac{x_\gamma^2}{4(t+1)}\Big)
			\omega\Big(\log (t+1), \frac{x_{\gamma}}{\sqrt{t+1}}\Big)
		\\&\qquad\qquad
		+ \Big(|h(\log(t+1))| + \frac{\sqrt{t+1}}{x_\gamma}|R(\tau,x_\gamma/\sqrt t)|\Big) x_{\gamma}e^{- \frac{x_\gamma^2}{6(t+1)}}.
\end{aligned}
\ee
Using that $h(\log(t+1)) \to 0$ and choosing $\gamma' \in(0,\gamma)$, we bound the second term on the very right side above by
\[
	\frac{\sqrt{t+1}}{x_\gamma}|R(\tau,x_\gamma/\sqrt t)|
		\leq C_{\gamma'} (t+1)^{-(\gamma-\gamma')}
		\to 0,\quad \text{ as }~t\to \infty.
\]
To handle the first term on the right side of (\ref{aug1002}), we simply notice that $\omega$ is bounded due to \Cref{lem:omega} and the time dependent terms in front of $\omega$ tend to zero.  This gives~\eqref{aug279}.
$\Box$

\subsection*{From the scales $x\sim O(t^\gamma)$ to $x\sim O(1)$}

The next step is to pass from the control of the
solution on the spatial scales $x\sim O(t^\gamma)$, provided by Lemma~\ref{lem:close},
to the spatial scales~$x\sim O(1)$ using the pulled nature of the Burgers-FKPP
equation for $\beta<2$.  

The arguments are a bit different for $\beta\le 0$ and $\beta\in(0,2)$.
We first discuss the latter case where the arguments deviate  from~\cite{NRR1}, as we have discussed in the outline of this section,
and later explain why the case~$\beta<0$ is quite similar to what was done in~\cite{NRR1}
for $\beta=0$, 
 
\subsubsection*{The case $0<\beta<2$}

We are going to use an argument inspired by~\cite{NRR1}
but we will only apply it after an
application of the weighted Hopf-Cole transform, and the conclusion is different as a result.
We take a traveling
wave $\phi_\beta(x)$, and shift it into the moving frame:
\be\label{jul2802}
\vphi_\alpha(t,x)=\phi_\beta(x+\zeta_\alpha(t)),
\ee
leaving the definition of $\zeta_\alpha(t)$ open for the moment.  
This function satisfies
\be\label{e.c6295}
	\begin{aligned}
		&\partial_t \varphi_\alpha - \partial_x^2\varphi_\alpha - \left(2 - \frac{3}{2(t+1)}\right) \partial_x \varphi_\alpha + \beta \varphi_\alpha \partial_x \varphi_\alpha - \varphi_\alpha + \varphi_\alpha^2
			= \left( \frac{3}{2(t+1)} + \dot \zeta_\alpha\right) \partial_x \varphi_\alpha.
	\end{aligned}
\ee
Next, we define its Hopf-Cole transform as
\be\label{e.c714}
\psi_\alpha(t,x) = e^{\Gamma_\alpha(t,x)}\vphi_\alpha(t,x),~~~
	\Gamma_\alpha(t,x) = x+\frac{\beta}{2}\int_x^{\infty} \vphi_\alpha(t,y) dy,
\ee
Noticing that
\be\label{e.c716}
	\partial_x \varphi_\alpha(t,x) = \phi_\beta'(x + \zeta_\alpha(t)),~~~
	\partial_t \varphi_\alpha(t,x) = \dot \zeta_\alpha(t) \phi_\beta'(x + \zeta_\alpha(t))
		= \dot \zeta_\alpha(t) \partial_x \varphi_\alpha(t,x),
\ee
we obtain an equation for the function $\psi$:
\be
	\begin{split}\label{e.c6283}
    &\partial_t \psi_\alpha - \partial_x^2\psi_\alpha + \frac{3}{2(t+1)}(\partial_x \psi_\alpha-\psi_\alpha) - \frac{\beta \psi_\alpha}{2}\int_x^{\infty} 
    \varphi_\alpha(1-\varphi_\alpha) dy +\psi_\alpha \varphi_\alpha\\
    & = \frac{\beta \dot{\zeta_\alpha}e^{\Gamma_\alpha} \varphi_\alpha }{2}\Big(\int_x^{\infty} \partial_x \varphi_\alpha dy\Big) 
    + e^{\Gamma_\alpha} \partial_x \varphi_\alpha\dot{\zeta_\alpha} 
    \\&~~-e^{\Gamma_\alpha} \Big(-\frac{\beta}{2} \partial_x \varphi_\alpha \varphi_\alpha
    + \left(1-\frac{\beta}{2}\varphi_\alpha\right)^2\varphi_\alpha + 2\left(1-\frac{\beta}{2}\varphi_\alpha\right) \partial_x \varphi_\alpha
    +\partial_x^2 \varphi_\alpha\Big)\\
    &~~+\frac{3e^{\Gamma_\alpha} }{2(t+1)}\Big(\partial_x \varphi_\alpha-\frac{\beta}{2}\varphi_\alpha^2\Big) - 
    \frac{\beta e^{\Gamma_\alpha} \varphi_\alpha }{2}\int_x^{\infty} \varphi_\alpha(1-\varphi_\alpha) dy 
    + e^{\Gamma_\alpha} \varphi_\alpha^2\\
    &= e^{\Gamma_\alpha} \Big(\dot{\zeta_\alpha}+\frac{3}{2(t+1)}\Big)\Big(-\frac{\beta}{2}\varphi_\alpha^2 + \partial_x \varphi_\alpha\Big) 
    + \frac{\beta e^{\Gamma_\alpha}\varphi_\alpha}{2} \Big(\partial_x \varphi_\alpha - \frac{\beta}{2}\varphi_\alpha^2
    + 2 \varphi_\alpha - \int_x^{\infty} 
    \varphi_\alpha(1-\varphi_\alpha) dy\Big)\\
    & = 
   - e^{\Gamma_\alpha} \Big(\dot{\zeta_\alpha}+\frac{3}{2(t+1)}\Big)\Big(2 \varphi_\alpha - \int_x^{\infty} 
    \varphi_\alpha(1-\varphi_\alpha) dy\Big).
	\end{split}
	\raisetag{7\normalbaselineskip}
\ee
We used (\ref{sep75}) in the last step above.  

For each $\alpha \in (\alpha_\infty/2, 2 \alpha_\infty)$, where $\alpha_\infty>0$ is defined in \Cref{lem:close}, we fix $\zeta_\alpha(t)$ by the
normalization 
\be\label{e.c6291}
	\psi_\alpha(t,(t+1)^\gamma)
	= \alpha (t+1)^\gamma.
\ee
We note that, differentiating~\eqref{e.c714} and using the asymptotics~\eqref{asymp1} of $\phi_\beta$, 
it is easy to check that~$\psi_\alpha$ is increasing in $x$ when $x\gg1$.  Hence, $\zeta_\alpha$ is well-defined for $t$ sufficiently large.  
In the sequel, we work with $\alpha = \alpha_\infty \pm \eps$, in which case, $\psi_\alpha$ approximately matches with $v$ at $x = (t+1)^\gamma$.  
Note a difference with~\cite{NRR1}:
the shift is determined by the value of $\psi_\alpha(t,x)$ at the point $x=(t+1)^\gamma$,
and not by the value of $\vphi_\alpha(t,x)$ as in~\cite{NRR1}. 

%

Next, using the asymptotics~\eqref{asymp1} of $\phi_\beta$, we find that
\be\label{e.zeta}
	\zeta_\alpha(t)
		= - \log(\alpha) - \frac{\log(\alpha)}{(t+1)^\gamma}
			+ o(t^{-\gamma}),
\ee
and 
\be\label{e.c6297}
|\dot \zeta_\alpha(t)|\leq \frac{C}{(t+1)^{1 + \gamma}},
\ee
where $C$ is independent of $\alpha$ over the interval $(\alpha_\infty/2,2\alpha_\infty)$.  
In addition, using~\eqref{asymp1} again, we find that
\be\label{e.c713}
e^{\Gamma_\alpha} \Big(\varphi_\alpha+\int_x^\infty \vphi_\alpha(1-\vphi_\alpha)dx\Big)
		\leq C (|x| + 1).
\ee
In view of~\eqref{e.c6283},~\eqref{e.c6297}, and~\eqref{e.c713}, we find, when $|x| \leq (t+1)^\gamma$,
\be\label{e.c6298}
	\left|\partial_t \psi_\alpha - \partial_x^2\psi_\alpha + \frac{3}{2(t+1)}(\partial_x \psi_\alpha-\psi_\alpha) - \frac{\beta \psi_\alpha}{2}\int_x^{\infty} 
    \varphi_\alpha(1-\varphi_\alpha) dy +\psi_\alpha \varphi_\alpha\right|
    	\leq \frac{C}{(t+1)^{1-\gamma}}.
\ee

Let us also recall that the Hopf-Cole transform $v(t,x)$ of the function $\tilde u$, defined in 
(\ref{dec332}), satisfies~(\ref{dec333})-(\ref{dec334}):
\begin{equation}\label{jul2804}
v_t - v_{xx}+\frac{3}{2(t+1)}(v_x-v) = 
-\Big(\tilde u(t,x)-\frac{\beta}{2}\int_x^{+\infty} \tilde u(t,y)(1-\tilde u(t,y))dy\Big)v.
\end{equation}
For $\alpha\neq \alpha_\infty$, we define $v_\alpha(t,x)$ as the solution to (\ref{jul2804}), thought of as a linear equation for $v$, with prescribed function $\tilde u(t,x)$:
\begin{equation}\label{jul2806}
\partial_tv_\alpha -\partial_x^2 v_{\alpha}+\frac{3}{2(t+1)}(\partial_x v_\alpha-v_\alpha) = 
-\Big(\tilde u(t,x)-\frac{\beta}{2}\int_x^{+\infty} \tilde u(t,y)(1-\tilde u(t,y))dy\Big)v_\alpha
\end{equation}
over $x < (t+1)^\gamma$ and $t> T_\eps$ with the boundary condition
\begin{equation}\label{jul2808}
v_\alpha(t,(t+1)^\gamma)=\alpha(t+1)^\gamma 
\end{equation}
and the initial condition
\be\label{jul2814}
v_\alpha(T_\eps,x)=v(T_\eps,x),
\ee
where $T_\eps$ is that of \Cref{lem:close} with $\eps = |\alpha - \alpha_\infty|$.

We now make a couple of observations.  
First, note that the normalization (\ref{e.c6291}) and
the boundary condition (\ref{jul2808}) imply that 
\be\label{jul2812}
v_\alpha(t,(t+1)^\gamma)=\psi_\alpha(t,(t+1)^\gamma),~~\hbox{ for all $t> T_\eps$ and $\eps>0$.}
\ee

Next, we claim that for any $\eps>0$, we have
\be\label{jul2816}
v_{\alpha_\infty + \eps}(t,x)
	>v(t,x)
	> v_{\alpha_\infty - \eps}(t,x),~~\hbox{ for all $t\ge T_\eps$ and $x<(t+1)^\gamma$.}
\ee
This follows from the fact that $v$ and $v_{\alpha_{\infty}\pm\eps}$ satisfy the same linear parabolic
equation (\ref{jul2804}) and~(\ref{jul2806}) and the same initial condition (\ref{jul2814}), but take ordered values at  the boundary.  
Indeed, with the help of Lemma~\ref{lem:close}, we obtain
\be\label{jul2818}
\bal
v(t,(t+1))^\gamma&=\exp\Big\{(t+1)^\gamma+\frac{\beta}{2}\int_{(t+1)^\gamma}^\infty\tilde u(t,y)dy\Big\}
\tilde u(t,(t+1)^\gamma)\\
&\le
(\alpha_\infty + \eps/2) (t+1)^\gamma\exp\Big\{ \frac{\beta}{2}\int_{(t+1)^\gamma}^\infty\tilde u(t,y)dy\Big\}\\
	&\le \Big(\alpha_\infty+\eps/2\Big) 
(t+1)^\gamma \exp\left\{C(t+1)^{3/2} e^{-(t+1)^\gamma}\right\}
	\le v_{\alpha_\infty+\eps}(t,(t+1)^\gamma),
\enbal
\ee
for $t>T_\eps$, up to increasing $T_\eps$ to deal with the different exponential factors in the second line above.  A similar argument gives $v_{\alpha_\infty - \eps}(t, (t+1)^\gamma) < v(t, (t+1)^\gamma)$.  Thus, comparison yields (\ref{jul2816}), as claimed.

Finally, we observe that
\be\label{e.c7291}
	v_{\alpha_\infty +\eps}(t,x)
		\leq \frac{\alpha_\infty+\eps}{\alpha_\infty-\eps} v_{\alpha_\infty-\eps}(t,x)
			\qquad\text{ for all } t \geq T_\eps \text{ and } x < (t+1)^\gamma.
\ee
The above is due to a comparison argument using the fact that the function in the right side of~(\ref{e.c7291}) 
solves the same linear equation as $v_{\alpha_\infty+\eps}$, with the same boundary condition, but with larger initial data.  The important consequence of~\eqref{e.c7291}, along with~\eqref{jul2816}, is that
\be\label{e.c7292}
	\begin{split}
	v_{\alpha_\infty + \eps}(t, - (t+1)^\gamma)
		&\leq \frac{\alpha_\infty+\eps}{\alpha_\infty-\eps} v_{\alpha_\infty-\eps}(t, - (t+1)^\gamma)
		\leq \frac{\alpha_\infty+\eps}{\alpha_\infty-\eps} v(t, - (t+1)^\gamma)
		\\&
		= \frac{\alpha_\infty+\eps}{\alpha_\infty-\eps} \exp\Big\{
			- (t+1)^\gamma + \frac{\beta}{2}\int_{-(t+1)^\gamma}^\infty \tilde u dy
		\Big\} \tilde u
		\leq C e^{-\frac{1}{2}(1-\beta/2)(t+1)^\gamma}.
	\end{split}
\ee
For the last inequality, we used~\eqref{calb}.

Let us define, for any $\eps \in(-1,1)$,
\be\label{jul2820}
s_\eps(t,x)=v_{\alpha_\infty+\eps}(t,x)-\psi_{\alpha_\infty+\eps}(t,x),
~~x<(t+1)^\gamma,~~t>T_\eps.
\ee
The following bound will be crucial for us.  
\begin{lemma}\label{lem:93}
There exists $\lambda>0$, so that we have, for any $\eps>0$:
\be\label{jul2821}
s_\eps(t,x)\le\farc{C}{(t+1)^\lambda}~~\hbox{ for all $t>T_\eps$ and $|x|\leq(t+1)^\gamma$.}
\ee
%
\end{lemma}
{\bf Proof.} 
We temporarily abuse notation and denote $s_\eps = s$ in this proof. Using \eqref{e.c6298} and~\eqref{jul2806}, we find, for $|x| \leq (t+1)^\gamma$,
\begin{equation}\label{sep76}
\begin{aligned}
s_t &- s_{xx}+ \frac{3}{2(t+1)}(s_x-s)+ s\Big(-\frac{\beta}{2}\int_x^{\infty}\tilde u(1-\tilde u) dy
+\tilde u\Big)\\
        &+ \psi\Big( \tilde u-\varphi+\frac{\beta}{2}\Big(\int_x^{\infty}\varphi(1-\varphi) dy-\int_x^{\infty}\tilde u(1-\tilde u)dy\Big)\Big)	\leq \frac{C}{(t+1)^{1-\gamma}}.
\end{aligned}
\end{equation}
Using~\eqref{e.c7292}, we see that
\be\label{jul1914}
	s(t, -(t+1)^\gamma)
		\leq C e^{- \frac{1}{2}(1 - \beta/2) (t+1)^\gamma},
\ee
and, by (\ref{jul2812}),
\be\label{jul2822}
s(t, (t+1)^\gamma)= 0.
\ee
The inequality (\ref{jul1914})  is another reminder of the importance of the condition $\beta < 2$ here.

Lemma~\ref{lem:neg} tells us that the zero order coefficient in (\ref{sep76}) is positive:
\be\label{jul2823}
G(t,x)=\tilde u(t,x)-\frac{\beta}{2}\int_x^{\infty}\tilde u(t,y)(1-\tilde u(t,y)) dy\ge 0.
\ee
In addition, we  claim that
\be\label{e.c717}
(\tilde u-\varphi)+\frac{\beta}{2}\Big(\int_x^{\infty}\varphi(1-\varphi) dy
-\int_x^{\infty}\tilde u (1 - \tilde u)dy\Big)\geq 0.
\ee
Indeed, recalling from~\eqref{e.c716} that $\varphi_x(t,x) = \phi_\beta'(x + \zeta(t))$,
and using the notation $E$ and $\bar E$ set in~(\ref{dec132}), 
we find
%
%
\be\label{21jun2120}
\begin{aligned} 
 &   \int_x^{\infty} \varphi(1-\varphi)dy - \int_x^{\infty} \tilde u(1 - \tilde u) dy 
    = \int_0^{\varphi(t,x)} z(1-z)\frac{dz}{|\bar{E}(z)|} - 
    \int_0^{\tilde u(t,x)}z(1-z)\frac{dz}{|E(t,z)|}\\
    &\geq \int_0^{\varphi(t,x)} z(1-z)\frac{dz}{|\bar{E}(z)|} - 
    \int_0^{\tilde u(t,x)}z(1-z)\frac{dz}{|\bar{E}(z)|}
 = -\int_{\varphi(t,x)}^{\tilde u(t,x)}z(1-z)\frac{dz}{|\bar{E}(z)|}.
\end{aligned}
\ee
As we show in (\ref{dec922}), we have
\be\label{21jun2121}
\frac{z(1-z)}{|\bar{E}(z)|}\leq \frac{2}{\beta},
\ee
and (\ref{e.c717}) follows from (\ref{21jun2120}) and (\ref{21jun2121}).

We deduce from (\ref{sep76}), (\ref{jul2823}), and (\ref{e.c717}) 
that $s(t,x)$ satisfies a differential inequality
\begin{equation}\label{sep761}
\begin{aligned}
s_t &- s_{xx}+ \frac{3}{2(t+1)}(s_x-s)+G(t,x)s\leq \frac{C}{(t+1)^{1-\gamma}},
\end{aligned}
\end{equation}
with $G(t,x)\ge 0$ and boundary conditions (\ref{jul1914})-(\ref{jul2822}). 
We can now argue as in \cite{NRR1} that we can find $\lambda,\gamma, \eps$ 
sufficiently small so that, for
\be\label{e.c723}
	\bar s(t,x)
		= \frac{1}{(t+1)^\lambda}\cos\Big(\frac{x}{(t+1)^{\gamma+\eps}}\Big),
\ee
we have
\be\label{e.c724}
	\bar s_t - \bar s_{xx} + \frac{3}{2(t+1)} (\bar s_x - \bar s)+G(t,x)\bar s
		\geq \frac{C}{(t+1)^{1-\gamma}}
\ee
for $t \geq T$, up to possibly increasing $T$.  
The choice of $\gamma$ occurs in this step.  
Hence, $\bar s$ is a super-solution of~\eqref{sep76} and, up to multiplying $\bar s$ by a large constant so that $s(T,\cdot) \leq \bar s(T,\cdot)$, we have
\[
	s(t,x) \leq \bar s(t,x)
		\qquad \text{ for all } t\geq T
			\text{ and } |x| \leq (t+1)^\gamma.
\] 
From this, (\ref{jul2821}) follows, finishing the proof of Lemma~\ref{lem:93}.~$\Box$ 

By virtue of \Cref{lem:93}, we have now established that, for $t > T_\eps$ and $|x| \leq (t+1)^\gamma$,
\be\label{jul2826}
	v
		\le v_{\alpha_\infty+\eps}
		= s_\eps + \psi_{\alpha_\infty + \eps}
		\le\psi_{\alpha_\infty+\eps}+\farc{C_\eps}{(t+1)^\lambda}.
\ee
The first inequality is due to~\eqref{jul2816} and the second is due to~\eqref{jul2821}.   We point out that we do not have a more precise lower bound.  This is another place where the steepness comparison will play a crucial role. 

Finally, note that the shifts corresponding to $\psi_{\alpha_\infty\pm\eps}$ satisfy, as in (\ref{e.zeta}):
\be\label{jul2835}
\zeta_{\alpha_{\infty}\pm\eps}(t)=-\log(\alpha_{\infty}\pm\eps)+O(t^{-\gamma}).
\ee

\subsubsection*{The end of the proof of Theorem~\ref{thm:main} for $\beta<2$}

 \subsubsection*{The case $0<\beta< 2$}
The definitions~\eqref{dec414} of~$\tilde u$ and~\eqref{jul2802} of the difference~$\varphi_{\alpha_\infty}$, and
the asymptotics~\eqref{e.zeta} of $\zeta$ 
imply that \Cref{thm:main} reduces to the uniform convergence of $\tilde u - \varphi_{\alpha_\infty}$ 
to zero on $\Rm$.  So far, we have only shown a weak version of closeness for their weighted Hopf-Cole transforms in (\ref{jul2826}). Indeed, notice that
\be
	\psi_{\alpha_\infty+\eps}(t,(t+1)^\gamma)
		- v(t, (t+1)^\gamma)
		= (\eps + o(1)) t^\gamma,
\ee
where $o(1)$ vanishes as $t\to\infty$.  It follows that the established inequalities are quite far apart due to the $t^\gamma$ factor.  We handle this issue now.

Before going to the proof, note that the uniform convergence to zero of $\tilde u - \varphi_{\alpha_\infty}$ reduces to showing that, for any $L>0$,
\be\label{e.c725}
	\lim_{t\to\infty} \sup_{x \in [-L, L]}	
	|\tilde u(t,x) - \varphi_{\alpha_\infty}(t,x)|
		=0.
\ee
This is sufficient due to the convergence in shape in Proposition~\ref{prop:aug262}. 

Fix any $\eps\in(0,\alpha_\infty/2)$.  We note that, due to~\eqref{jul2835}, we have
\be\label{e.c7293}
	|\vphi_{\alpha_\infty}(t,x) - \vphi_{\alpha_\infty \pm \eps}(t,x)|
		\leq C \eps
			\qquad\text{ for all $t$ sufficiently large and any $x$}.
\ee
Hence, it suffices to establish upper bounds on $\tilde u - \vphi_{\alpha_\infty+\eps}$ and $\vphi_{\alpha_\infty - \eps}-\tilde u$.

We establish the latter first,  as it is much simpler.  Indeed, note that
\[
\begin{split}
\varphi_{\alpha_\infty - \eps}(t, (t+1)^\gamma)
&= \exp\Big\{- (t+1)^\gamma - \frac{\beta}{2} \int_{(t+1)^\gamma}^\infty \vphi_{\alpha_\infty - \eps}(t,y) dy\Big\} 
\psi_{\alpha_\infty-\eps}(t, (t+1)^\gamma)
		\\&\leq (\alpha_\infty-\eps) (t+1)^\gamma \exp\left(- (t+1)^\gamma
		\right).
	\end{split}
\]
Thus, due to \Cref{lem:close}, we have $\varphi_{\alpha_\infty - \eps}(t,(t+1)^\gamma) < \tilde u(t,(t+1)^\gamma)$.  
Since $\tilde u$ is steeper than $\phi_\beta$ (recall \Cref{prop-jul1502}), we get that
\be\label{e.c7299}
	\varphi_{\alpha_\infty - \eps}(t,x) < \tilde u(t,x),
~~~\text{for all $t$ sufficiently large and $x < (t+1)^\gamma$},
\ee
whence, for $t$ sufficiently large,
\be\label{e.c72901}
	\sup_{x\in[-L,L]}
		\left(\vphi_{\alpha_\infty - \eps}(t,x)-\tilde u(t,x)\right)
			\leq 0.
\ee

Next, we consider the much more involved upper bound on $\tilde u - \vphi_{\alpha_\infty + \eps}$.  
Fix any $x \in [-L,L]$ and $t$ sufficiently large and write 
 \be\label{e.c7294}
 \begin{split}
 \tilde u(t,x) - &\varphi_{\alpha_\infty+\eps}(t,x)
= \exp\Big\{- x - \frac{\beta}{2} \int_x^\infty \tilde u dy\Big\} v(t,x)
			- \exp\Big\{- x - \frac{\beta}{2} \int_x^\infty \vphi_{\alpha_\infty+\eps} dy\Big\} \psi_{\alpha_\infty +\eps}(t,x)
		\\&= \exp\Big\{- x - \frac{\beta}{2} \int_x^\infty \vphi_{\alpha_\infty+\eps} dy\Big\} (v(t,x)- \psi_{\alpha_\infty +\eps}(t,x))
		\\&\qquad
			+ e^{-x} v(t,x)
				\Big[\exp\Big\{- \frac{\beta}{2} \int_x^\infty \tilde u  dy\Big\} - \exp\Big\{- \frac{\beta}{2} \int_x^{(t+1)^\gamma} \tilde u  dy\Big\}\Big]
		\\&\qquad 
			+ e^{-x} v(t,x)
				\Big[\exp\Big\{- \frac{\beta}{2} \int_x^{(t+1)^\gamma} \tilde u  dy\Big\} 
				- \exp\Big\{- \frac{\beta}{2} \int_x^{(t+1)^\gamma} \vphi_{\alpha_\infty+\eps} dy\Big\}\Big]
		\\&\qquad
			+ e^{-x} v(t,x)
				\Big[\exp\Big\{- \frac{\beta}{2} \int_x^{(t+1)^\gamma} \varphi_{\alpha_\infty +\eps}  dy\Big\}
					- \exp\Big\{- \frac{\beta}{2} \int_x^\infty  \vphi_{\alpha_\infty+\eps}  dy\Big\}\Big]\\
	&			=: I_1 + I_2 + I_3 + I_4.
	\end{split}
 \ee
As $\beta>0$, we immediately see that
\be\label{e.c7311}
	I_2 \leq 0.
\ee
Hence, we need only bound $I_1$, $I_3$, and $I_4$.

To handle the term $I_1$ in~\eqref{e.c7294}, we apply~\eqref{jul2826} to find, since $\beta\in(0,2)$:
\be\label{aug1006}
 	\begin{split}
 		I_1
		\leq \frac{C_\eps}{(t+1)^\lambda}
			\exp\Big\{- x - \frac{\beta}{2} \int_x^\infty \tilde u dy\Big\}\le \frac{C_\eps}{(t+1)^\lambda}e^L.
	\end{split}
 \ee
%
The term $I_4$ is the next simplest.  Using~\eqref{calb}, we see that
\be\label{aug1004}
	e^{-x} v(t,x)
		\leq C (t+1)^{3/2}.
\ee
Using (\ref{aug1004}) and also~\eqref{asymp1} to handle the tail integral over $((t+1)^\gamma, \infty)$, we obtain 
\be\label{e.c7296}
\begin{split}
I_4&\leq C (t+1)^{3/2} \Big[\exp\Big\{- \frac{\beta}{2} \int_x^\infty \tilde u  dy\Big\} -
 \exp\Big\{- \frac{\beta}{2} \int_x^{(t+1)^\gamma} \vphi_{\alpha_\infty+\eps}  dy\Big\}\Big]
		\\&= C (t+1)^{3/2} \exp\Big\{- \frac{\beta}{2} \int_x^{(t+1)^\gamma} \varphi_{\alpha_\infty+\eps}  dy\Big\} 
		\Big[1- \exp\Big\{-\frac{\beta}{2} \int_{(t+1)^\gamma}^\infty \vphi_{\alpha_\infty+\eps}  dy\Big\}\Big]
		\\&\leq C (t+1)^{3/2} \Big[1- \exp\Big\{-\frac{\beta}{2} \int_{(t+1)^\gamma}^\infty \vphi_{\alpha_\infty+\eps}  dy\Big\}\Big]
		\leq C (t+1)^{3/2} e^{-(t+1)^\gamma}.
	\end{split}
\ee

We now handle $I_3$.  If the bracketed term in the definition of $I_3$ 
is non-positive, there is nothing to prove as $I_3 \leq 0$.  If the bracketed term is positive, we write
\[
	\begin{split}
	I_3
		&= e^{-x} [v(t,x) - \psi_{\alpha_\infty+\eps} (t,x)]
				\Big[\exp\Big\{- \frac{\beta}{2} \int_x^{(t+1)^\gamma} \tilde u  dy\Big\} -
				\exp\Big\{- \frac{\beta}{2} \int_x^{(t+1)^\gamma} \vphi_{\alpha_\infty+\eps} dy\Big\}\Big]
		\\& \quad + e^{-x} \psi_{\alpha_\infty+\eps}(t,x)
				\Big[\exp\Big\{- \frac{\beta}{2} \int_x^{(t+1)^\gamma} \tilde u  dy\Big\} - 
				\exp\Big\{- \frac{\beta}{2} \int_x^{(t+1)^\gamma} \vphi_{\alpha_\infty+\eps} dy\Big\}\Big] := I_{31} + I_{32}.
	\end{split}
\]
The reason for the extra step here is that we do not know {\em a priori} that $e^{-x} v$ is bounded; however, we do know that $e^{-x} \psi_{\alpha_\infty + \eps}$ is.

For $I_{31}$, we use~\eqref{jul2826}, the positivity of the bracketed term, and the fact that the bracketed term is smaller than one (again recall that $\beta \in (0,2)$) to find
\[
	\begin{split}
	I_{31}
		\leq \frac{C_\eps}{(t+1)^\lambda} e^{-x}
				\Big[\exp\Big\{- \frac{\beta}{2} \int_x^{(t+1)^\gamma} \tilde u  dy\Big\} - \exp\Big\{- \frac{\beta}{2} \int_x^{(t+1)^\gamma} \vphi_{\alpha_\infty+\eps} dy\Big\}\Big]
		\leq e^{L} \frac{C_\eps}{(t+1)^\lambda}.
	\end{split}
\]
For $I_{32}$, we first notice that
\[
e^{-x} \psi_{\alpha_\infty + \eps}(t,x)
\leq C e^L.
\]
In addition, from~\eqref{e.c7299}, we see that
\[
\exp\Big\{- \frac{\beta}{2} \int_x^{(t+1)^\gamma} \tilde u dy\Big\}
		< \exp\Big\{- \frac{\beta}{2} \int_x^{(t+1)^\gamma} \varphi_{\alpha_\infty-\eps} dy\Big\}.
\]
Thus, we find
\[
	I_{32}
		\leq C e^L \Big[
				\exp\Big\{- \frac{\beta}{2} \int_x^{(t+1)^\gamma} \varphi_{\alpha_\infty-\eps} dy\Big\}
				- \exp\Big\{- \frac{\beta}{2} \int_x^{(t+1)^\gamma} \varphi_{\alpha_\infty+\eps} dy\Big\}
			 \Big].
\]
Taylor expanding the exponential and using the asymptotics~\eqref{jul2835} of $\zeta_{\alpha_\infty\pm \eps}$  and~\eqref{asymp1} of $\phi_\beta$, we see that
\[
	I_{32}
		\leq C_L \int_x^{(t+1)^\gamma} (\varphi_{\alpha_\infty+\eps} - \varphi_{\alpha_\infty-\eps}) dy
		\leq C_L \eps.
\]
We conclude that
\be\label{e.c7297}
	I_3
		\leq \frac{C_\eps e^L}{(t+1)^\lambda} + C_L \eps.
\ee

Combining~\eqref{e.c7294},~\eqref{e.c7311},~\eqref{aug1006},~\eqref{e.c7296}, and~\eqref{e.c7297}, and increasing $t$ suitably, we obtain 
\be\label{e.c7298}
	\sup_{x\in[-L,L]}
		\left(\tilde u(t,x) - \vphi_{\alpha_\infty+\eps}\right)
		\leq C_L \eps.
\ee
The claim~\eqref{e.c725} then follows from~\eqref{e.c7293},~\eqref{e.c72901}, and~\eqref{e.c7298}.  This concludes the proof for~$\beta\in(0,2)$.

\subsubsection*{The case $\beta\le 0$}

This case is essentially the same, except for some simplifications,
so we only highlight the changes that need to be made.  The key estimate  is to obtain an upper bound on
\[
	s_\eps(t,x) = e^x ( \tilde u(t,x) - \vphi_{\alpha_\infty+\eps}).
\]
Now, $\varphi_\alpha$ is defined such that
\[
	\varphi_\alpha(t,(t+1)^\gamma) = \alpha (t+1)^\gamma e^{-(t+1)^\gamma - t^{2\gamma-1}/4}.
\]
In this case $\zeta_\alpha$ has the same asymptotics as in~\eqref{jul2835}.  Then $s_\eps$ satisfies 
\begin{equation}\label{sep77}
\begin{aligned}
	s_t - s_{xx} &+ \frac{3}{2(t+1)}(s_x - s)
		+ \beta \tilde u (s_x - s) + s(\tilde u + \varphi_{\alpha_\infty+\eps})
		\\&= - \beta s \partial_x\varphi_{\alpha_\infty+\eps} - \partial_x\varphi_{\alpha_\infty+\eps} \left( \dot \zeta_{\alpha_\infty+\eps} + \frac{3}{2(t+1)}\right)
		\leq \frac{C}{(t+1)^{1-\gamma}}.
    \end{aligned}
\end{equation}
The last inequality holds when $|x| \leq (t+1)^\gamma$ for the same reasons as in~\eqref{e.c6298} and due to the fact 
that $\beta \leq 0$ and $\partial_x \varphi_{\alpha_\infty+\eps} = \phi_\beta' \leq 0$.   Thus, the same upper bound as in \Cref{lem:93} holds.

On the other hand, it is easy to check that $\tilde u > \varphi_{\alpha_\infty-\eps}$ on $x< (t+1)^\gamma$  for $t$ sufficiently large.  This is
due to \Cref{lem:close}, which yields the correct ordering at $x = (t+1)^\gamma$, and the fact that $\tilde u$ is steeper than $\varphi_{\alpha_\infty - \eps}$.

The combination of the above with a simpler version of the argument for the case $\beta \in(0,2)$ yields the desired convergence.  This finishes the proof of \Cref{thm:main} when $\beta < 2$.~$\Box$

\section{The case $\beta=2$: bounds on the front location}\label{sec:twp-bounds}

We now turn to the case $\beta=2$, where the analysis is particularly delicate.  
Let us define the shift~$\mu(t)$ by
\be\label{21apr1342}
u(t,2t-\mu(t))=\farc{1}{2}.
\ee 
Note that, due to the sign convention in (\ref{21apr1342}), 
an upper bound on $\mu(t)$ is a lower bound on the front location
and vice versa. 
Our goal in this section is to prove the following  upper and lower bounds on
$\mu(t)$. 
\begin{prop}\label{prop-apr1302}
Let $u(t,x)$ be the solution to the Burgers-FKPP equation (\ref{burgerskpp}),
with the initial condition as in Theorem~\ref{thm:main}, and let $\mu(t)$ be defined by
(\ref{21apr1342}).
There exists $m_0>0$ that depends on~$u(0,x)$  so that
\be\label{21apr1343}
\farc{1}{2}\log t-m_0\le \mu(t)\le \farc{1}{2}\log t+m_0,~~\hbox{ for all $t\ge 1$.}
\ee
\end{prop}
Later, in Section~\ref{sec:proof-beta=2}, 
we will improve these bounds to the precise asymptotics of \Cref{thm:main}:
\be\label{may2502}
\mu(t)=\farc{1}{2}\log t+x_0+o(1)\hbox{ as $t\to+\infty$.}
\ee
Proposition~\ref{prop-apr1302} is, however, a crucial step in the proof of (\ref{may2502}). 
Its proof occupies almost all of the rest of this section and is the heart of this paper. In order to explain
the outline of the proof, we will need to do some preliminary transformations, leading 
to Lemmas~\ref{lem-may2502} and~\ref{lem-may2504} below that imply the conclusion
of Proposition~\ref{prop-apr1302}.  We will discuss their
proofs when we come to their respective statements.
As a technical comment, we mention that without loss of generality we will, once again, take~$u(0,x)=\one(x\le 0)$. 

At the very end of this section, we use Proposition~\ref{prop-apr1302} to obtain a helpful bound in an intermediate
region in a short Section~\ref{sec:middle-end}.

\subsection{Outline of the proof of Proposition~\ref{prop-apr1302}}
 \label{sec:beta2_outline}

\subsubsection{An exponential moment}\label{sec:exp-moment} 


We first give a heuristic argument to explain the delay in the case $\beta=2$ and the role of $\beta$. 
Consider the Burgers-FKPP equation, in the moving frame $x\to x-2t+\mu(t)$, with
an unknown shift~$\mu(t)$: 
\begin{align}\label{jan1102}
 \tilde u_t - 2\tilde u_x+ \mu'(t) \tilde u_x+ \beta \tilde u \tilde u_x = \tilde u_{xx}+\tilde u-\tilde u^2.
\end{align}
To highlight the role of $\beta$, we have not yet specified it to the value $\beta =2$. 
A simple computation, using only~\eqref{jan1102} and several integration by parts, shows that the exponential moment
\begin{equation}\label{feb404}
I(t) = \int \tilde u(t,x) e^x dx
\end{equation}
satisfies an ODE
\begin{equation}\label{jan1300}
	\begin{split}
	\frac{dI(t)}{dt}
		&= \int (2\tilde u_x-\mu'(t) \tilde u_x- \beta \tilde u \tilde u_x +\tilde u_{xx}+\tilde u-\tilde u^2)e^x dx\\
		&= \int (-2\tilde u + \mu'(t) \tilde u - \frac{\beta}{2} (\tilde u^2)_x - \tilde u_x+\tilde u-\tilde u^2)e^x dx\\
		&= \int (-2\tilde u + \mu'(t) \tilde u + \frac{\beta}{2} \tilde u^2 + \tilde u+\tilde u-\tilde u^2)e^x dx\\
		&= \int \Big(\Big(\frac{\beta}{2} - 1\Big) \tilde u^2 + \mu'(t) \tilde u\Big)e^x dx
		= \mu'(t)\int \tilde u e^x dx = \mu'(t) I(t).
	\end{split}
\end{equation}
In the first equality, we simply used~\eqref{jan1102}. In the second-to-last inequality, we used that $\beta =2$.   In every inequality between, we only integrated by parts and cancelled terms.  
Therefore, 
\begin{align}\label{jan1301}
 I(t) = I(0)e^{\mu(t)},
\end{align}
as long as we choose $\mu(0)=0$. As is clear from~\eqref{jan1300}, this algebraic property is specific to $\beta=2$.

If $\mu(t)$ is the ``correct frame," 
in the sense that~$\tilde u(t,x)$ converges to a traveling wave, we expect from the philosophy of the self-similar variables that  $\tilde u(t,x)$
has the asymptotics 
\begin{align}\label{feb1204}
\tilde u(t,x) \sim \frac{1}{1+e^{x-x_0}} e^{-{x^2}/{(4t)}},~~\hbox{ for $t\gg1$ and $x \geq 0$.}
\end{align}
This indicates that if $I(t)$ is computed in the correct reference frame, then 
\begin{equation}\label{feb1202}
I(t)=\int e^x \tilde u(t,x)dx\approx \int_0^\infty \frac{e^x}{1+e^{x-x_0}} e^{-{x^2}/{(4t)}}dx \sim C_0\sqrt{t},~~\hbox{ for $t\gg 1$,}
\end{equation}
with an explicit constant $C_0$ that depends on $x_0$. 
For this to be consistent with (\ref{jan1301}), we should have 
\begin{equation}\label{feb408}
\mu(t)\sim \farc{1}{2}\log t+C_1,~~~~\hbox{ for $t\gg1$,}
\end{equation}
explaining the $(1/2)\log t$ shift for the front position. 
The rest of this section and the following one is 
a justification of (\ref{feb408}). 

\subsubsection{An inhomogeneous conservation law}

In order to explain the outline of the proof of Proposition~\ref{prop-apr1302}, we need
to introduce a change of variable related to the evolution of mass for the exponential moment
in (\ref{jan1300}). 
Let $u(t,x)$ be the solution to the Burgers-FKPP equation (\ref{burgerskpp})
with the specific value~$\beta =2$, and set 
\be\label{21jun2002}
\hat u(t,x) = u(t, x+ 2t),
\ee
which
satisfies
\begin{align}\label{apr1322}
 \hat u_t - 2\hat u_x+2 \hat u \hat u_x = \hat u_{xx}+\hat u(1-\hat u),
\end{align}
Note that we use a slightly different notation here for the shifted function $\hat u$ rather than~$\tilde u$.  
This is to denote the difference in shift: $\hat u$ is shifted into the moving frame $2t$, while $\tilde u$ is shifted to the moving frame $2t - \mu(t)$ matching the front, as,
for instance, in the proof of \Cref{prop:aug262}. 

The first key observation is that the function
\be\label{apr1327}
p(t,x)=e^x\hat u(t,x),
\ee
satisfies a viscous spatially inhomogeneous conservation law:
\begin{equation}\label{21apr1102}
p_t+\big(p^2e^{-x}\big)_x=p_{xx},
\end{equation}
with the initial condition
\be\label{21apr1335}
p(0,x)=e^x\one(x\le 0).
\ee
Recall that we assume, without loss of generality, that $u(0,x)=\one(x\le 0)$. 
Notice that~\eqref{21apr1102} conserves mass for $p(t,x)$:
\be\label{may2516}
\int p(t,x)dx=\int p(0,x)dx=1.
\ee
The normalization (\ref{21apr1342}) 
in terms of $\hat u(t,x)$  becomes
\be\label{21may2506}
\hat u(t,-\mu(t))=\farc{1}{2},
\ee
which translates into 
\be\label{21may2508}
p(t,-\mu(t))=\farc{1}{2}e^{-\mu(t)}.
\ee

A simple preliminary observation is that
the solution to~(\ref{21apr1102})-(\ref{21apr1335}) satisfies
\be\label{apr1337}
p(t,x)\le e^x,~~\hbox{ for all $t\ge 0$ and $x\in\Rm$,}
\ee
simply because $\hat u(t,x)\le 1$. Note also that $p(x)=e^x$ and
\be\label{apr1230}
p(x)=\farc{e^x}{1+e^{(x-\xi)}}
\ee
are exact solutions to (\ref{21apr1102}), for any shift $\xi\in\Rm$. 
 
The proof of Proposition~\ref{prop-apr1302} relies on the analysis of 
the solution to~(\ref{21apr1102})-(\ref{21apr1335}) and
proceeds in the following steps. First, we prove an upper bound on~$\mu(t)$.
\begin{lemma}\label{lem-may2502}
There exists $K_1>0$ so that $\mu(t)$ be defined by (\ref{21may2506}) satisfies
\be\label{21may2514}
\mu(t)\le\farc{1}{2}\log t+K_1,~~\hbox{ for all $t\ge 1$.} 
\ee
\end{lemma}
Lemma~\ref{lem-may2502} is proved in Section~\ref{s:upper}. Let us give a brief outline of the proof.
It  is based on the conservation (\ref{may2516})
of the total mass of $p(t,x)$, together with estimating the mass of $p(t,x)$ separately
in the regions $\{x<-\mu(t)\}$, $\{-\mu(t)<x<N\sqrt{t}\}$ and~$\{x>N\sqrt{t}\}$. 

First, we will show that the mass of $p(t,x)$ in the region $\{x>N\sqrt{t}\}$ is exponentially small in~$N$ 
by an argument that bounds
exponential moments of $p$.  Next, using the simple exponential bound~\eqref{apr1337}, we can check that 
the mass in the region $\{x < -\mu(t)\}$ is bounded above by~$\exp\{-\mu(t)\}$.  

To bound the mass in the middle region, recall that when $\beta=2$,  the traveling wave moving with the minimal speed $c_*=2$ is explicit
and is given by~(\ref{jul1502}):
\be \label{21apr1332}
\phi(x)=\farc{1}{1+e^x}.
\ee
From this and the steepness comparison in Proposition \ref{prop-jul1502}, we immediately deduce the following useful property. 
\begin{lemma}\label{lem:jun21}
If $u(0,\cdot)$ is steeper than $\phi$ given by (\ref{21apr1332}), then, for each $v \in (0,1)$ and $t>0$, we have
\be\bal
&u(t,x)\geq\ \Big(1 + \frac{1-v}{v} e^{x-x(t,v)}\Big)^{-1},~~\text{if } x < x(t,v),\\
&u(t,x)\leq\ \Big(1 + \frac{1-v}{v} e^{x-x(t,v)}\Big)^{-1},~~\text{if } x > x(t,v). 
\enbal
\ee
Here, $x(t,v)$ is defined implicitly by $v = u(t,x(t,v))$.
\end{lemma}
By \Cref{lem:jun21},
and (\ref{21may2508}), we know that
\be\label{21apr1330}
	p(t,x)
		= e^x \hat u(t,x)
		\le \frac{e^x}{1+e^{x+\mu(t)}}
		\leq e^{-\mu(t)},
			~~~\hbox{ for $x>-\mu(t)$}.
\ee
It follows that the mass in the region $-\mu(t) \leq x \leq N\sqrt t$ is bounded by $2N \sqrt t e^{-\mu(t)}$ 
(it is easy to show the weak bound $\mu(t) < N \sqrt t$, 
see~\eqref{21apr1340} below).  Combining all three bounds and recalling mass conservation~\eqref{may2516} of $p$, we deduce that $\mu(t)$ 
must satisfy the upper bound in \Cref{lem-may2502}.  The details are given in \Cref{s:upper}.

To prove a lower bound for $\mu(t)$ we use the following lemma.
\begin{lemma}\label{lem-may2504}
There exist $C>0$ so that, for all $x\in \R$ and $t>0$,
\be\label{21may2520}
	p(t,x)\le\farc{C}{\sqrt{t}}.
\ee
The constant $C$ depends on the initial data nontrivially.
\end{lemma}
We prove \Cref{lem-may2504} in \Cref{s:lower}, and its surprisingly delicate
proof is outlined there.  
When combined with (\ref{21may2508}), Lemma~\ref{lem-may2504} implies the lower bound
in Proposition~\ref{prop-apr1302}
\be\label{21may2521}
\mu(t)\ge \farc{1}{2}\log t-C,\hbox{ for all $t\ge 1$.}
\ee

Let us make a brief comment that the $t^{-1/2}$ decay rate in (\ref{21may2520}) 
is standard for parabolic equations in one dimension, and would be expected
for a solution to (\ref{apr1327}). Nevertheless, 
the proof of Lemma~\ref{lem-may2504} is much less straightforward than one would naively expect.
To illustrate the potential obstacles, notice that~\eqref{21apr1102} can be written as
\be\label{e.c643}
	p_t + (\hat up)_x = p_{xx}.
\ee
One might hope to ``forget'' the connection between $\hat u$ and $p$, 
and prove $t^{-1/2}$ decay for linear divergence form advection-diffusion equations
of the form (\ref{e.c643}) in general,
with an advection term~$\hat u(t,x)$ that, say, 
connects two constants on the left and on the right. This seems to
be a good cartoon for $\hat u(t,x)$. 
However, such decay cannot hold in general.  Indeed, consider the following explicit example.  Let
\[
q(x) = \frac{1}{e^x+e^{-x}} ,
 \quad\text{ and }\quad
v(x) = -\frac{e^x-e^{-x}}{e^x+e^{-x}},
\]
so that $q(x)$ is a steady solution to 
\[
q_t + (vq)_x = q_{xx}.
\]
The function $q(x)$ is rapidly decaying at infinity  but, of course, it
does not decay in time.  This shows that the boundedness or existence of the limits 
at infinity of $\hat u$ are not sufficient to determine the decay.
It is crucial that~$\hat u$ is negligible on $\R_+$ and its profile does not move
much -- these are, however, exactly the properties that we are trying to prove. Circumventing 
these difficulties requires interesting a priori estimates for our specific problem,
which will be discussed in Section~\ref{s:lower}.

To summarize, \Cref{prop-apr1302} reduces to Lemmas~\ref{lem-may2502} and~\ref{lem-may2504} and we prove these lemmas in the rest of this section.

\subsection{Preliminary weak bounds on $\mu$}  

Let us first give some very poor bounds on $\mu(t)$ that, at least, ensure that it does not behave too wildly. These are useful in the sequel.

For a very simple bound, a comparison to the standard KPP equation
with $\beta=0$ that uses monotonicity of $u(t,x)$ in $x$, implies that
\be\label{21apr1340}
\mu(t)\le\farc{3}{2}\log t+C,
\ee
with a universal constant $C$.  

We now use the conservation of mass for $p$~\eqref{may2516} to obtain a lower bound on $\mu$. 
Recalling \Cref{lem:jun21} and the definition~\eqref{21apr1342} of $\mu$, we find
\[
p(t,x)
	= e^x \hat u(t,x)
	\geq \farc{e^x}{1+e^{x+\mu(t)}},~~\hbox{ for $x<-\mu(t)$.}
\]
We emphasize that the proof of \Cref{lem:jun21} is independent of all other lemmas in this section -- it is simply a consequence of the steepness comparison.  Therefore, we have
\[
\bal
1&=\int_{-\infty}^{\infty} p(0,x)dx=\int_{-\infty}^\infty p(t,x)dx\ge \int_{-\infty}^{-\mu(t)} \farc{e^x}{1+e^{x+\mu(t)}}dx=e^{-\mu(t)}
\int_{-\infty}^0\farc{e^x}{1+e^x}dx\ge \farc{1}{2}e^{-\mu(t)}.
\enbal
\]
We conclude that 
\be\label{21apr1436}
\mu(t)\ge -\log 2.
\ee
The bounds (\ref{21apr1340}) and (\ref{21apr1436}) will be greatly improved below.

\subsection{An upper bound on the shift: the proof of Lemma~\ref{lem-may2502}}\label{s:upper}

As we have mentioned, the strategy of the proof is to show that if $\mu(t)$ is too large then
$p(t,x)$ can not have a mass larger than $1/10$ in any of the three regions 
\be\label{21may2522}
L=\{x<-\mu(t)\},~~M=\{-\mu(t)\le x\le N\sqrt{t}\},~~R=\{x>N\sqrt{t}\},
\ee
provided that $N$ is also chosen sufficiently large. This would contradict (\ref{may2516}).
Note that (\ref{21apr1436})  implies   $N\sqrt{t}\ge -\mu(t)$ for $t\ge 1$ and
$N$ sufficiently large, so that the regions above are well-defined.

For the left region $L$, we simply apply \eqref{apr1337} and write 
\be\label{21may2523}
	\int_L p(t,x) dx
		= \int_{-\infty}^{-\mu(t)}p(t,x)dx
		\le \int_{-\infty}^{-\mu(t)}e^x dx
		= e^{-\mu(t)}.
\ee
For the middle region $M$, we also have a simple estimate that uses \Cref{lem:jun21}:
\be\label{21may2524}
	\begin{split}
		\int_M p(t,x) dx
			&= \int_{-\mu(t)}^{N\sqrt{t}} e^x \hat u(t,x)dx
			\le \int_{-\mu(t)}^{N\sqrt{t}} \farc{e^x}{1+e^{x+\mu(t)}} dx\\
			&=e^{-\mu(t)} \int_0^{N\sqrt{t}+\mu(t)} \farc{e^x}{1+e^{x}}dx
			\le CNe^{-\mu(t)}\sqrt{t},
	\end{split}
\ee
for all $t\ge 1$. We used (\ref{21apr1340}) in the last step above. 

Now, we deal with the right region $R$.  
First we state a bound on an exponential moment of $p$: fix any $m\in(0,1/2)$ and let
\be\label{21apr1406}
I_m(t)=\int \big(e^{mx}+e^{-mx}\big)p(t,x)dx. 
\ee
We claim that
\be\label{e.c6152}
	I_m(t)
		\leq C e^{C m^2t}.
\ee
We postpone the proof of~\eqref{e.c6152} momentarily and show how to conclude the bound for
the integral over $R$.  For $N>1$ large, we estimate this mass by using (\ref{21apr1406})
with $m=1/\sqrt{t}$. This gives
\be\label{apr1133}
\bal
	\int_R p(t,x) dx
	&=\int_{N\sqrt t}^\infty p(t,x)dx
	\le \Big(\int _{N\sqrt{t}}^\infty e^{x/\sqrt{t}}p(t,x)dx\Big)^{1/2}\Big(\int _{N\sqrt{t}}^\infty e^{-x/\sqrt{t}}p(t,x)dx\Big)^{1/2}\\
	&\le C\Big(\int _{N\sqrt{t}}^\infty e^{-x/\sqrt{t}}p(t,x)dx\Big)
	\le  C \Big(\int _{N\sqrt{t}}^\infty e^{-N} p(t,x)dx\Big)^{1/2}\le C e^{-N/2}.
\enbal
\ee
We used the conservation of mass (\ref{may2516}) in the last step.

We may now put the estimates (\ref{21may2523}), (\ref{21may2524}) and (\ref{apr1133}) together to obtain
\be\label{21may2526}
\bal
	1
		&=\int_{-\infty}^\infty p(0,x)dx
		=\int_{-\infty}^\infty p(t,x)dx
		=\int_L p(t,x)dx+\int_M p(t,x)dx+\int_R p(t,x)dx\\
&\le e^{-\mu(t)}+CN\sqrt{t}e^{-\mu(t)}+C e^{-N/2}.
\enbal
\ee  
Choosing $N = 2 \log(2C)$ gives
\be\label{21may2527}
\bal
\farc{1}{2}& \le e^{-\mu(t)}+C\sqrt{t}e^{-\mu(t)},
\enbal
\ee  
so that
\be\label{21may2528}
e^{\mu(t)}\le C(1+\sqrt{t}),
\ee
and (\ref{21may2514}) follows, as desired.

All that remains to finish the proof of Lemma~\ref{lem-may2502} is to establish~\eqref{e.c6152}, which we do now.   Recalling $I_m$ from~\eqref{21apr1406}, multiplying~\eqref{21apr1102} by $e^{mx} + e^{-mx}$, and integrating yields
\be\label{21apr1106}
\bal
\frac{d}{dt} I_m
	=m^2\int (e^{mx}+e^{-mx})p(t,x)dx+m\int e^{-x}[e^{mx}-e^{-mx}]p^2(t,x)dx.
\enbal
\ee 
We decompose the last integral:
\be\label{21apr1312}
\bal
& J_m=m\int e^{-x}[e^{mx}-e^{-mx}]p^2(t,x)dx=m\int e^{-x}\tanh(mx)[e^{mx}+e^{-mx}]p^2(t,x)dx\\
	&\leq m \int_{-\infty}^{\log(2)} e^{-x} |\tanh(mx)| [e^{mx}+e^{-mx}]p^2(t,x) dx\\
		&\quad+ m \int^{\infty}_{\log(2)} e^{-x} |\tanh(mx)| [e^{mx}+e^{-mx}]p^2(t,x) dx
	=J_m^-+J_m^+.
\enbal
\ee 
For $J_m^-$, we use (\ref{apr1337}), the bound 
\be\label{jun1602}
|\tanh(mx)| \leq m|x|,
\ee
and the assumption $0<m <1/2$ to obtain
\be\label{21apr1316}
\bal
J_m^- 
	&\leq m \int_{-\infty}^{\log(2)} |\tanh(mx)| [e^{mx} + e^{-mx}] e^{x} dx
	\leq m^2 \int_{-\infty}^{\log(2)} |x| e^{(1-m)x} dx
	\leq C m^2.
\enbal
\ee 
On the other hand, as $x > \log(2) > -\mu(t)$, by~\eqref{21apr1436}, we may apply \Cref{lem:jun21} to obtain
\[
\bal
J_m^+
	&= m \int_{\log(2)}^{+\infty} |\tanh(mx)| [e^{mx} + e^{-mx}] \hat u(t,x) p(t,x) dx \\
	&\leq m \int_{\log(2)}^{+\infty} |\tanh(mx)| [e^{mx} + e^{-mx}] \frac{p(t,x)}{1 + e^{x+\mu(t)}} dx\\
	&\leq m^2 \int_{\log(2)}^{+\infty} \frac{|x|}{1 + e^{x+\mu(t)}} [e^{mx} + e^{-mx}] p(t,x) dx
	\leq m^2 e^{-\mu(t)} \int_{\log(2)}^{+\infty} [e^{mx} + e^{-mx}] p(t,x) dx.
\enbal
\]
In the second inequality, we used (\ref{jun1602}), 
and in the last inequality, we used that
\[
	|x| \leq e^x \leq e^{-\mu(t)}(1 + e^{x+\mu(t)}).
\]
Recalling that $-\mu(t) \leq \log(2)$ from~\eqref{21apr1436}, we find
\be\label{21apr1318}
	J_m^+
		\leq C m^2 I_m.
\ee 
The combination of~\eqref{21apr1316} and~\eqref{21apr1318} yields
\[
	\frac{dI_m}{dt}\le C_0m^2(I_m+1),
\]
for some universal constant $C_0$.  This implies that
\[
I_m(t)\le e^{C_0m^2t}[I_m(0)+1]\le 5e^{C_0m^2t},
\]
since, recalling that $m< 1/2$,
\[
I_m(0)=\int_{-\infty}^0 e^{x}(e^{mx}+e^{-mx})dx\le 4.
\]
This concludes the proof of~\eqref{e.c6152} and, thus, that of Lemma~\ref{lem-may2502}.~$\Box$

\subsection{A lower bound on the shift: the proof of Lemma~\ref{lem-may2504}}\label{s:lower}

The proof of Lemma~\ref{lem-may2504} is quite a bit more involved than that for Lemma~\ref{lem-may2502}.  
Let us first explain the main steps of the proof. Recall that the standard $L^\infty$-decay for diffusion equations
of the self-adjoint form 
\be\label{21jun1614}
z_t=\nabla\cdot(a(x)\nabla z),~~x\in\Rm^n,
\ee
with a uniformly positive and bounded diffusivity $a(x)$ is obtained as follows: first, one gets the dissipation inequality
\be\label{21jun1604}
\farc 12\farc{d}{dt}\int |z(t,x)|^2dx\le -C\int|\nabla z(t,x)|^2dz.
\ee
An application of the Nash inequality leads, after solving an elementary differential inequality,
to the $L^1-L^2$ decay estimate
\be\label{21jun1606}
\|z(t,\cdot)\|_{L^2}\le \farc{C}{t^{n/4}}\|z(0,\cdot)\|_{L^1}.
\ee
The self-adjoint form of (\ref{21jun1614}) and the estimate (\ref{21jun1606}) give the dual bound
\be\label{21jun1608}
\|z(t,\cdot)\|_{L^\infty}\le \farc{C}{t^{n/4}}\|z(0,\cdot)\|_{L^2}.
\ee
The last step is to apply the semi-group property and the above estimates to deduce that
\be\label{21jun1616}
\|z(t,\cdot)\|_{L^\infty}\le \farc{C}{t^{n/4}}\|z(t/2,\cdot)\|_{L^2}\le 
\farc{C}{t^{n/2}}\|z(0,\cdot)\|_{L^1}. 
\ee
See, e.g., \cite[Section~2.4]{Davies} for a full treatment of this.
 
It seems not possible to 
directly obtain a dissipation inequality for the $L^2$-norm of the function~$p(t,x)$,
starting with (\ref{21apr1102}), due to the spatial inhomogeneity  of the nonlinear term. 
Instead, to get an analogue of~(\ref{21jun1604}), we will use a 
suitably chosen weight~$\rho(t,x)$ that weighs~$\R_+$ more than~$\R_-$, and establish 
an $L_w^2$-dissipation inequality for the function 
\be\label{21jun1612}
\varphi=\farc{p}{\rho}.
\ee  
Here we use $w$ as a subscript to emphasize that we are working in weighted Lebesgue spaces.  Such   weight 
allows to focus on where advection is negligible and diffusion dominates the evolution 
of~\eqref{e.c643}. 

It turns out that the weight $\rho(t,x)=1-\hat u(t,x)$ 
can actually be used to produce a dissipation inequality because, as we will see,
the function $1-\hat u$ is a super-solution for (\ref{21apr1102}) satisfied by~$p(t,x)$. 
This property is not purely algebraic: it will, once again, use the steepness comparison of $\hat u(t,x)$ to the traveling wave. 
The dissipation computation is directly inspired by the relative entropy
arguments for linear advection-diffusion equations in~\cite{Const,MPP}. Here, however,
we compute the entropy relative not to a solution but a super-solution, and the nonlinear
nature of the present situation requires specific cancellations. This is the subject of \Cref{lem-may2506} below.

We also establish a weighted Nash inequality stated in Proposition~\ref{lem-may2508} below.
When adapted to our setting, it yields the appropriate long time decay of the $L_w^2$-norm of $\varphi$, 
up to a (potentially large) boundary layer in time.   An interesting complication is that the weighted Nash inequality
holds for a nontrivial class of weights satisfying certain assumptions, and we need to control the fact that our weight,
coming from the solution to a nonlinear evolution equation, satisfies these assumptions
for all $t>0$ in a uniform way.  This is done in the course of the proof of Lemma~\ref{lem-may3002}
below.

In contrast to the standard proof for diffusion equations, we cannot directly pass from $L^2_w$-decay to $L^\infty$-decay.  
Indeed, due to the weight, the aforementioned $L_w^2$-decay estimate
is {\em not} an~$L_w^1\to L_w^2$ estimate as the boundary layer depends 
on the initial $L^2_w$-norm of $\varphi$ and, hence, the usual adjointness trick in~(\ref{21jun1608}),
used to establish 
the $L^\infty$ decay is not available.  This complication is present even for the linear equation~\eqref{e.c643} when the relationship between $p$ and $\hat u$ is ``forgotten.''  
In fact, decay in~$L^\infty_w$ of $\varphi$ 
is not even expected in the setting in which we find ourselves. 
To pass from the $L_w^2$-bounds for the function $\vphi(t,x)$ to the $L^\infty$ decay for 
$p(t,x)$, we use time averages to find a particular intermediate 
time $T_g< T$, at which $\varphi$ satisfies ``good'' pointwise bounds in a region of interest.  
Those bounds, stated in Lemma~\ref{l.T_g}, can be transferred 
to show that $p$ is bounded both by $e^x$ and~$C/\sqrt T$.  We can then ``trap'' that estimate 
going forward in time, from the time $T_g$ until the time $T$, 
by breaking $p$ up into a small mass part, which necessarily stays small 
due to conservation of mass and parabolic regularity theory, and a part that sits under an explicit, 
small super-solution. This will conclude the proof.

As the reader will surely have noticed, there is a subtle, but serious, issue in the above outline. 
The $L^2_w$ decay requires that the initial $L^2_w$ norm be bounded, 
which, as can be immediately seen from (\ref{21jun1612}),
 is not true when $\hat u(0,\cdot)$ is the Heaviside function, 
the function $p(0,\cdot)$ is given by (\ref{21apr1335}), and $\rho=1-\hat u$.   
This requires an extra step where we choose an approximate initial condition
\be\label{21jun114bis}
\hat u_a(0,x)=\farc{1}{1+e^{\gamma(x-a)}}\one(x\le 0),
\ee
with $\gamma\in(1,2)$, to obtain new solutions $\hat u_a$ and $p_a$. Here, $a>0$ 
is a parameter depending on the final time $T$ at which we wish to establish the upper bound.  
A careful analysis shows that, with an appropriate choice of $a$,
the two solutions, $p$ and $p_a$, stay $O(1/\sqrt T)$ away from each other,
due to the error estimate in Lemma~\ref{lem-may2806}.

Below, we  state the upper bound on the modification $p_a$ and show how to bootstrap that 
bound to the decay of $p$ itself.  This is 
done in \Cref{ss:bootstrap}.  Then we give the proof of the upper bound on $p_a$ following the outline above, in \Cref{ss:decay}.

\subsubsection{The proof of \Cref{lem-may2504}} \label{ss:bootstrap}

\subsubsection*{A modified initial condition}

We first construct the modified solution $p_a(t,x)$. 
Recall that $\hat u(t,x)$ is the solution to the Burgers-FKPP equation
\be\label{e.21jun112}
	\hat u_t-2\hat u_x+2\hat u\hat u_x=\hat u_{xx}+\hat u-\hat u^2,
\ee
with the initial condition $\hat u_{\rm in}(x)=\one(x\le 0)$. For any $a>0$ and $\gamma \in (1,2)$, let
$\hat u_a(t,x)$ be the solution to (\ref{e.21jun112}) with the initial condition
\be\label{21jun114}
	\hat u_a(0,x)=\farc{1}{1+e^{\gamma(x-a)}}\one(x\le 0),
\ee
and set
\be\label{21jun1702}
p(t,x) = e^x \hat u(t,x),~~p_a(t,x) = e^x \hat u_a(t,x).
\ee
As in~\eqref{21apr1102} and~\eqref{e.c643}, the function $p_a$ satisfies
\be\label{e.c644}
	\partial_t p_a
		+ \partial_x (p_a^2 e^{-x})
			= \partial_x^2 p_a
\ee
and
\be\label{e.c645}
	\partial_t p_a
		+ \partial_x(\hat u_a p_a)
			= \partial_x^2 p_a.
\ee
We have made the switch to $\partial_{t,x}$ notation to avoid the awkward double sub-script.  
It is  easy to observe that $\hat u_a\leq 1$, and we have, by the comparison principle (notice that $e^x$ solves~\eqref{e.c644}):
\be\label{e.c6158}
	p_a(t,x)
		\leq e^x
			\qquad\text{ for all } t\geq 0,~~x \in \R.
\ee

Two important quantities of interest for us are  
\be\label{e.I_M}
	\begin{split}
&M_a(t)
		= \int \left( \frac{p_a(t,x)}{1- \hat u_a(t,x)}\right) (1 - \hat u_a(t,x)) \dx
		= \int p_a(t,x) \dx
		= M_a(0),\\
		&I_a(t)
		= \int \left( \frac{p_a(t,x)}{1- \hat u_a(t,x)}\right)^2 (1 - \hat u_a(t,x)) \dx
		= \int \frac{p_a(t,x)^2}{1- \hat u_a(t,x)}\dx.
\end{split}
\ee
The last equality in the first line above follows by conservation of mass.  Note that $M_a$
and $I_a$ are, respectively, the  weighted $L^1$ and $L^2$ norms of the function 
$\vphi=p_a/(1-\hat u_a)$ with weight $1-\hat u_a$. 
Two easy computations show that $M_a(0)$ is uniformly bounded in~$a>0$:
\be\label{e.Ma le 1}
	\log(2)
		= \int_{-\infty}^0 \frac{e^x}{1 + e^x} \dx 
		\leq M_a
		\leq \int_{-\infty}^0 e^x \dx
		= 1
\ee
and $I_a(0)$ is finite for all $a>0$:
\be\label{21may2820}
\bal
	I_a(0)
		&=\int\frac{e^{2x}\hat u_a^2(0,x)}{1-\hat u_a(0,x)}\dx
		\le \int_{-\infty}^0 \frac{e^{2x} }{1-\hat u_a(0,x)}\dx
		=\int_{-\infty}^0 e^{2x}\farc{1+e^{\gamma(x-a)}}{e^{\gamma(x-a)}}\dx\\
		&\leq 2 e^{\gamma a} \int_{-\infty}^0 e^{(2-\gamma)x}\dx
		= \farc{2}{2-\gamma}e^{\gamma a}
\enbal
\ee
but is not uniformly bounded as $a\to+\infty$. Note that an analogous estimate for $I(0)$ 
with $u_a$ and~$p_a$ replaced by $u$ and $p$, respectively, does not hold, as 
\be\label{21jun1706}
\farc{p^2(0,x)}{1-u(0,x)} = +\infty \cdot\one(x<0).
\ee
This is what prevents us from establishing the decay of $p$ directly.

As we see in the sequel, another crucial feature of this altered initial condition is that 
$\hat u_a(0,\cdot)$ is steeper than the traveling wave $1/ (1 + e^x)$. 
This is why we take $\gamma>1$ in~(\ref{21jun114}). 
The restriction~$\gamma<2$
comes from the upper bound on $I_a(0)$ in (\ref{21may2820}). 

The main estimate we establish on $p_a$ is:
\begin{lemma}\label{l:altered u_0 decay}
	Let $\hat u_a(t,x)$ be the solution to (\ref{e.21jun112}) with the initial condition (\ref{21jun114}) for some $a>0$ and~$1<\gamma<2$, 
and $p_a(t,x)=e^x\hat u_a(t,x)$. 
There exists a universal constant $C_{\ref{l:altered u_0 decay}}>0$ such that the following holds.  Set
\be\label{21may3004}
	t_1(a) =  C_{\ref{l:altered u_0 decay}} I_a(0),
\ee
with $I_a$ defined in~\eqref{e.I_M}. 
There exists $K>0$ that does not depend on~$a>1$ so that  
\be\label{21may3010}
p_a(t,x)\le\farc{K}{\sqrt{t}},~~\hbox{ for all } t\ge 4t_1(a).
\ee
\end{lemma}
We postpone the proof of \Cref{l:altered u_0 decay} until Section~\ref{ss:decay}.  
First, we obtain a closeness estimate on $p$ and $p_a$.  Its proof is succinct enough to give it immediately.
\begin{lemma}\label{lem-may2806}
Fix $a\in \R$ and $\gamma > 0$.  Let $\hat u(t,x)$ and $\hat u_a(t,x)$ be the solutions 
to (\ref{e.21jun112}) with the respective initial conditions~$\hat u(0,x)=\one(x\le 0)$ 
and~\eqref{21jun114}, respectively. 
There is a constant~$C>0$ that does not depend on $a$ or $\gamma$ such that~$p(t,x)=e^x \hat u(t,x)$ 
and $p_a(t,x)=e^x \hat u_a(t,x)$ satisfy
\be\label{21may2710}
		0\le p(t,x)-p_a(t,x)\le Ce^{-\gamma a}\hbox{ for all $x\in\Rm$ and $t\ge 1$.}
\ee
\end{lemma}
{\bf Proof.} 
Let $h(t,x) = p(t,x) - p_a(t,x)$ for all $t\geq 0$ and $x \in \R$.  Note that $h(0,\cdot) \geq 0$ 
and, due to~\eqref{21apr1102} and~\eqref{e.c644}, that $h$ satisfies the parabolic equation
\be
h_t - h_{xx}= - (p^2 e^{-x})_x + (p_a^2 e^{-x})_x
		= - (h w)_x,
\ee
with
\[
	w(t,x)
		=e^{-x}(p_a(t,x)+p(t,x))
		=\hat u(t,x)+\hat u_a(t,x).
\]
Hence, the comparison principle implies that $h\geq 0$.  This finishes the lower bound in~\eqref{21may2710}.

To conclude the upper bound, we use 
mass conservation for $p(t,x)$ and $p_a(t,x)$.  Indeed,
\be\label{21may2814}
\begin{split}
	\int h(t,x) \dx&= \int (p(t,x) - p_a(t,x)) \dx
		= \int (p(0,x) - p_a(0,x)) \dx
		= \int_{-\infty}^0 e^x  \frac{e^{\gamma(x-a)}}{1 + e^{\gamma(x-a)}}  \dx\\
		&= e^a \int_{-\infty}^{-a} \frac{e^{(1+\gamma)x}}{1 + e^{\gamma x}} \dx
		\le \farc{1}{1+\gamma}e^{a} e^{-(1+\gamma)a}
		=\farc{1}{1+\gamma}e^{-\gamma a}. 
\end{split}
\ee
Using parabolic regularity theory together with (\ref{21may2814}) and positivity of $h$, 
we get immediately that
\be\label{e.c5141}
	\sup_x h(t,x) \le Ce^{-\gamma a}
		\qquad\hbox{for all $t\ge 1$ and $x\in\Rm$,}
\ee
finishing the proof. The constant $C$ does not depend on $a$ or $\gamma$.~$\Box$

{\bf Proof of Lemma~\ref{lem-may2504}.} 
We are now in a position to combine \Cref{l:altered u_0 decay} with \Cref{lem-may2806}, 
in order to prove \Cref{lem-may2504}.  One delicate point is that Lemma~\ref{lem-may2806} 
requires us to take $a$ large.  On the other hand,
$I_a(0)$ in (\ref{21may2820}) blows up as $a\to+\infty$. This will require a careful balancing act.

We note that we need only prove \Cref{lem-may2504} for $T$ sufficiently large, as the claim follows for ``smaller'' $T$ by simply increasing the constant $C$.

Let $\gamma = 3/2$, and let $C_{\ref{l:altered u_0 decay}}$ be the universal constant from \Cref{l:altered u_0 decay}.  Given any 
\be\label{e.c110901}
	T>(16C_{\ref{l:altered u_0 decay}})^2+1,
\ee
choose
\be\label{e.c6153}
	a = \frac{1}{2\gamma}\log (T)
		> 0.
\ee
Recall from (\ref{21may2820}) and the choice of $\gamma$ that 
\be\label{21jun104}
\bal
I_a(0)&\le \farc{2}{2-\gamma}e^{\gamma a}
	= 4 e^{\gamma a}.
\enbal
\ee
By our choice of $a$ and~\eqref{21may3004}, it follows that
\be\label{e.21jun108}
	4t_1(a)
		\, {\color{red}
		= 4 C_{\ref{l:altered u_0 decay}} I_a(0)
		\leq 16 C_{\ref{l:altered u_0 decay}} e^{\gamma a}
		= 16 C_{\ref{l:altered u_0 decay}} e^{\gamma \left(\frac{1}{2\gamma} \log T\right)}
		= 16 C_{\ref{l:altered u_0 decay}} \sqrt T
		}
		\leq T,
\ee
where the last inequality follows due to~\eqref{e.c110901}. 
Due to~\eqref{e.21jun108}, we may apply \Cref{l:altered u_0 decay} to find
\be\label{e.c6154}
	\sup_{x}\, p_a(T,x)
		\leq \frac{K}{\sqrt T}.
\ee
We may also use Lemma~\ref{lem-may2806} to see that
\be\label{e.21jun114}
0\le p(t,x)-p_a(t,x)\le Ce^{-\gamma a},~~~~\hbox{ for all $x\in\Rm$ and $t>1$.}
\ee
Recalling the choice~\eqref{e.c6153} of $\gamma$ and $a$ in~\eqref{e.21jun114} and using~\eqref{e.c6154} to $t=T$, we deduce that
\[
	\sup_x \, p(T,x)
		\leq \sup_x \, p_a(T,x) + C e^{-\gamma a}
		\leq \frac{K}{\sqrt T} + \frac{C}{\sqrt T},
\]
finishing the proof of Lemma~\ref{lem-may2504}.~$\Box$

\subsubsection{The proof of \Cref{l:altered u_0 decay}} \label{ss:decay}

We now prove Lemma~\ref{l:altered u_0 decay},
following the outline from the beginning of \Cref{s:lower}.

\subsubsection*{Step one: decay in a weighted $L^2$ space}

Our first goal is to obtain an $L^2$-decay estimate, an analogue of (\ref{21jun1606}) in the
present situation. 
We will work with norms weighted by
\be\label{e.c641}
\rho_a = 1 - \hat u_a.
\ee
and make a change of function
\be\label{e.c642}
	\varphi_a = \frac{p_a}{1-\hat u_a} = \frac{p_a}{\rho_a}.
\ee
The main goal of this step is the decay rate of $I_a$, which we state here. 
Recall that $\hat u_a$ and $p_a$ are defined in \Cref{l:altered u_0 decay}.
\begin{lemma}\label{lem-may3002}
For $a>0$ and $\gamma \in (1,2)$, let  $t_1(a)$, $M_a$, and $I_a(0)$ be as in~\eqref{21may3004} and \eqref{21may2820}.  Then 
\be\label{21may3006}
	I_a(t_1)
		= \int\vphi_a^2(t_1,x)\rho_a(t_1,x) \dx\le 1,
\ee
and 
\be\label{21may3008}
	I_a(t)
		= \int \vphi_a^2(t,x)\rho_a(t,x) dx
		\leq \farc{C}{1+\sqrt{t-t_1(a)}},~~\hbox{ for $t\ge t_1(a)$.}
\ee
\end{lemma}

Note that $I_a$ may be thought of as $\|\vphi_a\|^2_{L^2(\rho_a)}$. To prove \Cref{lem-may3002}, we 
obtain a more general result for weighted $L^2$ decay, identifying assumptions on the weight under which 
decay holds.  Afterwards, we show that $\rho_a$ satisfies these assumptions and apply the general result
to our setting.

The first step is a dissipation inequality, inspired by the relative entropy arguments 
for linear equations in~\cite{Const,MPP}. As in~\cite{Const}, given a weight
$\rho(t,x)$ and an advection $v(t,x)$, we consider an
operator~$\cD_\rho$ defined by
\be\label{21jun1712}
\cD_\rho  q= \rho^{-2} \partial_x(\rho^2 \partial_x q) - v \partial_x q.
\ee
\begin{proposition}\label{lem-may2506}
Let $v(t,x)$ be a smooth bounded function, $q(t,x)$ be a solution to 
\be\label{21may2540}
q_t+(vq)_x=q_{xx},
\ee
and $\rho(t,x)$ be a super-solution to (\ref{21may2540}):
\be\label{21may2541}
\rho_t+(v\rho)_x\ge\rho_{xx},
\ee
with initial conditions $p(0,x)\ge 0$ and $\rho(0,x)>0$ such that 
\be\label{21may2542}
	\int \varphi(0,x)^2 \rho(0,x) dx
		<+\infty,
\ee
where $\vphi(t,x)=q(t,x)/\rho(t,x)$ for all $t\geq 0$ and $x\in \R$.  Then $\varphi$ satisfies
\be\label{e.c697}
	\vphi_t - \cD_\rho \vphi = - \vphi^2 \left( \frac{\rho_t - \rho_{xx} + (v\rho)_x}{\rho}\right),
\ee
and
\be\label{21may2544}
\frac{d}{dt}\int \ \vphi^2(t,x)\rho(t,x)dx\le -2\int\vphi_x^2(t,x)\rho(t,x)dx. 
\ee
\end{proposition}
The main difference with~\cite{Const,MPP} is that the function $\rho(t,x)$ is
not a solution to the advection-diffusion equation (\ref{21may2540}) but a super-solution. 

The second general result is an adaptation of the Nash inequality for weighted spaces.  This allows us to make use of the dissipation inequality in~\Cref{lem-may2506}.
\begin{proposition}\label{lem-may2508}
Let an increasing function $r \in L^1(\R_-) \cap C(\R)$ satisfy the following assumptions:
\be\label{21may2543}
\bal
&(i)~~~\lim_{x\to-\infty}r(x)=0,~~0<\lim_{x\to+\infty}r(x)=r_+<+\infty,
\\
&(ii)~~~ \doverline{r}(x)\le C_1\max\{1, \overline r(x)^2\} r(x),~~\hbox{ for all $x\in\Rm$},
\enbal
\ee
with
\be\label{21may2544bis}
\bar r(x)=\int_{-\infty}^x r(y)dy,~~\doverline r(x)=\int_{-\infty}^x \bar r(y)dy.
\ee
Then, for any $\theta>0$ and  any smooth non-negative function $\vphi(x)$ that is sufficiently rapidly decaying as $x\to+\infty$ and bounded as $x\to-\infty$ we have
\be\label{e.Nashish}
	\int \vphi^2(x)r(x)\, dx
		\le \frac{2}{\theta} \Big( \int \vphi(x)r(x) dx\Big)^2
			+ 8C_1\max\{1,\theta^2\} \int |\vphi_x(x)|^2 r(x) dx.
\ee
\end{proposition}
Let us point out that the weight 
\[
r(x)=1-\phi_2(x)=\farc{e^x}{1+e^x}
\]
satisfies assumptions of Proposition~\ref{lem-may2508}, and this is an example the reader may want to keep in mind.

Let us also briefly note the connection with the standard Nash inequality.  A key step in the proof 
of the latter is to establish~\eqref{e.Nashish} with $\max\{1,\theta^2\}$ replaced by $\theta^2$.   
The change here reflects the fact that the measure induced by $r$ is finite on $\R_-$.

The proofs of \Cref{lem-may2506} and \Cref{lem-may2508} are postponed until \Cref{s:abstract L2}. 
We first apply these results to establish the $L^2$ decay of $p_a$, that is, \Cref{lem-may3002}.
\subsubsection*{Proof of \Cref{lem-may3002}}

We first assume that the assumptions of the dissipation and Nash inequalities, 
that is, \Cref{lem-may2506} and \Cref{lem-may2508}, hold in a uniform way for $\rho_a(t,x)$, 
and show how to conclude.  Afterwards, we show how to verify those assumptions.

Assuming that we can apply \Cref{lem-may2506} and \Cref{lem-may2508}, we proceed as follows.  By \Cref{lem-may2506}, 
we have
\be\label{21may2842}
	\frac{d}{dt}\int \vphi_a^2(t,x)\rho_a(t,x)dx
		\le -2\int (\partial_x\vphi_a)^2(t,x) \rho_a(t,x)dx. 
\ee
Moreover, recalling the definition of $M_a$ from~\eqref{e.I_M} and that it is uniformly
bounded, as in~\eqref{e.Ma le 1}, and applying \Cref{lem-may2508}, we find, for any $\theta >0$,
\be\label{21may2841}
	\int \vphi_a^2(t,x)\rho_a(x)\, dx
		\le \frac{C}{\theta}
		+ C\max\{1,\theta^2\} \int |\partial_x\vphi_a(t,x)|^2 \rho_a(t,x) dx,
\ee
with a constant $C$ that does not depend on $t$, $\gamma$, or $a$.  

There are two cases to consider: first, assume that at some $t>0$ we have 
\be\label{21may2843}
	\int |\partial_x\vphi_a(t,x)|^2 \rho_a(t,x) \dx
		\geq 1.
\ee
Then, we get from (\ref{21may2842}): 
\be\label{21may2848}
	\farc{d}{dt}\int \vphi_a^2(t,x) \rho_a(t,x) \dx
		\leq-2.
\ee
On the other hand, if (\ref{21may2843}) fails, so that
\[
	\int |\partial_x\vphi_a(t,x)|^2 \rho_a(t,x) \dx
		< 1,
\]
we can use (\ref{21may2841}) with the constant 
\[
	\theta
		= 
			\Big( \int |\partial_x\vphi_a|^2 \rho_a\, dx\Big)^{-1/3}\geq 1.
\]
Keeping in mind the upper bound \eqref{e.Ma le 1} on $M_a$, this leads to
\[
	\int \vphi_a^2(t,x)\rho_a(x)\, dx
		\le  C \Big(\int |\partial_x\vphi_a(t,x)|^2 \rho_a(t,x) dx\Big)^{1/3}. 
\]
Therefore, in the second case we have 
\be\label{21may2847}
	\farc{d}{dt}\int \vphi_a^2(t,x) \rho_a(t,x) \dx
		\le - \farc{1}{C}\Big(\int \vphi_a^2(t,x) \rho_a(t,x) \dx \Big)^3.
\ee
Putting (\ref{21may2848}) and (\ref{21may2847})  together, we find
\be\label{may1710}
	\farc{d}{dt}\int \vphi_a^2(t,x) \rho_a(t,x) \dx
		\le - \frac{1}{C}\min\Big\{1, \Big(\int \vphi_a^2(t,x) \rho_a(t,x) \dx \Big)^3\Big\}.
\ee
Setting  
\be\label{may1714}
t_1(a)= C I_a(0)\geq C(I_a(0)-1)_+,
\ee
we see from (\ref{may1710}) that 
\be\label{21may2850}
\int\vphi_a^2(t_1,x)\rho_a(t_1,x) \dx\le 1. 
\ee
Therefore, for $t\ge t_1(a)$, we have
\be\label{21may2851}
	\farc{d}{dt}\int \vphi_a^2(t,x) \rho_a(t,x) \dx
		\le - C\Big(\int \vphi_a^2(t,x) \rho_a(t,x) \dx \Big)^3,
\ee
and
\be\label{may1702}
	\int \vphi_a(t,x)^2 \rho_a(t,x) dx
		\leq \farc{C}{1+\sqrt{t-t_1}}.
\ee
Hence, all that remains in the proof of Lemma~\ref{lem-may3002} is to verify 
the assumptions of \Cref{lem-may2506} and \Cref{lem-may2508}.

We consider first the assumptions of \Cref{lem-may2506}.  In the notation of \Cref{lem-may2506}, 
we have~$q = p_a$, $v = \hat u$, and $\rho = \rho_a$.  In view of~\eqref{e.c645} and~\eqref{21may2820}, 
the assumptions~\eqref{21may2540} and~\eqref{21may2542} are satisfied.  Hence,  the only assumption to 
check is~\eqref{21may2541}, a somewhat miraculous property that
\be\label{e.c646}
	\partial_t \rho_a  + \partial_x(\hat u_a \rho_a)
		\geq \partial_x^2 \rho_a.
\ee
First, using that $\hat u_a$ satisfies~\eqref{e.21jun112}, we find
\be\label{21jun1714}
	\begin{split}
		\partial_t \rho_a + \partial_x (\hat u_a \rho_a) - \partial_x^2 \rho_a
			&= - \partial_t \hat u_a + \partial_x\hat u_a - 2\hat u_a \partial_x\hat u_a + \partial_x^2 \hat u_a
			= - \partial_x\hat u_a - \hat u_a(1-\hat u_a).
	\end{split}
\ee
The initial condition for $\hat u_a$ in (\ref{21jun114}) was chosen to be steeper than
the  traveling wave $\phi$ given by~\eqref{21apr1332} -- this is
why we needed to take $\gamma>1$ in (\ref{21jun114}). 
Using Proposition~\ref{prop-jul1502}, we deduce that~$\hat u_a(t,x)$ 
is steeper than $\phi$ for all $t>0$.
Applying this property to (\ref{21jun1714}) gives 
\be\label{21jun1716}
\partial_t\rho_a + \partial_x(\hat u_a \rho_a) - \partial_x^2 \rho_a
= - \partial_x\hat u_a - \hat u_a(1-\hat u_a)
\geq - \phi_x(\phi^{-1}(\hat u_a))
- \phi(\phi^{-1}(\hat u_a)) (1 - \phi(\phi^{-1}(\hat u_a))).
\ee
Using the explicit expression \eqref{21apr1332} for $\phi$, 
it is straightforward to check that the right side of (\ref{21jun1716}) is non-negative everywhere:
\be\label{21may2821bis}
\begin{split}
-\phi_x - \phi (1-\phi)
&= - \left( - \frac{e^{x}}{(1+e^x)^2}\right)
				- \frac{e^x}{1+e^x} \frac{1}{1 + e^x}
			= 0.
	\end{split}
\ee
This can, of course, be also obtained from (\ref{sep75}).
Hence,~\eqref{e.c646} is established.

Next, we verify that $\rho_a$ satisfies assumptions (i) and (ii) of \Cref{lem-may2508} 
uniformly for all~$t>0$.    Assumption (i) holds automatically
since $\rho_a$ is increasing in $x$ (recall that $\hat u_a$ is decreasing in $x$), with $\rho_a(t,-\infty)=0$ and $\rho_a(t,+\infty)=1$ for all $t>0$. 
Assumption (ii) in  (\ref{21may2543})  requires that
\be\label{21may2823}
\bal
&\doverline{\rho}_a(t,x)\le C_1\max\{1, \overline \rho_a(t,x)^2\} \rho_a(t,x),~~\hbox{ for all $x\in\Rm$}.
\enbal
\ee
It is useful to (implicitly) define the analogue $\mu_a$ of $\mu$ for $\hat u_a$:
\be\label{e.mu_a}
	\hat u_a (t, -\mu_a(t)) = \frac{1}{2}.
\ee
Consider first the case when $x< -\mu_a(t)$.
As $\hat u_a(0,\cdot)$ is steeper than the traveling wave, so is $\hat u_a(t,\cdot)$,
by Proposition~\ref{prop-jul1502}, and we know from \Cref{lem:jun21} that, 
if $y \leq x$, then 
\be\label{21may2830}
	\rho_a(t,y)
		= 1- \hat u_a(t,y)
		\leq\farc{\rho_a(t,x)e^{y-x}}{1-\rho_a(t,x)+\rho_a(t,x)e^{y-x}}.
\ee
This gives  
\be\label{21may2831}
\bal
\bar\rho_a(t,x)&=\int_{-\infty}^x\rho_a(t,y)\dy\le \int_{-\infty}^x\farc{\rho_a(t,x)e^{y-x}}{1-\rho_a(t,x)+\rho_a(t,x)e^{y-x}}\dy\\
	&\le \farc{\rho_a(t,x)}{1-\rho_a(t,x)}\int_{-\infty}^x {e^{y-x}}\dy
	= \farc{\rho_a(t,x)}{1-\rho_a(t,x)}
	\le 2\rho_a(t,x). 
\enbal
\ee
In the last line we used that $x< -\mu_a(t)$, so that
\[
1-\rho_a(t,x) = \hat u_a(t,x) \geq 1/2.
\]
This implies (\ref{21may2823}) for $x$ such that $\rho_a(t,x)<1/2$:
\be\label{e.c5124}
\doverline \rho_a(t,x)= \int_{-\infty}^x \overline \rho_a(t,y) \dy
\le 2 \int_{-\infty}^x \rho_a(t,y) \dy
= 2\overline \rho_a(t,x)\le 4\rho_a(t,x).
\ee

On the other hand, when $x \geq -\mu_a(t)$, we have $\rho_a(t,x) \geq 1/2$, so that
\be\label{21may2832}
\overline \rho_a(t,x)
\leq 2 \overline \rho_a(t,x) \rho_a(t,x),~~\hbox{ for $x\ge -\mu_a(t)$.}
\ee
Since $\overline \rho_a(t,x)$ and $\rho_a$ are increasing, we find, 
using (\ref{21may2831}) for $y< - \mu_a(t)$ and (\ref{21may2832}) for $y>-\mu_a(t)$:
\be\label{21may2833}
\begin{split}
&\doverline \rho_a(t,x)
	= \int_{-\infty}^x \overline \rho_a(t,y) \dy
	=\int_{-\infty}^{-\mu_a(t)}\overline\rho_a(t,y)\dy+\int_{-\mu_a(t)}^x\overline\rho_a(t,y)\dy\\
	&\le 2	\int_{-\infty}^{-\mu_a(t)} \rho_a(t,y)dy	+	  2  \int_{-\infty}^x \overline \rho_a(t,y) \rho_a(t,y) \dy
	\le 2\overline \rho_a(t,-\mu_a(t))+ 2 \overline \rho_a(t,x) \int_{-\infty}^x \rho_a(t,y) \dy\\
	&\le 4\rho_a(t,-\mu_a(t))+2 \overline\rho_a^2(t,x)
	\le 4\rho_a(t,x)+4 \overline \rho_a(t,x)^2 \rho_a(t,x)
	\leq 8 \max\{1, \overline\rho_a(t,x)^2\} \rho_a(t,x).
\end{split}	
\ee
Taking the maximum of (\ref{e.c5124}) and (\ref{21may2833}), 
we arrive at~(\ref{21may2823}). We deduce that the function~$\rho_a(t,x)$
satisfies the assumptions of the weighted Nash inequality in \Cref{lem-may2508}, for all~$t>0$.  
This concludes the proof of  Lemma~\ref{lem-may3002}.~$\Box$

\subsubsection*{Step two: pointwise bounds at a particular time}

We now find a ``good time'' $T_g<T$ when $\varphi_a$ satisfies the desired bounds via a time averaging.  
The catch is that we do not have control over $T_g$ and, thus, a third step is required afterwards,
to control the solution on the time interval $T_g\le t\le T$.

\begin{lemma}\label{l.T_g}
For $a>0$ and $\gamma \in (1,2)$, let  $t_1(a)$ be as in~\eqref{21may3004}.  
Given any $T > 4(t_1+1)$, there exists~$T_g \in [T/2, 3T/4]$ such that
	\be\label{e.c696}
		\begin{split}
			(i) \quad \sup_{x \geq - \mu_a(T_g)} \varphi_a(t_T,x) \leq \frac{C}{\sqrt T},
				\qquad\text{and}\qquad
			(ii)\quad \mu_a(T_g) \geq \frac{1}{2} \log(T) - C,
		\end{split}
	\ee
	where $C$ does not depend on $a$ or $\gamma$.
\end{lemma}
{\bf Proof.}  We first note that~\eqref{e.c696}(ii) follows immediately from~\eqref{e.c696}(i).  
Indeed,~\eqref{e.mu_a} yields
\be\label{e.c691}
	\vphi_a(t,-\mu_a(t)) = e^{-\mu_a(t)},
\ee
which, in turn, implies, using~\eqref{e.c696}(i):
\[
\mu_a(T_g)=-\log\vphi_a(t,-\mu_a(t))\ge \farc{1}{2}\log T-C,
\]
which is~\eqref{e.c696}(ii).  Thus, we focus on proving~\eqref{e.c696}(i).

To this end, we use time averages in order to find the ``good time'' $T_g$ 
when the (weighted) $L^2$ bound of $\partial_x \varphi_a$ is small.  Since 
$T/2 \geq 2 t_1(a)$, \Cref{lem-may3002} yields
\be\label{may1802}
	\int \varphi_a^2(t,x)\rho_a(t,x) dx
		\le\farc{C}{\sqrt{t}},
\ee
for all $t \geq T/2$.  Integrating the dissipation inequality (\ref{21may2544}) and using~\eqref{may1802}, we obtain
\[
	\frac{4}{T}\int_{T/2}^{3T/4}\int |\partial_x\varphi_a(s,x)|^2\rho_a(s,x) dxds
		\le\farc{C}{T^{3/2}}. 
\]
As a result, there exists $T_g \in [T/2,3T/4]$ such that
\be\label{e.c5144}
	\int |\partial_x\varphi_a(T_g,x)|^2\rho_a(T_g,x) dx
		\le\farc{C}{T^{3/2}}. 
\ee

Heuristically, we conclude by arguing that, due to~\eqref{e.c5144}, if $\vphi_a$ is ``too big'' somewhere to the right of $-\mu_a(t)$, 
then it must be ``too big'' on a large set.  However, in this region, $\rho_a$ is bounded above and below and $p_a$ and $\varphi_a$ are comparable.  As a result, $p_a$ will be ``too big'' on a large set, violating mass conservation. The key point here that makes the above reasoning work is that, by the definitions \eqref{e.c641}, (\ref{e.c642}), and \eqref{e.mu_a}  of $\rho_a$,  $\vphi_a$, and $\mu_a$, we have
\be\label{e.c694}
\farc{1}{2}\le\rho_a(t,x)\le 1, 
		\quad\text{and}\quad
	\frac{1}{2} \vphi_a(t,x) \leq p_a(t,x) \leq \vphi_a(t,x),
	\qquad\text{ for $x>-\mu_a(t)$, } 
\ee
In particular, working on $x> - \mu_a(t)$ is crucial here as $\rho_a$ is potentially small to the left of $-\mu_a(t)$.



To make this reasoning rigorous, we take any $x_0 \geq -\mu_a(T_g)$, and
use the Newton-Leibniz formula for any $y \in [x_0, x_0 + \sqrt T]$ 
along with~\eqref{e.c5144}-\eqref{e.c694}, to obtain a lower bound for $\vphi_a(T_g,y)$:
\[
	\begin{split}
		\vphi_a(T_g,x_0)
			&\le \vphi_a(T_g,y)
				+\int_{x_0}^{y}|\partial_x\varphi_a|\dx
			\le \vphi_a(T_g,y)
				+\int_{x_0}^{x_0+ \sqrt T}|\partial_x \varphi_a|\dx\\
			&\le 
				\vphi_a(T_g,y)
				+ \sqrt 2 T^{1/4}\Big(\int_{x_0}^{x_0 + \sqrt{T}}|\partial_x\varphi_a|^2 \rho_a \dx\Big)^{1/2}
			\le 
				\vphi_a(T_g,y)
				+ \frac{C}{\sqrt T}.
	\end{split}
\]
Using again~\eqref{e.c694}, this yields
\be\label{e.c695}
	\inf_{x\in [x_0,x_0 + \sqrt T]} p_a(T_g,x)
		\geq \frac{1}{2}\Big(p_a(T_g, x_0) - \frac{C}{\sqrt T}\Big).
\ee
Recalling mass conservation~\eqref{may2516} and using~\eqref{e.c695} we arrive at
\[
	\begin{split}
	1
		&= \int p_a(T_g,x) \dx
		\geq \int_{x_0}^{x_0+\sqrt T} p_a(T_g,x) \dx
		\\&\geq \int_{x_0}^{x_0+\sqrt T} \frac{1}{2}\Big(p_a(T_g, x_0) - \frac{C}{\sqrt T}\Big) \dx
		= \frac{\sqrt T}{2} p_a(T_g,x_0)
			- \frac{C}{2}.
	\end{split}
\]
Rearranging the above and using the arbitrariness of $x_0>-\mu_a(t)$, as well
as (\ref{e.c694}), to translate this into a bound for $\vphi_a(t,x_0)$, finishes the proof
of Lemma~\ref{l.T_g}.~$\Box$.

\subsubsection*{Step three: preserving the $L^\infty$-smallness over $[T_g,T]$}

We are now in a position to combine the results in the first two steps
to finish the proof Lemma~\ref{l:altered u_0 decay}, 
establishing the upper bound (\ref{21may3010}) on $p_a$.

{\bf Proof of \Cref{l:altered u_0 decay}.}  
Combining~\eqref{e.c694}, \Cref{l.T_g}, and~\eqref{e.c6158}, we find $C>1$ such that
\be\label{e.c6145}
	p_a(T_g,x)
		\leq \min\Big( e^x, \frac{C}{\sqrt T}\Big).
\ee
The function
\be\label{21jun1721}
	P(x)
		= \frac{e^x}{1 + \frac{\sqrt T}{2C} e^x}
\ee
is a steady solution to the equation~\eqref{e.c644} that $p_a$ satisfies  as
\be\label{21jun1720}
\partial_x (P^2 e^{-x})
			= \partial_x^2P.
\ee
Moreover, it obeys a uniform bound 
\be\label{e.c6146}
\sup_x \, P(x) \leq \frac{2C}{\sqrt T}.
\ee
The next step is to split $p_a$ into a portion bounded above by $P$ and a small error part. 
Let us write
\be\label{e.c6149}
p_a(t,x)= \psi_P(t,x) + \psi_E(t,x)
			\qquad\text{ for $T_g\le t\le T$.}
\ee
Here, the ``$P$-portion'' $\psi_P\geq 0$ solves
\be\label{e.c6141}
\partial_t \psi_P+ \partial_x (\psi_P^2 e^{-x}) = \partial_x^2 \psi_P,~~T_g\le t\le T,
\ee
with the initial condition
\be\label{e.c6141bis}
\psi_P(T_g,x) = \min\{p_a(T_g,x), P(x)\},~~x\in\Rm.
\ee
The ``error part" $\psi_E\geq 0$ solves
\be\label{e.c6142}
\partial_t \psi_E
+ \partial_x ((2\psi_P + \psi_E)\psi_E e^{-x}) = \partial_x^2 \psi_E,~~T_g\le t\le T,
\ee
with the initial condition
\be\label{e.c6142bis}
\psi_E(T_g,x) = p_a(T_g,x) - \min\{p_a(T_g,x), P(x)\}.
\ee
We now estimate $\psi_P$ and $\psi_E$ at time $T$.
From~(\ref{21jun1720}) and \eqref{e.c6141bis}, it is clear that 
\[
\psi_P(T_g,x)\leq P(x),~~\hbox{ for all $x\in\Rm$},
\]
and both functions solve~\eqref{e.c6141}.  Hence, the comparison principle implies 
that $\psi_P(T,x) \leq P(x)$ for all $x\in\Rm$.  
Using this and~\eqref{e.c6146}, we find 
\be\label{e.c6143}
	\sup_x \, \psi_P(T,x)
		\leq \sup_x \, P(x)
		\leq \frac{C}{\sqrt T},
\ee
as desired.

On the other hand,~\eqref{e.c6142} implies that the total mass of $\psi_E$ is conserved:
\be\label{e.c6144}
\int \psi_E(T,x) dx= \int \psi_E(T_g,x) dx.
\ee
We now bound the right hand side.  First, by~\eqref{e.c6145} and~\eqref{21jun1721},
\[
	p_a(T_g,x)
		\leq \frac{C}{\sqrt T}
		\leq P(x),
		\quad \text{for all } x \geq -\frac{1}{2} \log(T) + \log(2C).
\]
From this and~\eqref{e.c6142}, it follows that
\[
	\psi_E(T_g,x) = 0,
		\quad\text{for all } x \geq -\frac{1}{2} \log(T) + \log(2C).
\]
Hence, using (\ref{e.c6142bis}) and~\eqref{e.c6158}, we obtain 
\be
\bal
\int \psi_E(T_g,x) dx
&= \int_{-\infty}^{-\frac{1}{2}\log(T) + \log(2C)} \psi_E(T_g,x) dx
\leq \int_{-\infty}^{-\frac{1}{2}\log(T) + \log(2C)} p_a(T_g,x) dx\\
&\leq \int_{-\infty}^{-\frac{1}{2}\log(T) + \log(2C)} e^x dx
= \frac{2C}{\sqrt T}.
\enbal
\ee
Invoking~\eqref{e.c6144}, we obtain
\[
	\int \psi_E(T,x) dx
		\leq \frac{2C}{\sqrt T}.
\]
Recall that $T_g \leq 3T/4$ and $T$ is sufficiently large.  
We may, thus, apply parabolic regularity theory (recall that $\psi_E$ solves~\eqref{e.c6142}) to conclude that, up to increasing $C$, we have
\be\label{e.c6147}
	\sup \psi_E(T,x)
		\leq C \int \psi_E(T,x) dx
		\leq \frac{C}{\sqrt T}.
\ee
Combining~\eqref{e.c6149},~\eqref{e.c6143}, and~\eqref{e.c6147}, finishes the proof 
of \Cref{l:altered u_0 decay}.
~$\Box$

\subsection{The weighted $L^2$ framework: the proof of Propositions \ref{lem-may2506} and \ref{lem-may2508}}\label{s:abstract L2}

First, we establish the dissipation inequality in \Cref{lem-may2506}.

{\bf Proof of \Cref{lem-may2506}.} 
Following~\cite{Const}, let us write an equation for $h =H(\vphi)$, with a given function $H$, and
not just for the cases $H(\vphi)=\vphi$ and $H(\vphi)=\vphi^2$ we use here. 
Setting 
\be
\mathcal D_\rho h= \frac{1}{\rho^2} \partial_x ( \rho^2 \partial_x h) - v \partial_xh,
\ee
we find
\be\label{e.c5201}
\bal
		\partial_t h - \mathcal D_\rho h
			&= H' (\vphi) \vphi_t
				- \frac{1}{\rho^2}( \rho^2 H'(\vphi) \vphi_x)_x
				+ v  H'(\vphi)\vphi_x\\
			&= H' (\vphi)\Big( \frac{q_t}{\rho} - \frac{q}{\rho} \frac{\rho_t}{\rho}\Big)
				- 2\frac{\rho_x}{\rho} H'(\vphi) \left( \frac{q_x}{\rho} - \frac{q}{\rho} \frac{\rho_x}{\rho}\right)
				- H''(\vphi) \vphi_x^2\\
				&\quad  - H'(\vphi) \left( \frac{q_{xx}}{\rho} - 2 \frac{q_x \rho_x}{\rho^2} - \frac{q}{\rho} \frac{\rho_{xx}}{\rho}
					+ 2 q \frac{\rho_x^2}{\rho^3}\right)
					+ v H'(\vphi)\left( \frac{q_x}{\rho} - \frac{q \rho_x}{\rho^2}\right).
\enbal
\ee
We now use (\ref{21may2540}) to obtain
\be\label{21may2545}
\bal
	\partial_t h - \mathcal D_\rho h
		&= - H' (\vphi) \frac{q \rho_t}{\rho^2}
				- 2\frac{\rho_x}{\rho} H'(\vphi) \Big( \frac{q_x}{\rho} - \frac{q}{\rho} \frac{\rho_x}{\rho}\Big)
				- H''(\vphi) \vphi_x^2\\
		&\quad +H'(\vphi) \Big(  2 \frac{q_x \rho_x}{\rho^2} + \frac{q}{\rho} \frac{\rho_{xx}}{\rho}
					- 2 q \frac{\rho_x^2}{\rho^3}\Big)
					+ v H'(\vphi) \Big(  - \frac{q \rho_x}{\rho^2}\Big)
					- v_x H'(\vphi) \frac{q}{\rho}\\
		&=
				- H''(\vphi)\vphi_x^2
				- H'(\vphi) \frac{q}{\rho^2}(\rho_t - \rho_{xx} + (v\rho)_x) .
\enbal
\ee
When $H(\vphi) = \vphi$, this establishes~\eqref{e.c697}.

Let us now specialize to the case  $h=H(\vphi) =\vphi^2$ to prove~\eqref{21may2544}.  Multiplying (\ref{21may2545}) by $\rho$ and integrating by parts, we find
\be\label{21may2546}
\begin{split}
\partial_t \int \rho \vphi^2&= \partial_t \int \rho h
= \int \rho_t h+\int \rho\Big(\frac{1}{\rho^2} ( \rho^2 h_x )_x- v h_x
				- 2 \varphi_x^2\Big)
				 - 2\int h (\rho_t - \rho_{xx} + (u\rho)_x)
				\\
			&= \int \rho_t h
				+ \int \rho_x h_x
				+ \int (\rho v)_x h
				- \int 2 \rho \varphi_x^2
				- 2\int h (\rho_t - \rho_{xx} + (u\rho)_x)\\
			&=
				\int h (\rho_t
				- \rho_{xx}
				+ (\rho v)_x)
				- \int 2 \rho \varphi_x^2
				- 2\int h (\rho_t - \rho_{xx} + (v\rho)_x)\\
			&= - \int 2 \rho \varphi_x^2
				- \int h (\rho_t - \rho_{xx} + (v\rho)_x)
			\leq - \int 2 \rho \varphi_x^2,
	\end{split}
\ee
finishing the proof.~$\Box$

Next, we prove the analogue of Nash's inequality.  Due to the inhomogeneity, we are not able to use the 
standard Fourier-based proof.  The proof is more similar to that of Carlen and Loss~\cite{CarlenLoss}, which is  
based on the Poincar\'e inequality.  A complication in our setting is that the decay of $r$ on the left 
makes our space more akin to $\R_+$ than $\R$; however, we only have a mild ``boundary condition'' on the 
left in that we only know that $\varphi^2 r$ is integrable.

{\bf Proof of \Cref{lem-may2508}.} Let us first assume that $\vphi$ is not only non-negative, bounded, and decays to zero sufficiently fast as $x\to+\infty$ but also monotonically decreasing.  Fix $\theta\in \Rm_+$ and notice that, due to the assumptions on $r$, we may find $L \in \Rm$ such that $\bar r(L) = \theta$.  Fix $\eps>0$ to be chosen, and write
\be\label{21may2548}
\begin{split}
0&\leq \int_{-\infty}^L \Big|\vphi_x \sqrt{\doverline r} + \eps \vphi \frac{\overline r}{\sqrt{\doverline r}}\Big|^2 dx
		= \int_{-\infty}^L\Big[ |\vphi_x|^2 \doverline r+ \eps (\vphi^2)_x \overline r + \eps^2 \vphi^2 \frac{\overline r^2}{\doverline r} \Big] dx\\
		&= \int_{-\infty}^L \Big[ |\vphi_x|^2 \doverline r - \eps \vphi^2 r + \eps^2 \vphi^2 \frac{\overline r^2}{\doverline r} \Big] dx
			+ \eps \overline r(L) \vphi^2(L).
\end{split}
\ee
Note that, since $\varphi(x)$ is decreasing, 
 \be\label{e.c563}
\vphi(L) \overline r(L)
\leq \int_{-\infty}^L \vphi r \, dx\leq \int_{-\infty}^\infty\vphi r\, dx.
\ee
Rearranging (\ref{21may2548}) and using~\eqref{e.c563}, we find
\be\label{e.c648}
\eps \int_{-\infty}^{\infty} \vphi^2 r \Big(1 - \eps \frac{\overline r^2}{r \doverline r}\Big)\,dx
		\leq \int_{-\infty}^L |\vphi_x|^2 \doverline r \,dx
			+  \frac{\eps}{\overline r(L)} \Big(\int_{-\infty}^\infty \vphi r \,dx\Big)^2.
\ee
Since $r$ is increasing, we have
\be\label{e.c647}
	\overline r^2(x)
		= 2 \int_{-\infty}^x \bar r(y) \bar r_y(y) \dy
		= 2 \int_{-\infty}^x \bar r(y) r(y) \dy
		\leq 2 r(x) \int_{-\infty}^x \bar r(y) dy
		= 2 r(x) \doverline r(x).
\ee
Hence, using~\eqref{e.c647} and choosing $\eps = 1/4$,~\eqref{e.c648} becomes
%
\be\label{e.c5102}
\bal
\int_{-\infty}^\infty \vphi^2 r\, dx
&\le 8\int_{-\infty}^L |\vphi_x|^2 \doverline r \,dx
			+ \frac{2}{\overline r(L)}\Big( \int_{-\infty}^\infty \vphi r\, dx\Big)^2.
\enbal
\ee
Next, using   assumption (ii)  in (\ref{21may2543}), together with monotonicity of $\bar r(x)$ gives
\be\label{21may2549}
\bal
\int_{-\infty}^\infty \vphi^2 r\, dx
&\le 8C_1\max(1,\bar r(L)^2)\int_{-\infty}^\infty |\vphi_x|^2  r \,dx
			+ \frac{2}{\overline r(L)}\Big( \int_{-\infty}^\infty \vphi r\, dx\Big)^2,
\enbal
\ee
proving (\ref{e.Nashish}) for monotonically decreasing functions $\vphi(x)$ as $\bar r(L) = \theta$.

In order to complete the proof for non-monotonically decreasing non-negative
functions $\vphi(x)$, 
we define an analogue of the decreasing rearrangement as follows.  Let $m_r$ be the measure  
\be
m_r(A)=\int_A r\, dx,
\ee
and for any measurable set $A$ define
\[
A^* = \{x: \bar r(x)<m_r(A)\}. 
\]
Note that $A^*$ is a half line of the form $A^*=(-\infty,x^*(A))$ such that
\be\label{may2602}
m_r(A)=\bar r(x^*(A))=m_r(A^*).
\ee

We may then define rearrangements of functions through the ``layer cake decomposition'':
\[
\vphi^*(x) = \int_0^\infty \1_{\{x : \vphi(x) > t\}^*}(t) dt.
\]
As with the usual symmetric decreasing rearrangement, we immediately find that $\vphi^*$ is 
decreasing, the functions~$\vphi$ and $\vphi^*$ have super-level sets of equal measure, and, for any $p>0$,
\be\label{e.c624}
	\int |\vphi^*(x)|^p r(x) dx
		= \int |\vphi(x)|^p r(x) dx.
\ee

We claim that an analogue of the P\'olya-Szeg\"o inequality holds:
\be\label{e.PolyaSzego}
	\int |\vphi^*_x(x)|^2 r(x) dx
		\leq 	\int |\vphi_x(x)|^2 r(x) dx.
\ee
Since we have already established~\eqref{e.Nashish} for decreasing functions, the fact that~\eqref{e.Nashish} holds for a general function $\vphi$ is an immediate consequence of~\eqref{e.c624} and~\eqref{e.PolyaSzego}.

We now prove~\eqref{e.PolyaSzego}.  In principle, 
it is an immediate consequence of a more general result in~\cite{Talenti}. We present a simpler proof
in our present one-dimensional setting for the convenience of the reader. We may assume without loss of generality that the function $\vphi(x)$ is 
positive everywhere, smooth and
takes each value finitely many times. 
A useful consequence is that any level set of $\vphi^*$ has at most one point.

First, we use the coarea formula to rewrite
\be\label{e.c5253}
	\int \vphi_x(x)^2 r(x) \, dx
		= \int_0^\infty \sum_{x \in \vphi^{-1}(t)} |\vphi_x(x)| r(x) \, dt.
\ee
Then, we note that
\be\label{e.c622}
	\sum_{x \in \vphi^{-1}(t)} r(x)
		\leq \Big( \sum_{x \in \vphi^{-1}(t)} r(x) |\vphi_x(x)|\Big)^{1/2}
			\Big(\sum_{x \in \vphi^{-1}(t)} \frac{r(x)}{|\vphi_x(x)|}\Big)^{1/2},
\ee
so that~\eqref{e.c5253} becomes
\be\label{e.c625}
	\int \vphi_x(x)^2 r(x) \, dx
		\geq \int_0^\infty \frac{\big(\sum_{x \in \vphi^{-1}(t)} r(x)\big)^2}{\sum_{x \in \vphi^{-1}(t)} \frac{r(x)}{|\vphi_x(x)|}} dt.
\ee
Next, we show that both the numerator and denominator in the right hand side of~\eqref{e.c625} 
may be replaced by their analogues with $\vphi^*$ in place of $\vphi$, up to an inequality.

First, since
\[
	\overline r ( (\vphi^*)^{-1}(t))
		= m_r(\{\vphi^* > t\})
		= m_r(\{\vphi > t\})
		\leq \overline r( \sup \vphi^{-1}(t)),	
\]
we conclude that $(\vphi^*)^{-1}(t) \leq \sup \vphi^{-1}(t)$.  Using that $r$ is increasing, we conclude that
\be\label{e.c623}
	\sum_{x \in \vphi^{-1}(t)} r(x)
		\geq r(\sup \vphi^{-1}(t))
		\geq r( (\vphi^*)^{-1}(t))
\ee
On the other hand, by the coarea formula, we have that, for any $t>0$
\[
	m_r(\{x : \vphi(x) >t\})
		= \int_{\{x : \vphi(x) >t\}} r(x) dx
		= \int_t^\infty \Big( \sum_{x \in \vphi^{-1}(s)} \frac{r(x)}{|\vphi_x(x)|} \Big)  ds.
\]
An immediate consequence is that
\[
	\frac{d}{dt} m_r(\{x : \vphi(x) >t\})
		= -\sum_{x \in \vphi^{-1}(t)} \frac{r(x)}{|\vphi_x(x)|}.
\]
The analogous formula for $\vphi^*$ holds.  Since $m_r(\{x : \vphi(x) >t\}) = m_r(\{x : \vphi^*(x) >t\})$ for all $t$, it follows that
\be\label{e.c626}
	\sum_{x \in \vphi^{-1}(t)} \frac{r(x)}{|\vphi_x(x)|}
		= \sum_{x \in (\vphi^*)^{-1}(t)} \frac{r(x)}{|\vphi^*_x(x)|}
		= \frac{r( (\vphi^*)^{-1}(t))}{\vphi^*_x( (\vphi^*)^{-1}(t))}.
\ee
The last equality above is due to the fact that $(\vphi^*)^{-1}(t)$ is a one point set. 

Including~\eqref{e.c623} and~\eqref{e.c626} in~\eqref{e.c625}, we conclude that
\[
	\int \vphi_x(x)^2 r(x) \, dx
		\geq \int_0^\infty \frac{
			r((\vphi^*)^{-1}(t))^2
			}{
			\frac{r( (\vphi^*)^{-1}(t))}{\vphi^*_x( (\vphi^*)^{-1}(t))}
			} dt
		= \int_0^\infty \vphi^*_x( (\vphi^*)^{-1}(t)) r((\vphi^*)^{-1}(t)) dt.
\]
Reapplying the coarea formula to the rightmost quantity above yields
\[
	\int \vphi_x(x)^2 r(x) \, dx
		\geq 	\int \vphi^*_x(x)^2 r(x) \, dx.
\]
which concludes the proof of~\eqref{e.PolyaSzego} and, thus, of \Cref{lem-may2508}.~$\Box$
%
%

\subsection{A lower bound on $p(t,x)$ in the middle region} \label{sec:middle-end} 

We finish this section with a lower bound on $p(t,x)$ in an intermediate region that will be useful to us later on.
\begin{lemma}\label{lem-jun2002}
Let $u(t,x)$ be the solution to the Burgers-FKPP equation (\ref{burgerskpp}) with the initial condition~$u(0,x)=\one(x\le 0)$, and $p(t,x)$ be given by 
(\ref{21jun2002}) and (\ref{apr1327}). 
There exist $k>0$, $T_0>0$ and~$c_0>0$ so that we have
\be\label{21jun2006}
p(t,x)\ge \farc{c_0}{\sqrt{t}},~~\hbox{ for all $t\ge T_0$ and $\disp -\farc{1}{2}\log t\le x\le k\sqrt{t}$.}
\ee
\end{lemma}
{\bf Proof.}  Let us take $t>0$ sufficiently large. We argue by contradiction.  For $\eps \in (0,1/10)$ and~$k>0$ to be chosen, suppose there exists
\be\label{21jun2018}
	x_0 \in \big(-\frac12\log t, k \sqrt t\big)
\ee
such that
\be\label{21jun2020}
p(t,x_0)\le \farc{\eps}{\sqrt{t}},
\ee
whence
\be\label{jul2102}
\hat u(t,x_0)=e^{-x_0}p(t,x_0)\le \farc{\eps}{\sqrt{t}}e^{-x_0}.
\ee
We will show that this violates mass conservation~(\ref{may2516}):
\be\label{21jun2008}
\int p(t,x)dx=\int p(0,x)dx,
\ee
if $\eps$ and $k$ are chosen sufficiently small. 

The first step is to notice that the mass to the left of $(-1/2)\log t$ is small.  
Indeed, due to the upper bound (\ref{apr1337}):
\be\label{21jun2010} 
\int_{-\infty}^{-(1/2)\log t}p(t,x)dx\le\int_{-\infty}^{-(1/2)\log t}e^xdx=\farc{1}{\sqrt{t}}.
\ee
We also know that the mass to the far right is small.  Indeed, recall from (\ref{apr1133}) that
there exist $T_0$ and $N_0$ so that for all $t>T_0$ and
$N>N_0$, we have
\be\label{21jun2012} 
\int_{N\sqrt{t}}^{\infty} p(t,x)dx\le Ce^{-N/2}. 
\ee

To estimate the mass of $p$ in the middle region, we split it into two parts.  
When $x \in [k\sqrt t, N\sqrt t]$, we use~(\ref{jul2102}) together with
the steepness estimate in Lemma~\ref{lem:jun21}  to find
\be\label{21jun2014}
\bal
p(t,x)&=e^x \hat u(t,x)\leq e^x\Big(1 +\Big(\farc{1}{\hat u(t,x_0)}-1\Big) e^{x-x_0}\Big)^{-1}\le 
 e^x\Big(1 +\Big(\farc{\sqrt{t}}{\eps}e^{x_0}-1\Big) e^{x-x_0}\Big)^{-1}\\
 &\le e^x\Big(1 +\farc{\sqrt{t}}{2\eps}e^{x_0} e^{x-x_0}\Big)^{-1}
 \le\farc{2\eps}{\sqrt{t}},~~\text{for}~~x > x_0.
\enbal
\ee
We used (\ref{21jun2018}) in the second inequality above.  It follows that
\be\label{21jun2023}
\int_{k\sqrt{t}}^{N\sqrt{t}}p(t,x)dx\le N\sqrt{t}\farc{2\eps}{\sqrt{t}}=2\eps N.
\ee
On the other hand, when $x \in [-(1/2)\log t, k \sqrt t]$, we use \Cref{lem-may2504} to find
\be\label{e.c7201}
\int_{-(1/2)\log t}^{k\sqrt{t}}p(t,x)dx
	\le 2k\sqrt{t}\farc{C}{\sqrt t}
	= 2k C.
\ee
Putting together~\eqref{21jun2010},~\eqref{21jun2012},~\eqref{21jun2023}, and~\eqref{e.c7201} yields
\[
	\int p(t,x) dx
		\leq \frac{1}{\sqrt t}+Ce^{-N/2}
			+ 2 \eps N
			+ 2k C.
\]
Taking $\eps$ and $k$ sufficiently small, and $N$ sufficiently large 
we obtain a contradiction to~\eqref{21jun2008}. The conclusion of Lemma~\ref{lem-jun2002} follows.~$\Box$

%
%

\section{The proof of Theorem~\ref{thm:main} for $\beta=2$}\label{sec:proof-beta=2}

We now prove Theorem~\ref{thm:main} in the critical case $\beta=2$. As in the case $\beta<2$
considered in Section~\ref{sec:beta<2}, the proof is based
on the analysis in the self-similar variables, and  using the pulled nature of the problem
to show that convergence on the diffusive scales implies convergence to a traveling wave on scales $x\sim O(1)$. 
The key difference with the situation for~$\beta<2$ is that, as we have mentioned previously, the Dirichlet boundary condition  at~$\eta=0$
for the function $v$, introduced by the weighted Hopf-Cole transform (\ref{dec332}), no longer approximately holds  in the self-similar
variables. Instead, the function $v$ has a positive but a priori unknown limit on the left: see Lemma~\ref{prop-june1802} below. 
The bounds in Proposition~\ref{prop-apr1302}
will be a crucial ingredient in establishing the correct boundary condition in the
self-similar variables.  
After we pass to the self-similar variables, the non-zero boundary condition changes the 
long time behavior in the self-similar variables. This is an algebraic reason for adjusting the logarithmic shift $(3/2)\log t$ 
in the front position to $(1/2)\log t$.   With these bounds in hand, the argument has many similarities
with the case $\beta<2$,  so we will omit many of the details, only highlighting the differences.  
As before, we will assume without loss of generality
that the initial condition is~$u(0,x)=\one (x\le 0)$.

\subsubsection*{The weighted Hopf-Cole transform}

Motivated by Proposition~\ref{prop-apr1302}, we consider the moving frame: $x\mapsto x-2t+(1/2)\log(t+1)$, and set
\begin{equation}\label{feb412}
\tilde u(t,x)=u(t,x+2t-\farc{1}{2}\log(t+1)).
\end{equation}
Here, $u(t,x)$ is the solution to the Burgers-FKPP equation (\ref{burgerskpp}),
in the non-shifted reference frame. We stress the difference between the function $\hat u(t,x)$ defined in (\ref{21jun2002}) and
used throughout Section~\ref{sec:twp-bounds}, and the function $\tilde u(t,x)$ used in the present section. The former is 
$u(t,x)$ in the reference frame~$x\to x-2t$,  while the latter  refers to the solution in the
reference frame used in (\ref{feb412}). In particular, the function $p(t,x)$ used throughout Section~\ref{sec:twp-bounds} and $\tilde u(t,x)$ are related by
\be\label{21jun2004}
p(t,x)=e^x\hat u(t,x)=e^x u(t,x+2t)=e^x\tilde u\Big(t,x+\farc{1}{2}\log (t+1)\Big).
\ee
The function
$\tilde u(t,x)$ satisfies
\begin{align}\label{21jun1820}
\tilde u_t - 2\tilde u_x+ \farc{1}{2(t+1)} \tilde u_x+ 2 \tilde u \tilde u_x 
=\tilde u_{xx}+\tilde u(1-\tilde u).
\end{align}
We will use the weighted Hopf-Cole transform, as in (\ref{dec332}) with $\beta=2$: 
\begin{align}\label{jan1302}
    v(t,x) = \exp\Big (x+\int_x^{\infty} \tilde u(t,y) dy \Big )\tilde u(t,x).
\end{align}
This gives
\begin{align}\label{jul1916}
    v_t - v_{xx} +\frac{1}{2(t+1)}(v_x-v) =  v\Big(\int_x^{\infty} \tilde u(1-\tilde u)dy - \tilde u\Big).
\end{align}
The key step in the proof will again be to analyze the long time behavior of $v$ at scales $O((t+1)^\gamma)$ for $\gamma \in (0,1/2)$, and this will be done through the use of self-similar variables.  However, the key difference with the case $\beta < 2$ is that $v$ does not decay to zero as $x\to - \infty$ (cf.~\eqref{e.c7305}).  This is quantified in the following:
\begin{lemma}\label{prop-june1802}
There exist $T_0$, $k_0$, $b_0$, $b_1>0$ so that   
\be\label{21jun1802}
v(t,x)\le b_0<+\infty,~~\hbox{ for all $t\ge T_0$ and any $x$}, 
\ee
and  
\be\label{21jun1804}
0<b_1\le v(t,x),~~\hbox{ for all $x\le k_0\sqrt{t}$ and $t\ge T_0$.}
\ee
\end{lemma}
 {\bf Proof.} Let us define $\tilde\mu(t)$ by 
\be\label{21jun1810}
\tilde u(t,\tilde\mu(t))=\farc{1}{2}.
\ee
Proposition~\ref{prop-apr1302} in terms of $\tilde\mu(t)$ says that
there exists a constant $K>0$ and $T_0>0$ so that
\be\label{21jun1812}
|\tilde \mu(t)|\le K
	\qquad\text{ for all } t \geq T_0.
\ee
Lemma~\ref{lem:jun21} and (\ref{21jun1812}) imply  that
\be\label{21jun1814}
\tilde u(t,x)\le \farc{1}{1+e^{x-K}},~~\hbox{ for $t\ge T_0$ and $x>K$,}
\ee
and
\be\label{21jun1816}
\tilde u(t,x)\ge \farc{1}{1+e^{x+K}},~~x<-K,~~t\ge T_0.
\ee
Hence, for $x\le 0$ we can estimate $v(t,x)$ from above by
\begin{equation}\label{dec1124}
\bal
v(t,x)&=\tilde u(t,x)\exp\Big(x+\int_x^\infty \tilde u(t,y)dy\Big)
\le  \exp\Big(\int_x^0(\tilde u(t,y)-1)dy+\int_0^\infty \tilde u(t,y)dy\Big)\\
&\le  
\exp\Big(\int_0^\infty \tilde u(t,y)dy\Big)\le \exp\Big(K+\int_K^\infty\frac{dx}{1+e^{x-K}}\Big)
\le C,
\enbal
\end{equation} 
which is (\ref{21jun1802}).  The proof for $x\geq 0$ is simpler and follows directly from~\eqref{21jun1814}.

To obtain (\ref{21jun1804}), we  first use (\ref{21jun1816}) to write for $x\le -K$:
\begin{equation}\label{dec1124bis}
\bal
v(t,x)&=\tilde u(t,x)\exp\Big(x+\int_x^\infty \tilde u(t,y)dy\Big)
\ge \farc{1}{2}
\exp\Big(\int_x^0(\tilde u(t,y)-1)dy+\int_0^\infty \tilde u(t,y)dy\Big)\\
&\ge \farc{1}{2}
\exp\Big(-K-\int_{-\infty}^{-K}\Big(1-\farc{1}{1+e^{x+K}}\Big)dy\Big)\ge C.
\enbal
\end{equation} 
In order to  extend this lower bound on the right, we recall the lower bound in Lemma~\ref{lem-jun2002}. Combined with the change of variables
(\ref{21jun2004}), it 
shows that we can find $c_0>0$ and $k_0>0$ so that for $0\le x\le k_0\sqrt{t}$ and $t\ge T_0$ we have 
\be\label{21jun2024}
\tilde u(t,x) =\sqrt{t+1}e^{-x}p(t,x-\farc{1}{2}\log(t+1))\ge c_0 e^{-x}.
\ee
We deduce that
\begin{equation}\label{21jun2025}
\bal
v(t,x)&=\tilde u(t,x)\exp\Big(x+\int_x^\infty \tilde u(t,y)dy\Big)
\ge  c_0,~~\hbox{ for all $t\ge T_0$ and $0\le x\le k_0\sqrt{t}$.}
\enbal
\end{equation} 
Finally, for  $x\in(-K,0)$, we see from the monotonicity of $\tilde u(t,x)$ and (\ref{21jun2024})   that
$\tilde u(t,x)\ge c_0$, so that
\[
v(t,x)=\tilde u(t,x)\exp\Big(x+\int_x^\infty \tilde u(t,y)dy\Big)
\ge  c_0 e^{-K},~~\hbox{ for all $t\ge T_0$.}
\]
This proves (\ref{21jun1804}).~$\Box$

\subsubsection*{Analysis in the self-similar variables }

We now outline the main ingredients required to establish the long-time behavior of $v$.  The thrust of the
argument is similar to the case $\beta < 2$, so the outline will be made in reference to the ideas used in Section~\ref{sec:beta<2}.

First, apply the familiar self-similar change of variables:
\begin{equation}\label{dec1120}
v(t,x)=w\Big(\log(t+1),\farc{x}{\sqrt{t+1}}\Big).
\end{equation}
Then, (\ref{jul1916}) leads to the evolution equation
\begin{align}\label{sep2710}
    w_{\tau} +\mathcal{L}w + \frac{1}{2}e^{-\tau/2}w_{\eta} =e^{\tau}w(\tau,\eta)  \Big(\int_{\eta e^{\tau/2}}^{\infty} \tilde u(1-\tilde u)dy -\tilde u\Big),
\end{align}
with the operator ${\cal L}$ defined by
\begin{align}\label{21jun1818}
    \mathcal{L}w := - w_{\eta \eta}- \frac{\eta}{2}w_{\eta} -\frac{1}{2}w.
\end{align}
This is different from the operator $L$ in (\ref{dec1128}) by a multiple of the identity.  
Considering $\mathcal L$ as an operator on $H^1(e^{\eta^2/4} d\eta; \R_+)$ augmented with Neumann boundary conditions, its spectrum  
consists of the eigenvalues $0,1,2,\dots$ with (unnormalized) principal eigenfunction
\[
	\psi_0(\eta)
		= e^{- \eta^2/4}.
\]
Importantly, Lemma~\ref{lem:neg} shows that the right side of~(\ref{sep2710}) is negative:
\begin{equation}\label{dec1122}
\int_x^{\infty} \tilde u(1-\tilde u)dy -\tilde u\le 0~~\hbox{ for all $x\in\R$.}
\end{equation}

Recall the major ingredients in establishing the long-time dynamics in self-similar variables when $\beta<2$:  
(1) a super-solution solving a tractable equation; (2) a sub-solution solving a (potentially different) tractable equation; 
(3) smallness of the right hand side when $\eta \gg e^{-\tau/2}$; (4) an approximate boundary condition when $\eta \ll - e^{-\tau/2}$.  
We now check that analogous ingredients are available in this case.

First, we find a super-solution.  This follows directly from \Cref{lem:neg}.  Indeed, as a result, we have that any solution to
\begin{align}\label{sep241}
\bar{w}_{\tau} +\mathcal{L}\bar{w} + \frac{1}{2}e^{-\tau/2}\bar{w}_{\eta}  = 0, 
~~\tau>0, ~~\eta>-e^{-(1/2-\gamma)\tau},
\end{align}
augmented with appropriate boundary conditions on $\bar w_\eta(\tau, - e^{-(1/2-\gamma)\tau})$ (see the fourth point, below), is a super-solution to~\eqref{sep2710}.

Second, we see that 
the solution $\under w$ to
\be
	\under{w}_{\tau} +\mathcal{L}\under{w} + \frac{1}{2}e^{-\tau/2}\under{w}_{\eta}  = -e^\tau \under w \tilde u, 
~~\tau>0, ~~\eta> e^{-(1/2-\gamma)\tau},
\ee
is a sub-solution to~\eqref{sep2710}, up to stating an appropriate boundary condition that we do now.  Due to the fact that $\tilde u$ is steeper than $\phi_2$,~\eqref{sep75},~\eqref{dec334}, and \Cref{lem:neg}, notice that
\be\label{e.c80402}
	w_\eta(\tau, \eta)
		= (\tilde u_x + \tilde u(1-\tilde u)) \exp\Big\{\tau/2 + x + \int_x^\infty \tilde u dy\Big\}
		\leq 0.
\ee
Hence, we may take Neumann boundary conditions for $\under w_\eta$:
\be
	\under w_\eta(\tau, e^{-(1/2-\gamma)\tau}) = 0
\ee
and we are guaranteed that, with ordered initial data, $\under w \leq w$.

Third, we establish the smallness of the right hand side of~\eqref{sep2710} whenever $\eta \geq e^{-(1/2-\gamma)\tau}$.  This follows directly from~\eqref{dec1122},  \Cref{prop-june1802}, and~\eqref{21jun1814}:
\be
	\begin{split}
	0 &\geq
		e^\tau \Big( \int_{\eta e^{\tau/2}}^\infty \tilde u(1-\tilde u) dy - \tilde u(e^\tau, \eta e^{\tau/2})\Big) w(\tau, \eta)
	\\&\geq - e^\tau \tilde u(e^\tau, \eta e^{\tau/2}) w(\tau, \eta)
	\geq - C e^{\tau - \eta e^{\tau/2}}
	\geq - C e^{\tau - e^{\gamma\tau}}.
	\end{split}
\ee
Hence, the right hand side of~\eqref{sep2710} is double exponentially small.

Finally, we address the approximate boundary conditions of $\omega$ at $\eta = -e^{-(1/2-\gamma)\tau}$.  Keeping \Cref{prop-june1802} in mind, an approximate Dirichlet boundary condition is not possible.  This is, of course, part of the reason for the different shift when $\beta =2$.  Instead, we have an approximate Neumann boundary condition.  First notice that, due to~\eqref{21jun1816}, we have that
\be\label{e.c80403}
	|1 - \tilde u(e^\tau, -e^{\gamma \tau})| \leq C e^{- e^{\gamma\tau}}.
\ee
By parabolic regularity theory, we find
\be\label{e.c80404}
	|\tilde u_x(e^\tau, -e^{\gamma\tau})|
		\leq C e^{- e^{\gamma\tau}}.
\ee
Using~\eqref{e.c80403} and~\eqref{e.c80404} in~\eqref{e.c80402}, we find
\be
	|w_\eta(\tau, -e^{-(1/2-\gamma)\tau})|
		\leq C e^{\tau - e^{\gamma \tau}}.
\ee
Hence, $w$ satisfies an approximate Neumann boundary condition with double exponentially small error.

Thus, the main ingredients to an analogous argument as in the case $\beta < 2$ are in place, with the only major difference being the change from 
approximate Dirichlet boundary conditions to approximate Neumann boundary conditions.  The change in boundary conditions changes the principal 
eigenfunction of $\mathcal{L}$ and suggests that the long-time behavior of $w$ should look like $\alpha_\infty e^{-\eta^2/4}$ for~$\eta \geq 0$.  
This is confirmed by the following key estimate, which is the analogue of \Cref{lem:omega}.  As the proof follows from similar arguments as in \Cref{lem:omega}, 
using the four ingredients above and our knowledge of the spectrum of $\mathcal{L}$, we omit the details.
\begin{lemma}\label{sep2711}
Given $w$ solving~\eqref{sep2710} with initial conditions $w(0,\eta) = \1(\eta\leq 0)$, there exists a constant $\alpha_\infty > 0$ and a function $R_\gamma$ such that
\be\label{sep245}
	w(\tau,\eta)
		= \alpha_\infty e^{-\eta^2/4} + R_\gamma(\tau,\eta) e^{-\eta^2/6}
		\qquad\text{ for } \eta \geq -e^{-(1/2 - \gamma)\tau},
\ee
for any $\gamma \in (0,1/2)$ and
\be
	\lim_{\tau\to\infty} |R_\gamma(\tau,\eta)| = 0
		\qquad\text{ uniformly for $\eta \geq -e^{-(1/2-\gamma)\tau}$}.
\ee
\end{lemma}
Before proceeding, we note that the role of \Cref{prop-june1802} in the proof of \Cref{sep2711} is to guarantee the positivity of $\alpha_\infty$ and the boundedness of $w$.  From \Cref{sep2711}, we immediately obtain the long-time behavior of $\tilde u$ at scales between $O(1)$ and $O(\sqrt t)$.

 \begin{cor}\label{cor:jun2002}
Let $\gamma\in(0,1/2)$, fix $\eps >0$, and set $x_{\gamma}= (t+1)^{\gamma}$. There exist~$\alpha_\infty>0$ and $T_{\eps,\gamma}$, so that
\begin{align}\label{jul2002}
|\tilde u(t,x_{\gamma}) &-\alpha_{\infty}e^{-x_{\gamma}}|
\le \eps e^{-x_{\gamma}}
	,\quad \text{ for all } t>T_{\eps,\gamma}.
\end{align}
\end{cor}

\subsubsection*{Convergence to a single wave: the proof of Theorem~\ref{thm:main} for $\beta=2$}

With \Cref{sep2711} and \Cref{cor:jun2002} in hand, the rest of the proof is quite straightforward; however, it deviates somewhat from that of the proof for $\beta < 2$ in \Cref{sec:beta<2}.  The main reason being that we cannot verify that an analogously defined $v_\alpha$ satisfies a small boundary condition at $x = -(t+1)^\gamma$ (cf.~\eqref{jul1914}).

As in (\ref{jul2802}), we take
\be\label{jul2006}
	\varphi_{\alpha}(t,x)
		=\phi_2(x + \zeta_{\alpha}(t)),
\ee
and perform the corresponding weighted Hopf-Cole transforms:
\be\label{jul2008}
\psi_{\alpha}(t,x) = e^{\Gamma_{\alpha}(t,x) }\varphi_{\alpha}(t,x),
		\qquad\text{ where }
	\Gamma_{\alpha}(t,x) = x+\int_x^{\infty} \varphi_{\alpha}(t,y) dy.
\ee
Similarly to (\ref{e.c6291}), we  define a shift of the traveling wave, fixing $\zeta_{\alpha}(t)$ by the normalization
\be\label{jul2004}
	\psi_{\alpha}(t, (t+1)^{\gamma}) = \alpha.
\ee
Recalling that $\phi_2(x) = 1/(1+e^x)$, we can easily compute that
\be\label{e.c80407}
	\int_x^\infty \varphi_\alpha(y) dy
		= \int_{x+\zeta_\alpha}^\infty \frac{1}{1+e^y} dy
		= -(x+\zeta_\alpha) + \log(1 + e^{x + \zeta_\alpha}),
\ee
so that
\be\label{e.c80408}
	\alpha
		= \psi_\alpha(t,(t+1)^\gamma)
		= e^{-\zeta_\alpha},
	\qquad\text{ or, equivalently}
	\quad \zeta_\alpha(t) = -\log(\alpha).
\ee

Fix $\eps>0$.  As in the proof of \Cref{thm:main} when $\beta < 2$, it is enough to obtain precise upper bounds on $\vphi_{\alpha_\infty - \eps} - \tilde u$ and $\tilde u - \vphi_{\alpha_\infty + \eps}$ because $|\vphi_{\alpha_\infty} - \vphi_{\alpha_\infty \pm \eps}| \leq C \eps$.  We do this now.

First consider $\varphi_{\alpha_\infty - \eps} - \tilde u$.  In exactly the same manner as we established~\eqref{e.c7299}, we see that
\be\label{e.c80405}
	\varphi_{\alpha_\infty-\eps} - \tilde u < 0,~~\text{ for all $t$ sufficiently large and $x < (t+1)^\gamma$}.
\ee
Here we used \Cref{cor:jun2002}.

We now turn to $\tilde u - \varphi_{\alpha_\infty +\eps}$.  We claim that
\be\label{e.c80406}
	\tilde u - \varphi_{\alpha_\infty + \eps} < C \eps,~~
	 \text{ for all $t$ sufficiently large and $|x| < (t+1)^\gamma$},
\ee
which is enough to conclude the proof.  To this end, we apply \Cref{sep2711} to take $t$ sufficiently large,
so  that $v \leq \alpha_\infty + \eps$ on $|x|< (t+1)^\gamma$.  Then, rewriting $\tilde u$ and using,~\eqref{e.c80407}, ~\eqref{e.c80408}, and~\eqref{e.c80405},
we find, on $|x| < (t+1)^\gamma$,
\be\label{e.c80409}
	\begin{split}
	\tilde u(t,x)
		&= v(t,x) \exp\Big\{-x - \int_x^\infty \tilde u dy\Big\}
		\leq (\alpha_\infty + \eps) \exp\Big\{-x - \int_x^\infty \tilde u dy\Big\}
		\\&
		\leq (\alpha_\infty + \eps) \exp\Big\{-x - \int_x^\infty \varphi_{\alpha_\infty-\eps} dy + \int_{(t+1)^\gamma}^\infty \vphi_{\alpha_\infty-\eps} dy\Big\}\\
		&= (\alpha_\infty + \eps) \frac{e^{\zeta_{\alpha_\infty - \eps}} }{1 + e^{x + \zeta_{\alpha_\infty - \eps}}} e^{Ce^{-(t+1)^\gamma}}
		\\&\leq (\alpha_\infty + \eps) \frac{e^{\zeta_{\alpha_\infty - \eps}} }{1 + e^{x + \zeta_{\alpha_\infty - \eps}}} + Ce^{-(t+1)^\gamma}
		= (\alpha_\infty + \eps)
			\frac{\frac{1}{\alpha_\infty-\eps}}{1 + e^x \frac{1}{\alpha_\infty-\eps}}
				+ Ce^{-(t+1)^\gamma}
		\\&= \frac{1}{\frac{\alpha_\infty-\eps}{\alpha_\infty + \eps} + \frac{1}{\alpha_\infty + \eps} e^x }
			+ Ce^{-(t+1)^\gamma}
		\leq \frac{1}{\frac{1}{\alpha_\infty + \eps} e^x + 1 - \frac{2\eps}{\alpha_\infty}} + Ce^{-(t+1)^\gamma}
		\\&\leq \frac{1}{\frac{1}{\alpha_\infty + \eps} e^x + 1} + C\eps + Ce^{-(t+1)^\gamma}.
	\end{split}
\ee
Using~\eqref{e.c80408}, it is easy to see that
\[
	\vphi_{\alpha_\infty+\eps}(x)
		= \frac{1}{1 + e^{x+\zeta_{\alpha_\infty+\eps}}}
		= \frac{1}{1 + \frac{1}{\alpha_\infty+\eps} e^x}.
\]
Using this along with~\eqref{e.c80409} yields~\eqref{e.c80406}, concluding the proof of \Cref{thm:main} when $\beta = 2$.~$\Box$

\section{Convergence to pushed fronts for $\beta>2$}\label{sec:beta>2} 

In this section, we consider convergence to the minimal speed traveling wave in the case $\beta>2$.
As we have mentioned in the introduction, the proof is quite standard and follows the classical 
approach of~\cite{rothe1981convergence,sattinger1976stability,sattinger1977weighted} for the
convergence to a traveling wave in the pushed front regime.
Let us recall that for all~$\beta>0$ there is a traveling wave solution to the Burgers-FKPP equation
(\ref{burgerskpp}) of the form 
\be\label{apr702}
\phi_{\beta}(x) = \frac{1}{1+e^{{\beta x}/{2}}}
\ee
that moves with the speed 
\be\label{apr704}
c_*(\beta) = \frac{\beta}{2}+\frac{2}{\beta}\ge 2.
\ee
A key point is that  for $\beta\ge 2$, the speed $c_*(\beta)$ given by (\ref{apr704})
also happens to be the minimal front
speed. In other words, for $\beta\ge 2$ the minimal speed wave profile is explicit and 
given by (\ref{apr702}). Another property that will be crucial for the analysis is that
for $\beta>2$ we have 
\be\label{apr706}
\phi_{\beta}(x)e^{c_*x/2}\in L^2(\R).
\ee
According to the criterion of~\cite{garnier2012inside}, this puts the traveling front (\ref{apr702}) into the category of pushed fronts, unlike the
fronts for $\beta\le 2$, for which (\ref{apr706}) does not hold. 
Accordingly, the proof of the large time 
convergence of the solutions to the initial value problem for (\ref{burgerskpp})
with~$\beta>2$ consists of three steps that are common in such results for pushed fronts: 
compactness, local stability and quasi-convergence.

\subsection*{Compactness}

As the first step, we show that the solution $u(t,x)$ to (\ref{burgerskpp}) with
$\beta>2$ can be trapped between an explicitly constructed 
super-solution $\bar{u}$ and a sub-solution $\underline{u}$.
Each of them will converge to a separate shift of the traveling front, exponentially fast in time. 
Here, we follow the construction in~\cite{rothe1981convergence}. We denote by~$c_*=c_*(\beta)$, as given by (\ref{apr704}), 
and drop the subscript $\beta$ in the notation
for the traveling wave profile $\phi_\beta(x)$ given by (\ref{apr702}).
Consider the
moving frame $z=x-c_*t$, setting 
\be\label{apr802}
\hat u(t,z) =u(t, z+c_*t).
\ee
This function satisfies
\begin{align}\label{apr746}
{\cal M}[\hat u]: = \hat u_t-c_*\hat u_x+\beta \hat u \hat u_x -\hat u_{xx} - \hat u(1-\hat u) = 0,
\end{align}
with the initial condition $\hat u(0,z) = u_{\rm in}(z)$. The sub- and super-solutions
in the moving frame are described by the following lemma. 
\begin{lemma}\label{lem:compact}
Fix $\lambda \in (2/\beta, \beta/2)$ and $\under z_0, \overline z_0 \in \R$.  For $\under q_0, \overline q_0, \mu > 0$, let
\be\label{apr710}
	\begin{split}
	&\under{v}(t,z)
		:=\phi(z+\under \xi(t))- \under q(t,z + \under z_0)
	\qquad\text{and}
	\\&\overline v(t,z)
		:= \phi(z-\overline \xi(t))+\overline q(t,z - \overline z_0)
		\qquad\text{ for all } t>0, z\in\R
	\end{split}
\ee
with
\[
	\under q(t,z)
		= \under q_0 e^{-\mu t}
			\min\{\exp\big(-\lambda z\big),1\},
	\quad
	\text{ and }
	\quad
	\overline q(t,z)
		= \overline q_0 e^{-\mu t}
			\min\{\exp\big(-\lambda z\big),1\},
\]
and, for $K>0$,
\be\label{e.c72203}
	\under \xi'(t)
		= K \under q_0 e^{-\mu t}
	\quad \text{ and }\quad
	\overline \xi'(t)
		= K \overline q_0 e^{-\mu t}.
\ee
Suppose, for some $D_0>0$,
\be\label{e.c72202}
	\frac{K \under q_0}{\mu},
		~~\frac{K \overline q_0}{\mu},
		~~ |\under \xi(0) + \under z_0|
		~~ |\overline \xi(0) + \overline z_0|
		\leq D_0.
\ee
Then, if $\under q_0, \overline q_0$ and $\mu$ are sufficiently small and $K$ is sufficiently large, all depending only on $\beta$, $\lambda$, and $D_0$, then $\under v$ and $\overline{v}$ 
are sub- and super-solutions of~\eqref{apr746}, respectively.
\end{lemma}
{\bf Proof.} 
We only prove the claim for $\under v$.  We drop the subscripts of $q_{1,2}$ and $\xi_{1,2}$ and superscripts of $q_0^{(1,2)}$, and $z_0^{(1,2)}$ in order
 to simplify
the notation. 
Let us insert the ansatz for $\under{v}(t,z) $ given by the left side of (\ref{apr710}) into the desired inequality  
\be\label{apr718}
{\cal M}[\under{v}]\leq 0,
\ee
that needs to hold for $\under v(t,x)$ to be a sub-solution to (\ref{burgerskpp}). 
This gives 
\be\label{apr720}
\begin{aligned}
{\cal M}[\under{v} ]&= \xi'(t)\phi' - q_t+\beta(\phi-q)(\phi'-q_z)-c_*(\phi'-q_z)-\phi''+q_{zz}-F(\phi-q)\\
    &=\xi'(t)\phi'-q_t+c_*q_z-\beta\phi q_z - \beta q\phi' + \beta q q_z + q_{zz} + F(\phi)-F(\phi-q).
\end{aligned}
\ee
Here, we have set
\[
F(u)=u(1-u).
\]
We fix a sufficiently small $\delta>0$, and consider (\ref{apr720}) in three regions of $z\in\Rm$
separately.

\textbf{The far right region:}
\be\label{apr726}
R=\{z:~\phi(z-\xi(t))\leq \delta\}.
\ee
In other words $z \geq \xi(t) + \phi^{-1}(\delta)$.  Thus,
\be
	z + z_0
		\geq \phi^{-1}(\delta) + \xi(t) + z_0.
\ee
By~\eqref{e.c72203}-\eqref{e.c72202},
\be
	|\xi(t) + z_0| \leq 2D_0.
\ee 
Hence, by choosing $\delta>0$ sufficiently small depending only on $D_0$, we have
\be
	z + z_0 \geq 0,
\ee
so that  
\be\label{apr731}
q(t,z+z_0)=q_0\exp\{-\mu t-\lambda(z-z_0)\},~~\hbox{ for all $z\in R$}.
\ee
It follows that 
\[
q_t = -\mu q, ~~ q_z= -\lambda q, ~~q_{zz} = \lambda^2 q,~~\hbox{ for all $z\in R$}.
\]
As 
$\xi'(t)>0$, $\phi'<0$, and $\phi(t,x)<\delta$ for 
all $x\in R$, we deduce that 
\be\label{apr722}
\begin{aligned}
{\cal M}[\under{v}]&= \xi'(t)\phi'+ (\mu-c_*\lambda+\lambda^2+\beta\lambda\phi-\beta \phi')q + F(\phi) - F(\phi-q) -\lambda\beta q^2\\
    &\leq (\mu-c_*\lambda+\lambda^2+\beta\delta\lambda-\beta \phi')q + F(\phi) - F(\phi-q).
\end{aligned}
\ee
Notice that the traveling wave satisfies
\begin{align}\label{aug16}
c_*\phi-\frac{\beta}{2}\phi^2 = -\phi' + \int_x^{\infty}\phi(1-\phi) dy,
\end{align}
we see that
\be\label{apr724}
\phi' = -c_*\phi+\frac{\beta}{2}\phi^2  + \int_x^{\infty}\phi(1-\phi) dy\geq -c_*\phi~~~\hbox{ for all $x\in \R$}.
\ee
In addition, we have
\be\label{apr725}
F(\phi)-F(\phi-q)=\phi-\phi^2-\phi+q+(\phi-q)^2=q-2\phi q+q^2\le (1+q_0)q. 
\ee
Using these inequalities, together with (\ref{apr726}), in (\ref{apr722}), gives 
\be\label{apr723}
\begin{aligned}
{\cal M}[\under v]&\leq (\lambda^2 - c_*\lambda + \beta\delta\lambda + \mu +1+q_0+c_*\beta \delta )q. 
\end{aligned}
\ee
Due to the assumption on $\lambda$, we have $\lambda^2 - c_*\lambda + 1 < 0$. 
Then take $\delta>0$, $\mu>0$ and $q_0>0$ all sufficiently small, depending on $\lambda$, so that
\be\label{apr728}
{\cal M}[\under v]\le 0.
\ee
We stress that here we use the assumption that $\beta>2$, which makes $c_*>2$.

\textbf{The middle region:} 
\be\label{apr729}
M=\{z:~\delta\leq \phi(z-\xi(t))\leq 1-\delta\}.
\ee
There exists $\alpha_\delta>0$, which depends on $\delta>0$, so that for all $z$ in the middle region we have  
\be\label{e.c72201}
\phi'(z-\xi(t))\leq -\alpha_\delta.
\ee
We need to consider separately the points where $z\ge z_0 $, or $z<z_0 $, as this changes the definition of~$q(t,z)$. 
If $z>z_0$, then $q$ satisfies (\ref{apr731}), and we use~\eqref{e.c72201} to obtain, using positivity
of $\xi'(t)$ once again,
\be\label{apr730}
\begin{aligned}
{\cal M}[\under{v}] &=\xi'(t)\phi'-q_t+c_*q_z-\beta\phi q_z - \beta q\phi' + \beta q q_z + q_{zz} + F(\phi)-F(\phi-q)
\\
&\leq -\alpha_\delta\xi'(t)+ (\mu+\lambda^2-c_*\lambda+\beta\phi\lambda+c_*\beta \phi)q + (1+q)q -\lambda\beta q^2
\leq - \alpha \xi'(t) + C q_0 e^{-\mu t}.
\end{aligned}
\ee
Above we used that $F(\phi) - F(\phi - q) \leq q(1+q)$.  Due to the assumption on $\lambda$, $\lambda \beta > 1$, which implies that $q^2 - \lambda \beta q^2 \leq0$.  It follows that $C$ does not depend on $q_0$.

On the other hand, if $z\leq z_0$, then
\be\label{apr732}
q_t=-\mu q,~~q_z=0,
\ee
and, as long as we choose $q_0<\delta$, we get, after another appplication of (\ref{apr724}):
\be\label{apr734}
	\begin{split}
	 {\cal M}[\under{v}]
	 	\le \xi'(t)\phi'+(\mu+c_*\beta\phi) q + F(\phi) - F(\phi-q)
		&\le -\alpha_\delta\xi'(t)  + (\mu +c_*\beta+1)q
 	\\&\leq - \alpha_\delta \xi'(t) + C q_0 e^{-\mu t}.
	\end{split}
\ee 
It follows that $M[\under v]\le 0$ both in (\ref{apr730}) and (\ref{apr734}) since, by~\eqref{e.c72203}
\begin{align}\label{eqn:xi}
\xi'(t) = K q_0 e^{-\mu t}
	\geq \frac{C}{\alpha_\delta} q_0e^{-\mu t},
\end{align}
after possibly increasing $K$ depending only on $\delta$, which, in turn, depends on $D_0$ and $\lambda$.

\textbf{The far left region:} 
\be\label{apr740}
L=\{z:~\phi(z-\xi(t))\geq 1-\delta\}.
\ee
Arguing as in the far right region setting, we see that, up to further decreasing $\delta$ depending on $D_0$, we have $z+z_0 \leq 0$.  
%
%
Therefore, $q(t)$ satisfies (\ref{apr732}), and we have 
\be\label{apr741}
\begin{aligned}
 {\cal M}[\under{v}] &=\xi'(t)\phi'-q_t+c_*q_z-\beta\phi q_z - \beta q\phi' + \beta q q_z + q_{zz} + F(\phi)-F(\phi-q)\\
 &=\xi'(t)\phi' + (\mu-\beta\phi') q + F(\phi) - F(\phi-q).
\end{aligned}
\ee
The explicit form (\ref{apr702}) of the profile $\phi(x)$ as $x\to-\infty$ implies that in the far left region we have
\be\label{apr742}
- \frac{\beta}{2}\delta\leq \phi'<0.
\ee
As we also have $F'(1) = -1$, we can assure that
\begin{align}\label{sep1711}
 M[\under{v}] \leq \Big(\mu+\farc{\beta^2\delta}{2}\Big) q + F(\phi) - F(\phi-q)\leq 0,
\end{align}
as long as we choose $\mu>0$ and $\delta>0$ sufficiently small. 

Thus, in all cases, we have show that $M[\under v]\leq 0$, which finishes the proof of Lemma~\ref{lem:compact}.~$\Box$ 

\subsection*{Local stability}

Lemma~\ref{lem:compact} implies a certain compactness for the solution $v(t,z)$ to (\ref{apr746}) 
that eventually will lead to the nonlinear stability result. To formulate it, we need to introduce the 
weighted Banach space 
\[
B_{\lambda} = \{u\in C(\R,\R): \|u\|_{\lambda}<+\infty\},
\]
with the norm
\begin{align}
\|v\|_{\lambda} = \max_{z\in \R}\frac{|v(z)|}{\eta(z)}
	\qquad \text{ where } \eta(z) = \min\{1, \exp(-\lambda z)\}.
\end{align}
We also introduce the norm 
\be\label{apr747}
\|v\|_{\lambda,1} = \|v\|_{\lambda} + \|v_z\|_{\lambda}
\ee
for a  later usage. 
Let us define the $\omega$-limit set of $u_{\rm in}\in B_{\lambda}$ with respect to the evolution (\ref{apr746}) as
\begin{align*}
    \omega(u_{\rm in}) := \{\psi \in B_{\lambda}\cap C^2(\R): ~~ \exists~ t_n\to +\infty ~\text{ s.t }~
    v(t_n,\cdot)\to \psi(\cdot) \text{ in } C_{loc}^2(\Rm) \}.
\end{align*}
By the standard parabolic regularity theory (see, e.g.,~\cite[Chapter IV]{Lieberman}), we know that the ``orbit" $\{v(t, \cdot): t\ge 1\}$ is relatively compact in $C^2_{loc}(\Rm)$. 
Thus, the set $\omega(u_{\rm in})$ is nonempty for a given initial condition~$u_{\rm in}\in B_{\lambda}$. 
Proposition~\ref{prop:aug262} implies that each element of $\omega(u_{\rm in})$ is either a constant
$0$ or $1$, or a shift of the 
traveling wave $\phi(x)$. Indeed, $v$ is simply $u$ in the moving frame $c_* t$.  Hence, if $m_\beta(t_n) - c_* t_n \to \infty$, then \Cref{prop:aug262} implies that $v(t_n,\cdot) \to 1$ as the front $m_\beta(t_n)$ is ``far ahead.''  Similarly, if $m_\beta(t_n) - c_* t_n \to -\infty$, then $v(t_n,\cdot) \to 0$.  Finally, if $m_\beta(t_n) - c_* t_n$ converges to a constant, then \Cref{prop:aug262} implies that $v(t_n,\cdot)$ converges to (a shift of) the traveling wave.

Lemma~\ref{lem:compact} rules out the possibility of convergence to a constant.  
Our goal is to prove that it consists of exactly one element, and that the solution
converges exponentially fast to that particular shift of a traveling wave. 

Lemma~\ref{lem:compact} implies the following local stability result. 
\begin{lemma}\label{lem-apr702}
Fix $\lambda \in (2/\beta, \beta/2)$.  For every $\eps>0$, there exists $\gamma>0$ such that whenever 
\be\label{apr749}
\|v(0,\cdot) - \phi(\cdot -s_0)\|_{\lambda}<\gamma,
\ee
we have
\be\label{apr748}
\|v(t,\cdot) - \phi(\cdot-s_0)\|_{\lambda}<\eps~~\hbox{ for all $t\ge 0$.}
\ee
\end{lemma}
{\bf Proof.} Let us assume that $s_0=0$. We see from (\ref{apr749}) that 
\be\label{apr751} 
\phi(z) - \gamma \eta(z)< v(0,z)< \phi(z) + \gamma \eta(z) ~~\text{ for all } z\in\R.
\ee
This allows us to apply Lemma \ref{lem:compact} with $\xi(t)$ satisfying (\ref{eqn:xi}) and $\xi(0) = 0$:
\[
\xi(t)=\farc{K\gamma}{\mu}\Big(1-e^{-\mu t}\Big),
\]
 leading to
\begin{align*}
    \phi(z-\xi(t)) - \gamma \eta(z)< v(t,z)< \phi(z + \xi(t)) + \gamma \eta(z),~~t\geq 0, ~z\in\R.
\end{align*}
Because 
\[
\phi(z)\sim e^{-{\beta z}/{2}}~~~\hbox{ as $x\to +\infty$,}
\]
and $\lambda<{\beta}/{2}$, there exists a constant $K_1$ such that
\begin{align*}
\|\phi(\cdot)- \phi(\cdot \pm \xi(t))\|_{\lambda}\le K_1 \xi(t) \le K_1 \gamma, ~~t\geq 0,
\end{align*}
finishing the proof.~$\Box$
 
\subsection*{Quasi-convergence and exponential stability} 
 
Now, convergence in shape in Proposition~\ref{prop:aug262} and the steepness comparison in \Cref{prop-jul1502}, together with Lemma~\ref{lem-apr702} 
imply the following quasi-convergence property.  
\begin{cor}\label{cor-apr802}
Fix $\lambda \in (2/\beta, \beta/2)$.  Let $v(t,x)$ be the solution to (\ref{apr746}) with the initial condition $v_0(x)=\one(x\le 0)$.
Then there exists a sequence $t_n\to+\infty$ and $s_0\in\Rm$ so that 
\be\label{apr802bis}
\lim_{n\to+\infty}\|v(t_n,\cdot) - \phi(\cdot -s_0)\|_{\lambda}=0.
\ee
\end{cor}

To improve the result of Corollary~\ref{cor-apr802} to the exponential convergence, we will use  the 
method of~\cite{sattinger1976stability, sattinger1977weighted}, 
which consider convergence to traveling waves for 
equations of the form 
\begin{align*}
    u_t = u_{xx} + f(u,u_x).
\end{align*}
In the Burgers-KPP case, we have
\[
f(u,u_x) = -\beta u u_x + u(1-u).
\]
If one looks for solutions that are perturbations of a traveling wave, that is
\begin{align*}
    u(t,x) = \phi(x-c_* t) + w(t,x),
\end{align*}
then in the moving coordinate $z = x-c_*t$, the equation of $w$ can be decomposed as
\begin{align*}
    \frac{\d w}{\d t} = {\cal L}w + R(w).
\end{align*}
Here, the linearized operator is
\begin{align}\label{apr910}
{\cal L}w = w_{zz} +c_* w_z + \frac{\d f}{\d\phi'} w_z + \frac{\d f}{\d \phi} w = w_{zz}+(c_*-\beta \phi) w_z + (-\beta \phi'+1-2\phi)w.
\end{align}
The remainder $R$ is a nonlinear operator whose Fr\'echet derivative vanishes at $w= 0$. 

The main result of~\cite{sattinger1977weighted} is that a local stability result, as we have 
in Corollary~\ref{cor-apr802}, implies exponential convergence to a traveling wave,
as long as the operator ${\cal L}$ satisfies certain spectral assumptions
that we will now recall and specify to the Burgers-FKPP equation.  
%
%
%
Let us put the linearized operator in the from
\begin{align}\label{apr912}
 {\cal L}w = w'' - 2b w' +q w,~~b=\farc{1}{2}(\beta\phi-c_*),~~q=-\beta\phi'+1-2\phi,
\end{align}
and set
\[
B(z) = \int_0^z b(s) ds.
\] 
To remove the drift term in the operator ${\cal L}$ in (\ref{apr912}), 
Sattinger introduces the operator 
\[
\widetilde{\cal M} = e^{-B} {\cal L} e^{B},
\]
so that
\begin{align}
 \widetilde{\cal M}w = w'' + pw,~~p= b'-b^2+q.
\end{align}
In our case,  with $\phi(x)$ given explicitly by (\ref{apr702}), and $c_*$ by (\ref{apr704}), we have 
the left and right limits
\[
 \begin{aligned}
	&b_{+}=\lim_{x\to+\infty}b(x)=-\farc{c_*}{2}=-\frac{\beta}{4}-\farc{1}{\beta} ,
	\qquad
	b_{-}=\lim_{x\to-\infty}b(x)=\farc{1}{2}(\beta-c_*)=\farc{1}{2}\Big(\farc{\beta}{2}-\farc{2}{\beta}
\Big),\\
	& q _\pm=\lim_{x\to\pm\infty}q(x)= \pm 1.
\end{aligned}
\]
We see that
\[
p_{\pm} = \lim_{z\to \pm\infty} p (z)
\]
are given by
\[
p_+ = -b_+^2 + q_+ = -\Big(\frac{\beta}{4}+\frac{1}{\beta}\Big)^2+1=-\Big(\frac{\beta}{4}-\frac{1}{\beta}\Big)^2 <0,
\]
as $\beta>2$, and 
\[
p_- =-b_-^2+q_-= - \Big(\frac{1}{\beta}+\frac{\beta}{4}\Big)^2-1<0 .
\]
%
%
%
As both $p_+<0$ and $p_-<0$, the operator $\widetilde{\cal M}$ is stable both as $x\to+\infty$ and $x\to-\infty$,
and Theorems 1 and 2 in \cite{sattinger1977weighted} imply that convergence
in Corollary~\ref{cor-apr802} is actually exponential in time.
\begin{cor}\label{cor-apr804}
Let $v(t,x)$ be the solution to (\ref{apr746}) with the initial condition $v_0(x)=\one(x\le 0)$.
There exist $C>0$ and $\omega_0>0$, as well as $s_0\in\Rm$ so that 
\be\label{apr804}
\|v(t,\cdot) - \phi(\cdot -s_0)\|_{\lambda}\le C e^{-\omega_0t}.
\ee
\end{cor}
Using that $\|\cdot\|_{L^\infty} \leq \|\cdot\|_\lambda$, this completes the proof of Theorem~\ref{thm:main} for $\beta>2$.

\section{An informal derivation of the higher corrections}\label{sec:higher-order}

We explain in this section how the non-rigorous but extremely interesting
methodology  of \cite{berestycki2018new} can be used to predict the higher 
order corrections to the logarithmic shift in the front position. This strategy was applied in~\cite{berestycki2018new}
for the classical Fisher-KPP equation, and  the first two extra 
terms in the expansion were rigorously confirmed in~\cite{Graham,NRR2},
leading to the long-time front position asymptotics
\begin{align}\label{jul2238}
m(t) = 2t-\frac{3}{2}\log(t+1)+x_{\infty} -\frac{3\sqrt{\pi}}{\sqrt{t+1}}+\frac{9}{8}(5-6\log2)\frac{\log (t+1)}{t+1}+o\Big(\farc{\log t}{t}\Big),
\end{align} 
a significant refinement of Theorem~\ref{thm:main}. 
The terms that appear in (\ref{jul2238}), except for $x_\infty$, do not depend on the initial conditions.
Moreover, they 
are expected to be universal for a large class of Fisher-KPP type problems. That is, this expansion 
has been shown to hold, with exactly the same coefficients,
for all equations of the form
\be\label{jul2240}
u_t=u_{xx}+f(u),
\ee
with a nonlinearity $f(u)$ of the Fisher-KPP type, normalized so that $f'(0)=1$ as long as $f \in C^{1,\delta}$ near $0$ (see~\cite{BouinHenderson} for the treatment of the less regular case). 

An analysis similar to what we do in this section for $\beta=2$ would show that (\ref{jul2238}) holds
with {\em exactly} the same coefficients for the Burgers-FKPP equation (\ref{burgerskpp})
for all $\beta<2$, confirming the universality prediction of~\cite{berestycki2018new}.  We omit the details.   
The goal of this section is to show that when~$\beta=2$ this expansion changes the coefficients
but not the form of the individual terms to 
\begin{align}\label{jul2243}
m(t)=2t-\frac{1}{2}\log (t+1)-x_{\infty} -\frac{\sqrt{\pi}}{2\sqrt{t+1}}
+\frac{(1-\log 2)}{4} \frac{\log (t+1)}{t+1}+o\Big(\farc{\log t}{t}\Big).
\end{align}
It is tempting to conjecture
that the expansion (\ref{jul2243}) is also universal for problems that combine the pulled nature with traveling waves that decay as $e^{-x}$ rather
than the Fisher-KPP asymptotics~$xe^{-x}$. This regime is referred to as Case (II) in Chapter 2 of 
\cite{Leach-Needham} for reaction-diffusion equations of the form~(\ref{jul2240}).


We begin with the following observation, inspired by \cite{berestycki2018new}, that we state for all $\beta\le 2$. Let $u(t,x)$ be the solution
to (\ref{burgerskpp}) with an initial condition $u_{\rm in}(x)$  and set 
\begin{align}\label{aug110}
    \varphi(t,r) = \int_\R u^2(t,z)e^{rz}dz,~~\Phi_1(r)=(1-r\beta/2) \int_0^{\infty} \varphi(t,r) e^{- (r^2+1)t} dt,
\end{align}
and 
\begin{align}\label{aug112}
    \Phi(r) = \int_{\R} u_{\rm in}(x) e^{rx} dx.
\end{align}
Note that the function $\Phi(r)$ is  smooth for all $r>0$, as long as the initial  condition~$u_{\rm in}(x)$ is compactly supported on the right. 
On the other hand, as we will see below, the function $\Phi_1(r)$ may potentially blow up as $r\to 1^-$. This possibility is removed by the
following identity. 
\begin{proposition}\label{prop-jul2102}
Assume that the initial condition $u_{\rm in}(x)$ for (\ref{burgerskpp}) satisfies
\begin{equation}\label{dec102}
u_{\rm in}(x)\leq \frac{1}{1+e^{x-L_0}},~~\hbox{ for all $x\in\Rm$},
\end{equation} 
with some $L_0>0$. Then, we have
\begin{align}\label{relation}
\Phi(r) =\Phi_1(r), ~~\hbox{for all~$\beta\leq 2$ and $r\in(0,1)$}. 
\end{align}
\end{proposition}
{\bf Proof.} The proof is a modification of the argument in  \cite{berestycki2018new}. Without loss of generality we assume that $L_0=0$.
First, we recall that, as we have shown in (\ref{dec321})-(\ref{dec704}), for all $\beta\le 2$, we have the upper bound
\begin{equation}\label{dec104}
u(t,x)\leq \bar u(t,x):=\frac{1}{1+e^{x-2t}},~~\hbox{ for all $x\in\Rm$}.
\end{equation} 
It follows that $\varphi(t,r)$ is defined and differentiable in $r$ for all $r<2$.
This also shows that for all~$r\in(0,1)$ we  have
\begin{equation}\label{dec108}
 g(t,r): = \int_{\R} u(t,x) e^{rx} dx\le\int_\Rm\farc{e^{rx}}{1+e^{x-2t}} dx=I_0(r)e^{2rt},
\end{equation}    
with
\[    
I_0(r)=\int_\Rm\farc{e^{rx}}{1+e^{x}} dx<\infty,
\]
since $r\in(0,1)$. Next, 
differentiating $g(t,r)$ in $t$, we obtain  
\begin{equation}
\begin{aligned}
g_t(t,r) &= \int_{\R} u_t (t,x)e^{rx} dx = \int_{\R} \Big(-\frac{\beta}{2}(u^2)_x + u_{xx} + u-u^2\Big) e^{rx} dx\\
    & = (1+r^2) g(t,r)-(1-r\beta/2)\varphi(t,r).
\end{aligned}
\end{equation}
Integrating this identity from $0$ to $t$ and using the definition of $\Phi(r)$, we see that
\begin{align*}
    g(t,r) e^{-(1+r^2)t} = \Phi(r)-(1-r\beta/2)\int_0^{t}\varphi(s,r) e^{-(1+r^2)s} ds.
\end{align*}
As $2r<1+r^2$, we may use (\ref{dec108}) to pass to the limit $t\to+\infty$ and obtain (\ref{relation}).~$\Box$


An immediate consequence of Proposition~\ref{prop-jul2102} is that the function $\Phi_1(r)$ remains regular 
and even infinitely differentiable, with bounded derivatives as~$r\to 1^-$. A surprising discovery of~\cite{berestycki2018new}
is that for the classical Fisher-KPP equation one can perform a careful analysis of that limit in terms of the
front location $m(t)$ and this regularity alone 
can be used  to obtain the asymptotics~(\ref{jul2238}) of~$m(t)$ as~$t\to+\infty$.   
As noted above, this approach can be applied nearly verbatim for $\beta<2$, so we only consider~$\beta=2$ below.  

\subsubsection*{Assumptions on the rate of convergence}

Let us now formalize the assumptions that go into the derivation of (\ref{jul2238}) and
(\ref{jul2243}). 
We know from Proposition~\ref{prop:aug262} that there is a reference frame $m(t)$ such that 
\begin{equation}\label{dec110}
u(t,x+m(t))\to \phi_2(x),
\end{equation}
and from  Theorem~\ref{thm:main} that, when $\beta=2$:
\be\label{jul2016}
m(t)=2t-\gamma(t),~~\gamma(t)=a\log (t+1)-\alpha(t),~~\alpha(t)=x_\infty+o(1)~~\hbox{as $t\to+\infty$}.
\ee
Of course, we already know that $a=1/2$ when $\beta=2$ but we leave this coefficient undetermined for now, to show how this value can be discovered
by the arguments below. 
One assumption of~\cite{berestycki2018new} and
later proved in~\cite{Graham}, 
is that the analogue of~(\ref{dec110}) with~$\beta=0$
holds at the rate $O(1/t)$.
For the Burgers-FKPP equation this translates into  the assumption that 
\begin{equation}\label{jul2120}
u(t,x+m(t))= \phi_2(x)+\farc{1}{t}\eta(t,x),~~\hbox{ as $t\to+\infty$,}
\end{equation}
with a rapidly decaying in space (but not necessarily in time) function $\eta(t,x)$. A result of~\cite{Graham} for~$\beta=0$ is that $\eta(t,x)$
has a positive limit as $t\to+\infty$, and the rate~$O(t^{-1})$ cannot be improved.  We stress that (\ref{jul2120}) is an assumption and
not a rigorous claim, even though we believe that it holds, as it does for the classical Fisher-KPP equation. 

Let us now  re-write the function $\Phi_1(r)$ in terms of $m(t)$. We introduce
\begin{align}\label{aug111}
\varphi_m(t,r) = \int_\R u(t,z+m(t))^2e^{rz}dz=e^{-rm(t)}\int_\R u(t,z)^2e^{rz}dz=e^{-rm(t)}\vphi(t,r),
\end{align}
and write
\be\label{jul2018}
\Phi_1(r)=(1-r) \int_0^{\infty} \varphi_m(t,r) e^{rm(t)- (r^2+1)t} dt.
\ee
 %
Given the exponential decay of $u(t,x)$ and $\phi_2(x)$, and assumption (\ref{jul2120}), we deduce that there exists a function $\tilde\eta(t,r)$, so that 
\begin{equation}\label{dec112}
\varphi_m(t,r)= \int_{\R} u(t, z+m(t))^2 e^{rz} dz=\tilde\vphi(r)+\farc{1}{t+1}\tilde\eta(t,r),
 \end{equation}
with
\be\label{jul2121}
\tilde\varphi(r):= \int_\Rm\phi_2(z)^2e^{rz}dz,
\ee
and
\be\label{jul2123}
|\tilde\eta(t,r)|\le K,~~\hbox{ for all $t>0$ and $r\in(1/2,3/2)$. }
\ee
Inserting (\ref{dec112}) into (\ref{jul2018}) gives
\be\label{jul2222}
\Phi_1(r)=(1-r) \int_0^{\infty} \Big(\tilde\vphi(r)+\farc{1}{t+1}\tilde\eta(t,r)\Big) e^{rm(t)- (r^2+1)t} dt.
\ee

\subsubsection*{Higher order corrections from the limit $r\to 1^-$}

We now pass to the limit $r\to 1^-$ in (\ref{jul2222}) and use the regularity of $\Phi_1(r)$ in this limit, implied by Proposition~\ref{prop-jul2102},
to get extra terms in the asymptotic expansion for $m(t)$ as $t\to+\infty$.
We take~$r=1-\eps$ in (\ref{jul2222}), with~$\eps\in(0,1)$ and use expression
(\ref{jul2016}) for $m(t)$:
\begin{equation}\label{dec116}
\begin{aligned}
\Phi_1(1-\eps)& =\eps\int_0^{\infty} \Big(\tilde\vphi(1-\eps)+\farc{1}{t+1}\tilde\eta(t,1-\eps)\Big) 
e^{(1-\eps)(2t-\gamma(t)) - ((1-\eps)^2+1)t} dt\\
&=\tilde\vphi(1-\eps)I(\eps)+  E(\eps).
\enbal
\ee
Here, the main term is
\be\label{jul224}
I(\eps)=\eps\int_0^{\infty}  
e^{-(1-\eps)\gamma(t)-\eps^2t} dt,
\ee
and
\be\label{jul225}
E(\eps)=\eps\int_0^{\infty} \farc{1}{t+1}\tilde\eta(t,1-\eps)   e^{-(1-\eps)\gamma(t)-\eps^2t}  dt
\ee
is the error term. 
%

In order to avoid additional technicalities in an argument that is not rigorous (because assumption (\ref{jul2120})
has not been justified), we use (\ref{jul2120}) to replace the rigorous claim that $\Phi_1(1-\eps)$ remains regular as $\eps\to 0$ by the 
assumption that
\be\label{jul2226}
I(\eps)\hbox{ remains strictly positive and regular as $\eps\to 0$,}
\ee
neglecting the error term $E(\eps)$. Note that the function $\tilde\vphi(r)$ is regular near $r=1$, as can be seen immediately from its definition (\ref{jul2121}).
The positivity  of $I(\eps)$ holds because we know from (\ref{relation}) that $\Phi_1(1-\eps)$ is finite and  not small for any $\eps>0$, and so is $\tilde\vphi(1)$. 

The surprising fact is that (\ref{jul2226}) by itself leads to the asymptotic 
expansion (\ref{jul2243}) for $m(t)$. 
We first explain how we can find from (\ref{jul2226}) that the coefficient $a$ that appears in (\ref{jul2016}) equals to $1/2$, as expected.  
Using the expression for $\gamma(t)$ in (\ref{jul2016}), we write
\begin{equation}\label{dec118}
\begin{aligned} 
I(\eps)&
	=\eps \int_0^{\infty} \farc{e^{-\eps^2t+(1-\eps)\alpha(t)}}{(t+1)^{a(1-\eps)}}  dt
	=\eps \int_0^{\infty} \farc{e^{-s+(1-\eps)\alpha(s/\eps^2)}}{(s/\eps^2+1)^{a(1-\eps)}} \farc{ds}{\eps^2}\\
	&=\eps\big(1+o(1)\big)e^{(1-\eps)x_\infty} \int_0^{\infty} \farc{\eps^{2a(1-\eps)-2}}{(s+\eps^2)^{a(1-\eps)}}e^{-s} ds
	\\&
	=\eps^{2a(1-\eps) - 1}\big(1+o(1)\big)e^{(1-\eps)x_\infty} 
	\int_{\eps^2}^{\infty} \farc{e^{-r + \eps^2}}{r^{a(1-\eps)}} dr,
\end{aligned}
\end{equation}
where the $o(1)$ is according to the limit $\eps\to0$.  We replaced $\alpha(s/\eps^2)$ by its limit $x_\infty$ in second-to-the last step above.
If $a> 1$, then the non-integrability of $r^{-a}$ near the origin will cause the integral to grow like $\eps^{-2(a-1)}$.  In this case, we deduce that $I(\eps) \sim \eps$, which cannot happen as $I(\eps)$ is uniformly positive~\eqref{jul2226}.  A similar argument applies to the case $a=1$.  Hence, $a<1$.  In this case, the integral is finite, so that $I(\eps) \sim \eps^{2a-1}$.  Since $I$ is positive and bounded (again, by~\eqref{jul2226}), the only choice is $a=1/2$ as in \Cref{thm:main}.

Next, we show how the remaining terms in (\ref{jul2243}) come about. 
Taking $a=1/2$ in (\ref{dec118}) and writing
\be\label{jul2228}
\alpha(t) = x_{\infty} + p(t),~~p(t)\to 0\hbox{ as $t\to+\infty$}, 
\ee
leads to
\begin{equation}\label{jul2112}
\begin{aligned}
I(\eps)& =\eps e^{(1-\eps)x_\infty}\int_0^{\infty} \farc{1}{\sqrt{t+1}}e^{(\eps/2)\log (t+1)+(1-\eps)p(t)-\eps^2t}dt.
\end{aligned}
\end{equation}
Let us assume without loss of generality that $x_\infty=0$, and note that 
$I(\eps)$ has the limit
\begin{align}\label{jul2233}
I(\eps)\to \int_0^{\infty}  \farc{e^{-t}}{\sqrt{t}}dt=\sqrt{\pi},~~\hbox{ as $\eps\to 0$},
\end{align}
so the positivity requirement in (\ref{jul2226})  holds automatically. Expanding further gives 
\be\label{jul2229}
\bal
I(\eps)&=\sqrt{\pi}+\eps\int_0^{\infty} \farc{1}{\sqrt{t+1}}\Big(\farc{\eps}{2}\log (t+1)+ p(t)\Big)
e^{-\eps^2t}dt+o(\eps\log\eps).
\enbal
\ee
Note that the first term in the integral above has the asymptotics 
\be\label{jul2231}
\bal
	\farc{\eps^2}{2}&\int_0^{\infty} \farc{\log (t+1)}{\sqrt{t+1}}   e^{-\eps^2t}dt
	= \farc{\eps^2}{2}\int_1^{\infty} \farc{\log t}{\sqrt{t}}e^{-\eps^2t}dt+o(\eps\log\eps)
	\\&=\frac{\eps}{2}\int_{\eps^2}^\infty\farc{\log(t/\eps^2)}{\sqrt{t}}e^{-t}dt+o(\eps\log\eps)
=-\sqrt{\pi}\eps\log\eps+o(\eps\log\eps).
\enbal
\ee
We now make an ansatz 
\be\label{jul2230}
p(t)=\farc{b}{\sqrt{t}}+o\big(t^{-1/2}\big),~~\hbox{ as $t\to+\infty$},
\ee
with the constant $b$ to be determined, and insert it into the second integral in (\ref{jul2229}):
\be\label{jul2232}
\bal
\eps\int_0^{\infty} \farc{p(t)}{\sqrt{t+1}}e^{-\eps^2t}dt =\eps b\int_1^\infty \farc{e^{-\eps^2t}}{t}dt+o(\eps\log\eps)=\eps b\log(\eps^{-2})+o(\eps\log\eps).
\enbal
\ee
Taking into account the regularity of $I(\eps)$ together with (\ref{jul2229})-(\ref{jul2231}) 
and (\ref{jul2232}), we conclude that, since the overall coefficient in front of $\eps\log\eps$ 
has to cancel, we must have
\be\label{jul2234}
b=\farc{\sqrt{\pi}}{2}.
\ee
We deduce that when $\beta=2$ the front location has the asymptotics 
\be\label{jul2235}
m(t)=2t-\farc{1}{2}\log t-x_\infty-\farc{\sqrt{\pi}}{2\sqrt{t}}+o(t^{-1/2}),~~\hbox{ as $t\to+\infty$,}
\ee
recovering all but the last term in the right side of (\ref{jul2243}).
 
A similar computation, expanding $I(\eps)$ further in $\eps$, and
observing that the terms of the order~$\eps^2\log\eps$ must cancel, improves (\ref{jul2235}) to
\begin{align}\label{jul2236}
m(t)=2t-\frac{1}{2}\log t-x_{\infty} -\frac{\sqrt{\pi}}{2\sqrt{t}}+\frac{(1-\log 2)}{4} \frac{\log t}{t}+....
\end{align}
which is (\ref{jul2243}). 
We leave the computational details to an interested reader.

  \begin{appendix}

\section{Traveling waves for the Burgers-FKPP equation}\label{sec:phase-plane}

In this appendix, we recall some basic facts on the Burgers-FKPP traveling waves.
Most of them can be found in Section 13.4 of \cite{murray2007mathematical}, at least on a formal level,
and also in~\cite{vladimirova2006flame} for $\beta<0$. 

We will use the notation 
\be\label{aug406}
f(u)=u(1-u)
\ee
for the reaction term, and write the traveling wave Burgers-FKPP equation as
\be\label{aug402}
-cU'+\beta UU'=U''+f(U),~~~U(-\infty)=1,~U(+\infty)=0.
\ee

\subsection*{Existence of the traveling waves}

We first state the result on the existence of the traveling waves.
\begin{prop}\label{prop-aug402}
A traveling wave solution   to (\ref{aug402}) exists for all $c\ge c_*(\beta)$, with
\begin{equation}\label{speed}
c_* (\beta)= \left\{\begin{matrix}
2\sqrt{f'(0)} =2, ~~~&\text{if } \beta<2,\cr
\dfrac{\beta}{2}+\dfrac{2}{\beta},
&\text{if }\beta\geq 2.\cr\end{matrix}\right.
\end{equation} 
\end{prop}
{\bf Proof.} First, we note that the boundary conditions in (\ref{aug402}) allow us to use the 
sliding method to deduce that any traveling wave solution $U(x)$ to (\ref{aug402}) is monotonically
decreasing.   
Introducing~$V=-U'>0$, we write (\ref{aug402}) as a system
\begin{equation}\label{aug404}
\frac{dU}{dx} = -V,~~\frac{dV}{dx} = (\beta U-c)V + f(U).
\end{equation}
This leads to 
\begin{align}\label{ratio}
    \frac{dV}{dU} = c-\beta U -\frac{f(U)}{V}.
\end{align}
Consider the plane formed by the horizontal $U$-axis and the vertical $V$-axis. 
The explicit form of~$f(u)$ in (\ref{aug406}) shows that the system (\ref{aug404}) has two equilibrium points 
\[
\hbox{$E_1=(0,0)$ and~$E_2=(1,0)$.}
\]
A traveling wave is a heteroclinic orbit of (\ref{aug404}) that goes from $E_2$ to $E_1$. 

The linearization of~(\ref{aug404}) around $E_1$ is 
\begin{align*}
\frac{d}{dx} \begin{pmatrix}
\tilde U \\
\tilde V 
\end{pmatrix} = \begin{pmatrix}
0 & -1\\
f'(0) & -c
\end{pmatrix}
\begin{pmatrix}
\tilde U \\
\tilde V 
\end{pmatrix}.
\end{align*}
The eigenvalues of this matrix are 
\be\label{aug408}
\lambda_{1,2} = \frac{-c\pm\sqrt{c^2-4}}{2},
\ee
and are both real and negative if $c\geq 2$. Thus, the point $E_1$ is a stable equilibrium of $c\ge 2$. 

On the other hand, the linearization of (\ref{aug402}) around $E_2$ is
\begin{align*}
\frac{d}{dx} \begin{pmatrix}
\tilde U \\
\tilde V 
\end{pmatrix} = \begin{pmatrix}
0 & -1\\
f'(1) & \beta-c
\end{pmatrix}
\begin{pmatrix}
\tilde U \\
\tilde V 
\end{pmatrix}.
\end{align*}
The eigenvalues of this matrix are 
\be\label{aug410}
\lambda_{\pm} = \frac{\beta-c\pm\sqrt{(\beta-c)^2+4}}{2}.
\ee
They are real 
and have opposing signs, so $E_2$ is a saddle point.   

\subsubsection*{Existence of traveling waves for $c\ge c_*(\beta)$}

We prove existence of a heteroclinic orbit by constructing an invariant region in the $(U,V)$-plane
that contains
the unstable manifold of the point $(1,0)$. The regions are different for $\beta\ge 2$ and~$\beta<2$,
so we consider them separately. 

When $\beta\geq 2$, we consider the region $D_1$ formed by the interval
$\ell_1= [0,1]\times\{0\}$ along the~$U$-axis, 
and the curve 
$\ell_2 = \{V = (\beta/{2})f(U)\}$ that connects the equilibrium points $E_1$ and $E_2$
and lies in the upper half-plane $\{V>0\}$. Note that the slope of $\ell_2$ at the point $E_2$ 
is $(-\beta/2)$, while the slope~$e_u$ of the unstable orbit at $E_2$ is $(-\lambda_+)$. Hence, the
unstable orbit starting at $E_2$ enters the region~$D_1$ at the point $E_2$ if $\lambda_+<\beta/2$. 
This condition is satisfied if
\be\label{aug416}
\lambda_{+} = \frac{\beta-c+\sqrt{(\beta-c)^2+4}}{2}<\farc{\beta}{2},
 \ee
or 
\be\label{aug418}
c>c_*(\beta)=\farc{\beta}2+\farc{2}{\beta}.
\ee
Now, we check that the region $D_1$ is invariant. 
Along the interval $\ell_1$, we have, from (\ref{aug404}):  
\be 
\frac{dU}{dx} = 0, ~~ \frac{dV}{dx} = f(U)>0,
\ee
so the trajectories point upward, into $D_1$. Along the curve $\ell_2$, we have, from (\ref{ratio}):  
\be\label{aug2426}
\frac{dV}{dU} = c-\beta U -\frac{2}{\beta}.
\ee
The slope of the curve $\ell_2$ itself is 
\be\label{aug2425}
\frac{\beta}{2}f'(U) = \beta\Big(\frac{1}{2}-U\Big).
\ee
Thus,  the trajectories along $\ell_2$ point into $D_1$ if
\begin{align}\label{ineq1}
    c-\beta U -\frac{2}{\beta}> \beta\Big(\frac{1}{2}-U\Big),~~\hbox{ for all $U\in[0,1]$,}
\end{align}
or, equivalently, 
\be\label{aug412}
c> c_*(\beta)=\frac{\beta}{2}+\frac{2}{\beta}.
\ee
This shows that traveling waves exist for all $c>c_*(\beta)$ when $\beta\ge 2$.
It is easy to see that the curve~$\ell_2$ itself forms a heteroclinic connection between $E_1$ and $E_2$ when $c = c_*(\beta)$.

When $0<\beta<2$, we consider a different region $D_2$ formed by the same interval $\ell_1 = [0,1]\times\{0\}$ 
and the curve $\ell_3= \{V = f(U)\}$, also connecting the points $E_1$ and $E_2$ through the upper half-plane.  
The slope of $\ell_3$ at $E_2$ is $(-1)$, hence the condition that the unstable direction $e_u$ points into the region~$D_1$ at
the point $E_2$ is $\lambda_+<1$, which is
\be\label{aug419}
\lambda_{+} = \frac{\beta-c+\sqrt{(\beta-c)^2+4}}{2}<1.
\ee
This condition holds as soon as $c>\beta$.

To check that the region $D_2$ is invariant, we first note that,
as for $\beta\ge 2$, the trajectories point upward and inside $D_2$ along the interval $\ell_1$. The slope of the curve $\ell_3$ is $1-2U$, and we have,
from~(\ref{ratio}):  
\begin{align}\label{ineq2}
 \frac{dU}{dx}<0,~~   \frac{dV}{dU} = c-\beta U - 1 \geq 1-2U,~~\hbox{for all $U\in[0,1]$},
\end{align}
as long as $c\geq c_*=2$. We use here the fact that $\beta<2$. Thus, the region $D_2$ is invariant. As the unstable direction $e_u$ enters $D_2$,
existence of the heteroclinic orbit connecting $E_2$ and $E_1$ follows. 

\subsubsection*{Non-existence of traveling waves for $c<c_*(\beta)$}

Consider first the case $\beta\le 2$.  Suppose that $c< c_*(\beta) = 2$. Given a wave $U(x)$ that satisfies (\ref{aug402}), define the normalized translates 
\[
U_n(x)=\farc{U(x+n)}{U(n)},
\]
that satisfy
\be\label{aug420}
-cU_n'+\beta U(n)U_nU_n'=U_n''+U_n-U(n)U_n^2,~~U_n(0)=1.
\ee
As $U(n)\to 0$ as $n\to+\infty$, the Harnack inequality and normalization $U_n(0)=1$ imply that the sequence $U_n(x)$ converges locally uniformly to a limit
$\bar U(x)>0$ that satisfies
\be\label{aug421}
-c\bar U' =\bar U''+\bar U,~~\bar U(0)=1.
\ee
Since $c<2$, solutions to (\ref{aug421}) are exponentials of the form
\[
\bar U(x)=A_1e^{\lambda_1x}+A_2e^{\lambda_2x},
\]
with $\lambda_{1,2} \in \C \setminus \R$ given by (\ref{aug408}).  On the other hand, $\bar U$ is real and positive.  This is only possible if~$A_1 = A_2 = 0$, which violates the normalization $\bar U(0) = 1$.  This is a contradiction, and, thus,~$c<2$ is not possible when $\beta \leq 2$.

In the case $\beta\ge 2$ we argue as follows to show that no wave can exist for $c<c_*(\beta)$. First, the same argument as for $\beta<2$ implies that
there exist no waves for $c<2$. Let us assume that there there exists a heteroclinic orbit connecting $E_2$
to $E_1$ for some $c<c_*(\beta)$.
Then, (\ref{aug408}) implies that its slope at the point $E_1=(0,0)$ is
\be\label{aug422}
\farc{dV}{dU}\Big|_{(0,0)}=\frac{c-\sqrt{c^2-4}}{2}.
\ee
On the other hand, its slope at the point $E_2$ is still given by $(-\lambda_+)$, and we see from (\ref{aug416}) that it goes out of the region
$D_1$ as it leaves $E_2$. Moreover, looking at (\ref{aug2426}) and (\ref{aug2425}) we see that if~$c<c_*(\beta)$, 
then the heteroclinic orbit can not re-enter the domain $D_1$ along the curve $\ell_2$ at a point with $0<U<1$.
Therefore, if it arrives to $E_1$,  it has to do that above the curve $\ell_2$, and its slope at this point has to satisfy
\be\label{aug423}
\farc{dV}{dU}\Big|_{(0,0)}>\frac{\beta}{2}. 
\ee 
Combining (\ref{aug422}) and (\ref{aug423}) gives
\be\label{aug424}
\frac{c-\sqrt{c^2-4}}{2}>\frac{\beta}{2}. 
\ee 
It follows that
\be\label{aug427}
\beta<\farc{4}{c+\sqrt{c^2-4}}\le\farc{4}{c}\le 2,
\ee
as we already know that $c\ge 2$. This is a contradiction to the assumption that $\beta\ge 2$, finishing the proof of Proposition~\ref{prop-aug402}.~$\Box$ 


\paragraph{Remark.} When $c=2$ and $\beta \in (0,2)$, the above proof shows that for $\beta\in (0,2)$, all trajectories are trapped in the region bounded by the curves $\ell_2 = \{V = ({\beta}/{2})f(U)\}$ 
and $\ell_3 = \{V = f(U)\}$, that is $D_2 \setminus D_1$.  The reasons for this are that the inequality in (\ref{ineq1}) is reversed and (\ref{ineq2}) holds, so this region is invariant, and $\beta/2 < \lambda_+ < 1$, so trajectories enter $D_2 \setminus D_1$ from $E_2$. This implies that for any $0<z<1$ we have 
\begin{equation}\label{dec922}
\frac{z(1-z)}{|\bar{E}(z)|}\leq \frac{2}{\beta},
\end{equation}
with $\bar E(z)$ defined in (\ref{dec132}).

\subsubsection*{The asymptotic profile of the traveling waves} 

We briefly summarize asymptotic profiles of  traveling wave solutions. The details are essentially identical to the Fisher-KPP equation,
see, for instance,~\cite{Leach-Needham} for a detailed analysis.  

When $\beta < 2$, the minimal speed traveling wave  has the following asymptotics on the right:
\begin{align}\label{asymp1}
    U(z) = (Az+B)e^{-z} + O(e^{-(1+\delta)z}), ~~\text{as } z\to +\infty. 
\end{align}
with $A>0$ and $\delta>0$. 
%
When $\beta \geq 2 $, the critical front has an explicit form
\begin{align}
    U(z) = \frac{1}{1+e^{\beta z/2}}.
\end{align}
%

\subsubsection*{Steepness comparison of the waves}

Finally, we present a steepness comparison for the traveling waves that is used to show convergence of the solution in shape to the minimal speed traveling wave.
\begin{lemma}\label{lem-apr902}
Fix any $\beta\in\Rm$.  Let $U_1$ and $U_2$ be two traveling waves with speeds $c_1$ and $c_2$, respectively.  
If $c_1\leq c_2$, then $U_1$ is steeper than $U_2$.
\end{lemma}
{\bf Proof.} The proof follows an approach from \cite{fife1977approach}. If we define 
\be
	\bar E_i(z) = - U_i'( U_i^{-1}(z)),
\ee
it is enough to show that 
\be\label{aug434}
\bar E_1(z) \leq \bar E_2(z), ~~\hbox{ for all $z\in(0,1)$.}
\ee  
Using \eqref{ratio}, we obtain
\begin{align}
	\bar E_i'(z) + \beta z + \frac{f(z)}{\bar E_i(z)} = c_i,
\end{align}
so that 
\begin{align}
	\bar E_1'(z)-\bar E_2'(z) -\frac{f(z)}{\bar E_1(z)\bar E_2(z)}(\bar E_1(z)-\bar E_2(z)) = c_1-c_2.
\end{align}
It follows that  the function 
\begin{align*}
	F(z) = (\bar E_1(z)-\bar E_2(z))\exp\Big(-\int_{1/2}^z \frac{f(z')}{\bar E_1(z')\bar E_2(z')}dz'\Big),
\end{align*}
satisfies 
\begin{align}\label{aug436}
	F'(z) = (c_1-c_2)\exp\Big(-\int_{1/2}^z \frac{f(z')}{\bar E_1(z')\bar E_2(z')}dz'\Big).
\end{align}
Observe that $\bar E_1-\bar E_2\to 0$ as $z\to 1$. It follows that $F(z)\to 0$, as well.   
On the other hand, since~$c_1\leq c_2$, we see from (\ref{aug436}) that 
$F(z)$ is decreasing.  It follows that $F(z)\leq 0$ for all $z$, which implies (\ref{aug434}) and concludes the proof.~$\Box$

\end{appendix}

\bibliographystyle{plain}

\begin{thebibliography}{10}

\bibitem{Angenent} S. Angenent, {The zero set of a solution of a parabolic equation}. J. Reine Angew. Math. {\bf 390}, 79-96, 1988.

\bibitem{ABK1} L.-P. Arguin, A. Bovier, and N. Kistler, {Poissonian statistics in the extremal process of branching Brownian motion}. Ann. Appl. Probab. {\bf 22},  
1693--1711, 2012. 


\bibitem{ABK2} L.-P. Arguin, A. Bovier, and N. Kistler, {The extremal process of branching Brownian motion}. Probab. Theory Relat. Fields {\bf 157}, 535--574, 2013.

\bibitem{AG}
M. Avery and L. Gar\'enaux, Spectral stability of the critical front in the extended Fisher-KPP equation,  Preprint arXiv:2009.01506v1, 2020.

\bibitem{AS1} M. Avery and A.  Scheel,  
Asymptotic stability of critical pulled fronts via resolvent expansions near the essential spectrum, 
SIAM J. Math. Anal. {\bf 53}, 2206--2242, 2021.  

\bibitem{AS2} M. Avery and A.  Scheel,  
Universal selection of pulled fronts, Preprint arXiv:2012.06443, 2020. 



 


\bibitem{berestycki2017exact}
J. Berestycki, {\'E}. Brunet, and B. Derrida.
  Exact solution and precise asymptotics of a Fisher--KPP type front, 
 {Jour. Phys. A: Math. Theor.} {\bf 51}, 035204, 2018.


\bibitem{berestycki2018new}
J. Berestycki, {\'E}. Brunet, and B. Derrida.
  A new approach to computing the asymptotics of the position of
  {F}isher-{KPP} fronts,
 {EPL (Europhysics Letters)}  {\bf 122}, 10001, 2018.


\bibitem{BouinHenderson}
E. Bouin and C. Henderson.
 The Bramson delay in a Fisher-KPP equation with log-singular non-linearity,
 {Nonlinear Anal.} {\bf 213}, 112508, 2020.

\bibitem{BramburgerHenderson}
J. Bramburger and C. Henderson.
The speed of traveling waves in a FKPP-Burgers system,
{Arch. Ration. Mech. Anal.} {\bf 241} (2021), no. 2, 643–681.
 

\bibitem{Bramson1} 
M. D. Bramson, {Maximal displacement of branching Brownian motion}, 
{Comm. Pure Appl. Math.} {\bf 31}, 531--581, 1978.

\bibitem{Bramson2} 
M. D. Bramson, {\it Convergence of solutions of the Kolmogorov equation
  to travelling waves}, { Mem. Amer. Math. Soc.} {\bf 44}, 1983. 

\bibitem{BD1} E. Brunet and B. Derrida. {Statistics at the tip of a branching random walk and the delay of traveling waves}. Eur. Phys. Lett. {\bf 87}, 60010, 2009.

\bibitem{BD2} E. Brunet and B. Derrida.  A branching random walk seen from the tip, 
Jour. Stat. Phys. {\bf 143}, 420--446, 2011. 

\bibitem{CarlenLoss}
E. Carlen and M. Loss,
 Sharp constant in Nash's inequality,
 {Int. Math. Res. Not.}, 1993, no. 7, 213--215.

\bibitem{Carlen1996optimal}
E. Carlen and M. Loss,
Optimal smoothing and decay estimates for viscously damped conservation laws, with applications to the 2-D Navier-Stokes equation,
{Duke Math. Jour.} {\bf 31}, 135--157, 1996.

\bibitem{Const}
P. Constantin, Generalized relative entropies and stochastic representation, 
Int. Math. Res. Not. 2006, Art. ID 39487, 9 pp.


\bibitem{Davies}
E.B. Davies, 
{\em Heat kernels and spectral theory},
Cambridge Tracts in Mathematics, 92. Cambridge University Press, Cambridge, 1989.


\bibitem{Ebert-vanSaarlos}
U. Ebert and W. van Saarloos, Front propagation into unstable states: universal algebraic convergence towards uniformly pulled fronts, Physica D {\bf 146}, 1-99,  2000.


\bibitem{fife1977approach}
P.C. Fife and J.B. McLeod,
The approach of solutions of nonlinear diffusion equations to travelling front solutions,
{Arch. Rat. Mech. Anal.} {\bf 65},
335--361,
  1977.
  
  

\bibitem{Fisher}
R. A. Fisher, The wave of advance of advantageous genes, {Ann. Eugen.} {\bf 7}, 355--369, 1937. 


\bibitem{FuGrietteMagal1}
X. Fu, Q. Griette, and P. Magal,
A cell-cell repulsion model on a hyperbolic Keller-Segel equation.
J. Math. Biol. {\bf 80}, 2257--2300, 2020. 


\bibitem{FuGrietteMagal2}
X. Fu, Q. Griette, and P. Magal,
Sharp discontinuous traveling waves in a hyperbolic Keller-Segel equation.
Math. Models Methods Appl. Sci. {\bf 31}, 861--905, 2021. 


\bibitem{garnier2012inside}
J. Garnier, T. Giletti, F. Hamel, and L. Roques.
 Inside dynamics of pulled and pushed fronts.
Jour. Math. Pures Appl.,
  {\bf 98}, 428--449, 2012.
  
\bibitem{Giletti}
T. Giletti, Monostable pulled fronts and logarithmic drifts, Preprint arXiv:2105.12611, 2021. 
  
\bibitem{GM}
T. Giletti and H. Matano, Existence and uniqueness of propagating terraces, 
 Comm. Contemp. Math. {\bf 22}, 1950055, 38 pp., 2020.  
  
  

\bibitem{Graham}
C. Graham, Precise asymptotics for Fisher--KPP fronts, {Nonlinearity }
{\bf 32}, 1967--1998, 2019.  

\bibitem{HadelerRothe}
K. Hadeler and F. Rothe, Travelling fronts in nonlinear diffusion equations, J. Math. Biol. {\bf 2},
251--263, 1975. 



\bibitem{HNRR}
F. Hamel, J.  Nolen, J.-M. Roquejoffre and L. Ryzhik, 
A short proof of the logarithmic Bramson correction in Fisher-{KPP} equations,
{Netw. Heterog. Media} {\bf 8}, 275--289, 2013.



%

\bibitem{henderson2016}
C. Henderson,
Population stabilization in branching Brownian motion with absorption and drift,
{Commun. Math. Sci.} {\bf 14}, 973--985, 2016. 


\bibitem{henderson_chemotaxis}
C. Henderson,
Slow and fast minimal speed traveling waves of the FKPP equation with chemotaxis,
arXiv preprint arXiv:2102.06065,
2021.


\bibitem{howard2002}
 P. Howard,
 Pointwise estimates and stability for degenerate viscous shock waves
{J. reine angew. Math.} {\bf 545}, 19--65, 2002. 


\bibitem{kolmogorov1937etude}
A.N. Kolmogorov, I.G. Petrovskii, and N.S. Piskunov,
{\'E}tude de l'{\'e}quation de la diffusion avec croissance de la
  quantit{\'e} de mati{\`e}re et son application {\`a} un probl{\`e}me
  biologique.
{Bull. Univ. Moskow, Ser. Internat., Sec. A} {\bf 1}, 1--25, 1937.

\bibitem{Lau} K.-S. Lau, {On the nonlinear diffusion equation of Kolmogorov, Petrovskii and Piskunov}, {J.~Diff. Eqs.} {\bf 59}, 44--70, 1985. 

\bibitem{Lieberman}
G.M. Lieberman,
{\em Second order parabolic differential equations}. World Scientific Publishing Co., Inc., River Edge, NJ, 1996.

\bibitem{Lofting}
H. Lofting, {\it The Story of Doctor Dolittle}, New York: Frederick A. Stoke Co., 1920.

\bibitem{LeachHanac2016} 
J. Leach and E. Hana{\c{c}},
{On the evolution of travelling wave solutions of the Burgers-Fisher equation}, 
{Quart. Appl. Math.} {\bf 74}, 337--359, 2016. 

\bibitem{Leach-Needham}
J.A. Leach and D.J. Needham, {\it Matched Asymptotic Expansions in Reaction-Diffusion Theory},
Springer-Verlag London, 2004. 



\bibitem{slimemolds}
M.T. Keating and J.T. Bonner,
Negative chemotaxis in cellular slime molds,
J. Bacteriology {\bf 130}, 144--147, 1977. 



\bibitem{McK} H.P. McKean, Application of Brownian motion to the equation of Kolmogorov-Petrovskii-
Piskunov, Comm. Pure Appl. Math. {\bf 28}, 323--331, 1975. 

\bibitem{MPP}
P. Michel,  S. Mischler,and B. Perthame, General entropy equations for structured population models 
and scattering, C. R. Math. Acad. Sci. Paris {\bf 338}, 697--702, 2004.  

\bibitem{murray2007mathematical}
J.~D Murray.
\newblock {\em Mathematical biology: I. An introduction}, Third ed., Vol.~17, 
Interdisciplinary Applied Mathematics Series, Springer New York, 2002.


\bibitem{MRR}
L. Mytnik, J.-M. Roquejoffre and L. Ryzhik, Fisher-KPP equation with small data and the extremal process of branching Brownian motion,
arXiv:2009.02042, 2020.



\bibitem{NRR1}
J.  Nolen, J.-M. Roquejoffre and L. Ryzhik,
Convergence to a single wave in the Fisher--KPP equation, {Chin. Ann. Math. Ser. B} {\bf 38},
629--646, 2017. 
 
\bibitem{NRR2}
J.  Nolen, J.-M. Roquejoffre and L. Ryzhik,
Refined long-time asymptotics for Fisher--KPP fronts, {Comm. Contemp. Math.}, 2018, 1850072. 
 
 
  
\bibitem{Roberts} M. Roberts, 
{A simple path to asymptotics for
the frontier of a branching Brownian motion}, {Ann. Prob.} {\bf 41}, 
3518--3541, 2013. 


 

\bibitem{rothe1981convergence}
F. Rothe, 
Convergence to pushed fronts, Rocky Mount. Jour. Math. {\bf 11}, 617--633,
  1981.

\bibitem{sattinger1976stability}
D.H Sattinger,
On the stability of waves of nonlinear parabolic systems.
{Adv. Math.}, {\bf 22}, 312--355, 1976. 

\bibitem{sattinger1977weighted}
D.H.~Sattinger,
Weighted norms for the stability of traveling waves,
Jour. Diff. Eqs. {\bf 25}, 130--144, 1977. 

%
\bibitem{Talenti}
G. Talenti
A weighted version of a rearrangement inequality.
Ann. Univ. Ferrara Sez. VII (N.S.), {\bf 43}, 121--133, 1997. 

\bibitem{vanSaarlos}
W. van Saarlos, Front propagation into unstable states, Physics Reports,
{\bf 386}, 29-222, 2003. 



\bibitem{Uchiyama} K. Uchiyama, The behavior of solutions of some nonlinear 
diffusion equations for large time, {J.~Math. Kyoto Univ.} {\bf 18}, 453--508,
1978. 


\bibitem{vladimirova2006flame}
N. Vladimirova, G. Weirs, and L. Ryzhik,
Flame capturing with an advection--reaction--diffusion model.
{ Comb. Theory Model.}, {\bf 10}, 727--74, 2006. 


\bibitem{bacteria}
M. Zaki, N. Andrew, and R.H. Insall,
Entamoeba histolytica cell movement: a central role for selfgenerated chemokines and chemorepellents, Proc. Nat. Acad. Sci. {\bf 103}, 18751--18756, 2006. 

\end{thebibliography}

\end{document}